\documentclass[11pt]{amsart}
\usepackage[latin1]{inputenc}
\usepackage{amsmath}
\usepackage{amssymb}
\usepackage[all]{xy}
\usepackage{mathrsfs}
\usepackage{amsthm}

\theoremstyle{plain}
\newtheorem{thm}[equation]{Theorem}
\newtheorem{lem}[equation]{Lemma}
\newtheorem{prop}[equation]{Proposition}
\newtheorem{cor}[equation]{Corollary}

\theoremstyle{definition}
\newtheorem{defin}[equation]{Definition}
\newtheorem{remark}[equation]{Remark}

\numberwithin{equation}{subsection}

\def\sheafEnd{\mathcal{E} \hspace{-1pt} \mathit{nd}}
\def\sheafHom{\mathcal{H} \hspace{-1pt} \mathit{om}}
\def\sheafExt{\mathcal{E} \hspace{-1pt} \mathit{xt}}

\newcommand{\bk}{\Bbbk}

\newcommand{\bbX}{\mathbb X}
\newcommand{\bbZ}{\mathbb Z}
\newcommand{\calA}{\mathcal{A}}
\newcommand{\calB}{\mathcal{B}}
\newcommand{\calC}{\mathcal{C}}
\newcommand{\calD}{\mathcal{D}}

\newcommand{\calF}{\mathcal{F}}
\newcommand{\calG}{\mathcal{G}}
\newcommand{\calH}{\mathcal{H}}
\newcommand{\calI}{\mathcal{I}}
\newcommand{\calJ}{\mathcal{J}}
\newcommand{\calK}{\mathcal{K}}
\newcommand{\calL}{\mathcal{L}}
\newcommand{\calM}{\mathcal{M}}
\newcommand{\calN}{\mathcal{N}}
\newcommand{\calO}{\mathcal{O}}
\newcommand{\calP}{\mathcal{P}}

\newcommand{\calS}{\mathcal{S}}

\newcommand{\calU}{\mathcal{U}}

\newcommand{\wcalN}{\widetilde{\mathcal{N}}}
\newcommand{\wfrakg}{\widetilde{\mathfrak{g}}}
\newcommand{\wcalD}{\widetilde{\mathcal{D}}}

\newcommand{\frakn}{\mathfrak{n}}
\newcommand{\frakg}{\mathfrak{g}}
\newcommand{\frakt}{\mathfrak{t}}
\newcommand{\frakb}{\mathfrak{b}}

\newcommand{\frakL}{\mathfrak{L}}

\newcommand{\frakh}{\mathfrak{h}}
\newcommand{\frakZ}{\mathfrak{Z}}

\newcommand{\frakp}{\mathfrak{p}}

\newcommand{\frakX}{\mathfrak{X}}
\newcommand{\frakM}{\mathfrak{M}}
\newcommand{\fraksl}{\mathfrak{sl}}
\newcommand{\mbfA}{\mathbf{A}}
\newcommand{\mbfB}{\mathbf{B}}
\newcommand{\mbfC}{\mathbf{C}}

\newcommand{\mbfF}{\mathbf{F}}
\newcommand{\mbfG}{\mathbf{G}}
\newcommand{\lotimes}{{\stackrel{_L}{\otimes}}}
\newcommand{\Rtimes}{{\stackrel{_R}{\times}}}

\newcommand{\Gm}{{\mathbb{G}}_{\mathbf{m}}}

\newcommand{\Hom}{{\rm Hom}}
\newcommand{\Ext}{{\rm Ext}}
\newcommand{\Tor}{{\rm Tor}}

\newcommand{\Ind}{{\rm Ind}}

\newcommand{\For}{{\rm For}}
\newcommand{\Id}{{\rm Id}}

\newcommand{\Coh}{{\rm Coh}}
\newcommand{\QCoh}{{\rm QCoh}}

\newcommand{\Mod}{{\rm Mod}}

\newcommand{\qis}{{\rm qis}}

\newcommand{\pro}{{\rm prop}}
\newcommand{\rk}{{\rm rk}}
\newcommand{\qc}{{\rm qc}}
\newcommand{\fg}{{\rm fg}}

\newcommand{\SL}{{\rm SL}}

\newcommand{\rmi}{\rm{(i)}}
\newcommand{\rmii}{\rm{(ii)}}
\newcommand{\rmiii}{\rm{(iii)}}
\newcommand{\rmiv}{\rm{(iv)}}

\newcommand{\rmvi}{\rm{(vi)}}
\newcommand{\aff}{{\rm aff}}

\newcommand{\reg}{\mathrm{reg}}
\newcommand{\rs}{\mathrm{rs}}
\newcommand{\bX}{{\bf X}}
\newcommand{\bY}{{\bf Y}}
\newcommand{\bZ}{{\bf Z}}

\newcommand{\wS}{\widetilde{S}}
\newcommand{\wf}{\widetilde{f}}
\newcommand{\rmS}{\mathrm{S}}
\newcommand{\scS}{\mathscr{S}}
\newcommand{\bB}{\mathbb{B}}
\newcommand{\Z}{\mathbb{Z}}
\newcommand{\F}{\mathbb{F}}
\newcommand{\C}{\mathbb{C}}

\author{Roman Bezrukavnikov}
\address{Department of Mathematics \\ Massachusetts Institute of Technology \\ Cambridge MA \\ 02139 \\ USA.}
\email{bezrukav@math.mit.edu}

\author{Simon Riche}
\address{Clermont Universit{\'e}, Universit{\'e} Blaise Pascal, Laboratoire de Math{\'e}ma\-tiques, BP 10448, F-63000 Clermont-Ferrand. \newline
\indent CNRS, UMR 6620, Laboratoire de Math{\'e}matiques, F-63177 Aubi{\`e}re.}
\email{simon.riche@math.univ-bpclermont.fr}

\title[Affine braid group actions]{Affine braid group actions on derived categories of Springer resolutions}

\begin{document}

\begin{abstract}
In this paper we construct and study an action of the affine braid group associated to a semi-simple algebraic group on derived categories of coherent sheaves on various varieties related to the Springer resolution of the nilpotent cone. In particular, we describe explicitly the action of the Artin braid group. This action is a ``categorical version'' of Kazhdan--Lusztig--Ginzburg's construction of the affine Hecke algebra, and is used in particular in \cite{BM} in the course of the proof of Lusztig's conjectures on equivariant $K$-theory of Springer fibers.
\end{abstract}

\maketitle


\section*{Introduction}

\subsection{}

The goal of this paper is to introduce an action of the affine braid group
on the derived category of coherent sheaves on the Springer resolution (and
some related varieties) and prove some of its properties.

The most direct way to motivate this construction is via the well-known heuristics of
Springer correspondence theory. Let $\frakg$ be a semi-simple Lie algebra over
$\C$,\footnote{In the body of the paper we work over a finite localization
of $\Z$ or over a field of arbitrary characteristic rather than over
$\C$. Such details are ignored in the introduction.}
 let $\pi:\wfrakg \to \frakg$ be the Grothendieck--Springer map and $\pi':
\wcalN \to \calN \subset \frakg$ be the Springer map; here $\calN$ is the nilpotent cone
and $\wcalN$
is the cotangent bundle to the flag variety. Let $\frakg_{\reg}\subset \frakg$ be
the subset of regular elements and $\wfrakg_{\reg}=\pi^{-1}(\frakg_{\reg})$. Then
$\pi|_{\wfrakg_{\reg}}$ is a ramified Galois covering with Galois group $W$, the
Weyl group of
$\frakg$. Thus $W$ acts on $\wfrakg_{\reg}$ by deck transformations. Although the
action does not extend to an action on $\wfrakg$, it still induces various
interesting structures on the Springer resolution $\wcalN$. The most well-known example is the Springer action of $W$ on (co)homology of a fiber of
$\pi'$, called a Springer fiber. The procedure of passing from the action
of $W$ on $\wfrakg_{\reg}$ to the Springer action can be performed using
different tools, such as minimal (Goresky--MacPherson) extension of a
perverse sheaf (see \cite{LUSGr} for the original idea of this construction, and \cite{JANNil} for a detailled treatment and further references), nearby cycles, or degeneration of the correspondence cycle (see \cite{G1, CG}).

The main result of this paper can also be viewed as a realization of that
general idea. Namely, we show that a ``degeneration" of the action of $W$ on
$\wfrakg_{\reg}$ provides an action of the corresponding Artin braid group $\bB$
on the
derived categories of coherent sheaves $\calD^b \Coh(\wcalN)$, $\calD^b \Coh(\wfrakg)$.
More precisely, we consider the closure $Z_w$ of the graph of the action of $w\in W$
on $\wfrakg_{\reg}$. Using
this as a correspondence we get a functor $\calD^b \Coh(\wfrakg) \to
\calD^b \Coh(\wfrakg)$; we then prove that there exists an action of $\bB$ on
$\calD^b \Coh(\wfrakg)$ where
 a minimal length representative $T_w \in \bB$ of $w\in W$ acts by the
resulting functor. It also induces a compatible action on $\calD^b \Coh(\wcalN)$. 
 
 The fact that functors admitting such a simple description give an action
of the braid group is perhaps surprising; it implies that the closures
 of the graphs are Cohen--Macaulay.
 
 Furthermore, the categories in question carry an obvious action of the
weight lattice
 $\bbX$ of $G$ which is identified with the Picard group of the flag variety; here
an element $\lambda \in \bbX$ acts by twist by the corresponding line
bundle. We prove that
 this action of $\bbX$ together with the above action of $\bB$ generate
an action of the {\em extended affine braid group}\footnote{In the standard
terminology (see \cite{BLie}) this  is the extended affine braid group of the
\emph{Langlands dual} group $^L \hspace{-1pt} G$.
 In fact, this may be viewed as the simplest manifestation of the relation
of our $\bB_{\aff}$-action to Langlands duality mentioned below.} $\bB_{\aff}$.

In fact we construct a structure stronger than just an action of $\bB_{\aff}$
on the two derived categories of coherent sheaves; namely, we show the
existence of a (weak) {\em geometric action} of this group. Informally,
this means that the action of elements of the group come from ``integral
kernels," i.e.~complexes of sheaves on the square of the space, and
relations of the group come from isomorphisms between convolutions of the
kernels. Formal definition of this convolution requires basic formalism of
differential graded schemes. On the other hand, this geometric action induces a usual action
on the derived categories
of varieties obtained from $\wcalN$, $\wfrakg$ by base change. In the simplest
case of base
change to the transversal slice to a subregular nilpotent orbit we recover
the action of $\bB_{\aff}$  on the derived category of the minimal resolution
of a Kleinian singularity considered e.g.~in \cite{ST, Bd}.

\subsection{}

We now list some contexts where the action of $\bB_{\aff}$ appears and plays
an essential role (see \cite{BEZICM} for a complementary discussion).

The work \cite{BM} uses the geometric theory of representations of
semi-simple Lie algebras in positive characteristic developed in
\cite{BMR, BMR2} to deduce
Lusztig's conjectures on numerical properties of such representations (see \cite{LUSBas, LUSBas2}). It
uses the action considered in this paper, which is related to the action of $\bB_{\aff}$
on the derived category of modular representations by {\em intertwining} functors
(or shuffling functors, or Radon transforms, in a different terminology). In
fact, \cite{BM} and the present paper are logically interdependent. This application was our main motivation for considering this action over a localization of $\Z$. The action studied in this paper (and in particular its version for certain differential graded schemes) also plays a technical role in the study of Koszul duality for representations of
semi-simple Lie algebras in positive characteristic (see \cite{R2}).

The induced action of $\bB_{\aff}$ on the Grothendieck group of
$\C^*$-equivariant coherent sheaves factors
through an action of the {\em affine Hecke algebra}, i.e.~the action of
simple reflections satisfies a certain quadratic relation. A weak form of
the categorical counterpart of  the quadratic relation is used in
\cite{BM}; a more comprehensive
development of the idea that our action induces an action of the
``categorical affine
Hecke algebra" is the subject of a forthcoming work by Christopher Dodd.

In fact, in view of the work of Lusztig and Ginzburg, the monoidal
category
$\calD^b \Coh^{G\times \C^*}(\wfrakg \times_\frakg \wfrakg)$ (or the $\wcalN$-version; the
monoidal structure on these categories is defined below) can be considered
a categorification of the affine Hecke algebra.
See \cite{BEZICM} for a discussion of an equivalence between this
categorication and another one coming from perverse sheaves on the affine
flag variety of the dual group. 
Such an equivalence, inspired by the ideas of local geometric Langlands
duality theory, also implies the existence of the $\bB_{\aff}$-action
constructed in this paper (at least over $\C$).

Another approach to the construction of the $\bB$-action (over $\C$) relates it
to the well-known action of $\bB$ on the category of $D$-modules on the flag variety by
Radon transforms (see e.g.~\cite{BBM}). Passing
from $D$-modules to coherent sheaves on the cotangent bundle is achieved
by means of the Hodge $D$-modules formalism. We plan to develop this
approach in a future publication.

Finally, we would like to mention that in the $\wcalN$-version of the
construction, (inverses of) simple
generators act by reflection at a spherical functor in the sense of \cite{An, Ro} (see Remark \ref{rk:spherical}), and
that
in the particular case of groups of type $\mbfA$ the action (in its
non-geometric form, and over $\C$ rather than a localization of $\Z$)
has been constructed in \cite{KhT} and more recently, as a part of a
more general picture, in \cite{CaK}.

\subsection{Contents of the paper}

In Section \ref{sec:braidrelations} we prove (under very mild restrictions on the characteristic $p$ of our base field) that there exists an action of $\bB_{\aff}$ on $\calD^b \Coh(\wfrakg)$ and $\calD^b \Coh(\wcalN)$ where generators associated to \emph{simple reflections} in $W$ and elements of $\bbX$ act as stated above. This result was already proved under stronger assumptions and by less satisfactory methods in \cite{RAct}. We also extend this result to the schemes over a finite localization of $\Z$.

In Section \ref{sec:kernels} we prove that, if $p$ is bigger than the Coxeter number of $G$, the action of the element $T_w \in \bB$ ($w \in W$) is the convolution with kernel $\calO_{Z_w}$. This proof is based on representation theory of semi-simple Lie algebras in positive characteristic. We also extend this result to the schemes over a finite localization of $\Z$, and (as an immediate consequence) over an algebraically closed field of characteristic zero.

In Section \ref{sec:dgschemes} we prove generalities on dg-schemes, extending results of \cite[\S 1]{R2}. (Here, we concentrate on \emph{quasi-coherent} sheaves.) In particular, we prove a projection formula and a (non-flat) base change theorem in this context.

In Section \ref{sec:convolution} we use the results of Section \ref{sec:dgschemes} to show that the action of $\bB_{\aff}$ induces actions on categories of coherent sheaves on various (dg-)varieties related to $\wcalN$ and $\wfrakg$, in particular inverse images of Slodowy slices under the Springer resolution.

Finally, in Section \ref{sec:equivariant-version} we prove some equivariant analogues of the results of Section \ref{sec:convolution} which are needed in \cite{BM}.

\subsection{Acknowledgements}

Part of this work was done while both authors were participating the special year of the Institute for Advanced Study on ``New Connections of Representation Theory to Algebraic Geometry and Physics.'' Another part was done while the second author was visiting the Massachusetts Institute of Technology.

We thank Michel Brion and Patrick Polo for helpful comments.

R.~B.~is supported by NSF grant 
DMS-0854764. S.~R.~is supported by the ANR project RepRed (ANR-09-JCJC-0102-01).

\section{Existence of the action} \label{sec:braidrelations}

\subsection{Notation}
\label{ss:notation}

Let $G_{\Z}$ be a split connected, simply-connected, semi-simple algebraic group over $\mathbb{Z}$. Let $T_{\Z} \subset B_{\Z} \subset G_{\Z}$ be a maximal torus and a Borel subgroup in $G_{\Z}$. Let $\frakt_{\Z} \subset \frakb_{\Z} \subset \frakg_{\Z}$ be their respective Lie algebras. Let also $U_{\Z}$ be the unipotent radical of $B_{\Z}$, and $\frakn_{\Z}$ be its Lie algebra. Let $\Phi$ be the root system of $(G_{\Z},T_{\Z})$, and $\Phi^+$ be the positive roots, chosen as the roots of $\frakg_{\Z} / \frakb_{\Z}$. Let $\Sigma$ be the associated system of simple roots. Let also $\bbX:=X^*(T_{\Z})$ be the weight lattice. We let $\frakg^*_{\Z}$ be the coadjoint representation of $G_{\Z}$.

Let $W$ be the Weyl group of $\Phi$, and let $\scS=\{s_{\alpha}, \, \alpha \in \Sigma\}$ be the set of Coxeter generators associated to $\Sigma$ (called \emph{simple reflections}). Let $W_{\aff}^{\mathrm{Cox}}:= W \ltimes \mathbb{Z}R$ be the affine Weyl group, and $W_{\aff}:= W \ltimes \bbX$ be the extended affine Weyl group. Let $\bB \subset \bB_{\aff}^{\mathrm{Cox}} \subset \bB_{\aff}$ be the braid groups associated with $W \subset W_{\aff}^{\mathrm{Cox}} \subset W_{\aff}$ (see e.g.~\cite[\S 2.1.1]{BMR2} or \cite[\S 1.1]{RAct}). Note that $W$ and $W_{\aff}^{\mathrm{Cox}}$ are Coxeter groups, but not $W_{\aff}$. For $s,t \in \scS$, let us denote by $n_{s,t}$ the order of $st$ in $W$. Recall (see \cite{BRApp}) that $\bB_{\aff}$ has a presentation with generators $\{T_s, \, s \in \scS\}$ and $\{\theta_x, \, x \in \bbX\}$ and the following relations:
\begin{align*}
\rmi \ & T_s T_t \cdots = T_t T_s \cdots \quad \text{(} n_{s,t} \ \text{elements on each side)}; \\
\rmii \ & \theta_x \theta_y = \theta_{x+y}; \\
\rmiii \ & T_s \theta_x = \theta_x T_s \quad \text{if} \ s(x)=x; \\
\rmiv \ & \theta_x = T_s \theta_{x-\alpha} T_s \quad \text{if} \ s=s_{\alpha} \ \text{and} \ s(x)=x-\alpha. 
\end{align*}
Relations of type $\rmi$ are called \emph{finite braid relations}.

Recall that there exists a natural section $W_{\aff} \hookrightarrow \bB_{\aff}$ of the projection $\bB_{\aff} \twoheadrightarrow W_{\aff}$, denoted $w \mapsto T_w$ (see \cite[\S 2.1.2]{BMR2}, \cite{BRApp}). We denote by $\scS_{\aff}$ the set of simple affine reflections, i.e.~the natural set of Coxeter generators of $W_{\aff}^{\mathrm{Cox}}$ (which contains $\scS$), and by $\Sigma_{\aff}$ the set of affine simple roots. By definition, $\Sigma_{\aff}$ is in bijection with $\scS_{\aff}$ via a map denoted $\alpha \mapsto s_{\alpha}$. We denote by $\ell$ the length function on $W_{\aff}^{\mathrm{Cox}}$, and extend it naturally to $W_{\aff}$. Note that $\bB$ is generated by the elements $T_s$ ($s \in \scS$), $\bB_{\aff}^{\mathrm{Cox}}$ is generated by the elements $T_s$ ($s \in \scS_{\aff}$), and $\bB_{\aff}$ is generated by the elements $T_s$ ($s \in \scS_{\aff}$) and $T_{\omega}$ ($\omega \in \Omega:=\{w \in W_{\aff} \mid \ell(\omega)=0 \}$).

In this paper we study an action\footnote{As in \cite{BMR, BEZICM, RAct}, here we consider the \emph{weak} notion: an action of a group $\Gamma$ on a category $\calC$ is a group morphism from $\Gamma$ to the group of \emph{isomorphism classes} of auto-equivalences of $\calC$. See Remark \ref{rk:thm-description-kernels}(4) for comments on a stronger structure.} of the group $\bB_{\aff}$ on certain derived categories of coherent sheaves. Let us introduce the varieties we will consider.

Let $R$ be any (commutative) algebra. We replace the index $\Z$ by $R$ in all the notations introduced above to denote the base change to $R$. Let $\calB_{R}:=G_{R}/B_{R}$ be the flag variety. Let $\wcalN_{R}:=T^*\calB_{R}$ be its cotangent bundle. We have the more concrete description
\[
\wcalN_{R} = G_{R} \times^{B_{R}} (\frakg_{R}/\frakb_{R})^* = \{(X,gB_R) \in \frakg_{R}^* \times_R \calB_{R} \mid X_{|g \cdot \frakb_R}=0\}.
\]
Let also $\wfrakg_{R}$ be the Grothendieck resolution, defined as
\[
\wfrakg_{R} = G_{R} \times^{B_R} (\frakg_R/\frakn_R)^*=\{(X,gB_R) \in \frakg_R^* \times_R \calB_R \mid X_{|g \cdot \frakn_R}=0\}. \]
There is a natural inclusion $i : \wcalN_R \hookrightarrow \wfrakg_R$, and a natural morphism $\wfrakg_R \to \frakg^*_R$ induced by the projection on the first summand.

The varieties $\wfrakg_R$ and $\wcalN_R$ are endowed with an action of $G_R \times_R (\Gm)_R$, where $G_R$ acts by the natural (diagonal) action, and $(\Gm)_R$ acts via
\[
t \cdot (X,gB_R) \ = \ (t^2 X, gB_R).
\]
We will also consider the diagonal action of $G_R \times_R (\Gm)_R$ on $\wfrakg_R \times_R \wfrakg_R$ and $\wcalN_R \times_R \wcalN_R$.

For any $R$-scheme $X$ endowed with an action of $(\Gm)_R$, we denote by
\[
\langle 1 \rangle : \calD^b \Coh^{(\Gm)_R}(X) \to \calD^b \Coh^{(\Gm)_R}(X)
\]
the shift functor, i.e.~the tensor product with the free rank one $(\Gm)_R$-module given by the natural identification $(\Gm)_R \cong \mathrm{GL}(1,R)$. We denote by $\langle j \rangle$ the $j$-th power of $\langle 1 \rangle$.

If $\lambda \in \bbX$, we denote by $\calO_{\calB_R}(\lambda)$ the line bundle on $\calB_R$ naturally associated to $\lambda$. And if $X \to \calB_R$ is a variety over $\calB_R$, we denote by $\calO_X(\lambda)$ the inverse image of $\calO_{\calB_R}(\lambda)$.

For any $R$-scheme $X$, we denote by $\Delta X \subset X \times_R X$ the diagonal copy of $X$.

\subsection{Convolution}
\label{ss:convolution1}

Let $X$ be a smooth scheme over $R$, and let $p_1,p_2 : X \times_R X \to X$ be the projections on the first and on the second component. We define the full subcategory
\[
\calD^b_{\mathrm{prop}} \Coh(X \times_R X) \subset \calD^b \Coh(X \times_R X)
\]
as follows: an object of $\calD^b \Coh(X \times_R X)$ belongs to $\calD^b_{\mathrm{prop}} \Coh(X \times_R X)$ if it is isomorphic to the direct image of an object of $\calD^b\Coh(Y)$ for a closed subscheme $Y \subset X \times_R X$ such that the restrictions of $p_1$ and $p_2$ to $Y$ are proper. Any $\calF \in \calD^b_{\mathrm{prop}} \Coh(X \times_R X)$ gives rise to a functor
\[
F^{\calF}_X : \left\{
\begin{array}{ccc}
\calD^b \Coh(X) & \to & \calD^b \Coh(X) \\
\calM & \mapsto & R(p_2)_*(\calF \, \lotimes_{X \times_R X} \, (p_1)^* \calM)
\end{array}
\right. .
\]
Let $p_{1,2}, p_{2,3}, p_{1,3} : X \times_R X \times_R X \to X \times_R X$ be the natural projections. The category $\calD^b_{\mathrm{prop}} \Coh(X \times_R X)$ is endowed with a convolution product, defined by
\[
\calF \, \star \, \calG \ := \ R(p_{1,3})_* \bigl( (p_{1,2})^* \calG \, \lotimes_{X \times_R X \times_R X} \, (p_{2,3})^* \calF \bigr).
\]
With these definitions, for $\calF,\calG \in \calD^b_{\mathrm{prop}} \Coh(X \times_R X)$ we have a natural isomorphism (see e.g.~\cite[Lemma 1.2.1]{RAct} or \cite{Huy})
\[
F^{\calF}_X \circ F^{\calG}_X \ \cong \ F_X^{\calF \star \calG}.
\]

One can define similarly convolution functors for equivariant coherent sheaves. We will use the same notation in this setting also.

The convolution formalism is compatible with base change: if $R'$ is an $R$-algebra, and if $X'=X \times_{\mathrm{Spec}(R)} \mathrm{Spec}(R')$, then the (derived) pull-back under the morphism $X' \times_{R'} X' \to X \times_R X$ is monoidal, and the (derived) pull-back under the morphism $X' \to X$ is compatible with the actions.

\subsection{Statement}
\label{ss:statement-R}

Let now $R$ be a finite localization of $\Z$. For $s \in \scS$, we denote by
\[
Z_{s,R}
\]
the closure of the inverse image of the $G_R$-orbit of $(B_R/B_R,s B_R/B_R) \in \calB_R \times_R \calB_R$ (for the diagonal action) under the morphism $\wfrakg_R \times_{\frakg_R^*} \wfrakg_R \hookrightarrow \wfrakg_R \times_R \wfrakg_R \twoheadrightarrow \calB_R \times_R \calB_R$. It is a reduced closed subscheme of $\wfrakg_R \times_{\frakg_R^*} \wfrakg_R$. We also define the following closed subscheme of $(\wcalN_R \times_R \wfrakg_R)$:
\[
Z_{s,R}':= Z_{s,R} \cap (\wcalN_R \times_R \wfrakg_R).
\]
(Note that here we consider the \emph{scheme-theoretic} intersection.) It is easy to prove (see e.g. Lemma \ref{lem:Z_w'} below for a more general claim in the case of a field) that $Z_{s,R}'$ is in fact a closed subscheme of $\wcalN_R \times_{\frakg_R^*} \wcalN_R$.

We define $n_G \in \{1,2,3,6\}$ to be the smallest integer such that $2 \, | \, n_G$ if $G_{\Z}$ has a component of type $\mbfF_4$, and $3 \, | \, n_G$ if $G_{\Z}$ has a component of type $\mbfG_2$. The main result of this section is the following.

\begin{thm}
\label{thm:existenceaction-R}

Let $R=\Z[\frac{1}{n_G}]$.

There exists an action $b \mapsto \mathbf{J}_b$, respectively $b \mapsto \mathbf{J}'_b$, of $\bB_{\aff}$ on the category $\calD^b \Coh(\wfrakg_R)$, respectively $\calD^b \Coh(\wcalN_R)$, such that
\begin{enumerate}
\item
For $s \in \scS$, $\mathbf{J}_{T_s}$, respectively $\mathbf{J}'_{T_s}$, is isomorphic to $F^{\calO_{Z_{s,R}}}_{\wfrakg_R}$, respectively $F^{\calO_{Z_{s,R}'}}_{\wcalN_R}$;
\item
For $x \in \bbX$, $\mathbf{J}_{\theta_x}$, respectively $\mathbf{J}'_{\theta_x}$, is isomorphic to $F^{\calO_{\Delta \wfrakg_R}(x)}_{\wfrakg_R}$, respectively $F^{\calO_{\Delta \wcalN_R}(x)}_{\wcalN_R}$, or equivalently to the functor of tensoring with the line bundle $\calO_{\wfrakg_R}(x)$, respectively $\calO_{\wcalN_R}(x)$.
\end{enumerate}

Moreover, these two actions are compatible (in the obvious sense) with the inverse and direct image functors 
\[
\xymatrix{ \calD^b \Coh(\wcalN_R) \ar@<0.5ex>[rr]^-{Ri_*} & & \calD^b \Coh(\wfrakg_R). \ar@<0.5ex>[ll]^-{Li^*} }
\]

\end{thm}

In fact, we will also prove an equivariant analogue of this result:

\begin{thm}
\label{thm:existenceaction-R-equivariant}

Let $R=\Z[\frac{1}{n_G}]$.

There exists an action $b \mapsto \mathbf{J}_b^{\mathrm{eq}}$, respectively $b \mapsto \mathbf{J}^{\prime,\mathrm{eq}}_b$, of $\bB_{\aff}$ on the category $\calD^b \Coh^{G_R \times_R (\Gm)_R}(\wfrakg_R)$, respectively $\calD^b \Coh^{G_R \times_R (\Gm)_R}(\wcalN_R)$, such that
\begin{enumerate}
\item
For $s \in \scS$, $\mathbf{J}^{\mathrm{eq}}_{T_s}$, respectively $\mathbf{J}^{\prime,\mathrm{eq}}_{T_s}$, is isomorphic to $F^{\calO_{Z_{s,R}}\langle 1 \rangle}_{\wfrakg_R}$, respectively $F^{\calO_{Z_{s,R}'}\langle 1 \rangle}_{\wcalN_R}$;
\item
For $x \in \bbX$, $\mathbf{J}^{\mathrm{eq}}_{\theta_x}$, respectively $\mathbf{J}^{\prime,\mathrm{eq}}_{\theta_x}$, is isomorphic to $F^{\calO_{\Delta \wfrakg_R}(x)}_{\wfrakg_R}$, respectively $F^{\calO_{\Delta \wcalN_R}(x)}_{\wcalN_R}$, or equivalently to the functor of tensoring with the line bundle $\calO_{\wfrakg_R}(x)$, respectively $\calO_{\wcalN_R}(x)$.
\end{enumerate}

Moreover, these two actions are compatible (in the obvious sense) with the inverse and direct image functors
\[
\xymatrix{ \calD^b \Coh^{G_R \times_R (\Gm)_R}(\wcalN_R) \ar@<0.5ex>[rr]^-{Ri_*} & & \calD^b \Coh^{G_R \times_R (\Gm)_R}(\wfrakg_R). \ar@<0.5ex>[ll]^-{Li^*} }
\]

\end{thm}

\subsection{Preliminary results}

Let $R$ be a finite localization of $\Z$, and let $A$ be a finitely generated $R$-algebra, which is flat over $R$. For any prime $p \in \Z$ which is not invertible in $R$, consider the specialization $A_p:=A \otimes_R \F_p$, and the extension of scalars $A_{\overline{p}}:=A \otimes_R \overline{\F_p} \cong A_p \otimes_{\F_p} \overline{\F_p}$. For any object $M$ of $\calD^b \Mod(A)$ we set
\begin{align*}
M_p \, & := \, M \, \lotimes_R \, \F_p \quad \text{in } \ \calD^b \Mod(A_p), \\
M_{\overline{p}} \, & := \, M \, \lotimes_R \, \overline{\F_p} \quad \text{in } \ \calD^b \Mod(A_{\overline{p}}).
\end{align*}

\begin{lem}
\label{lem:specialization}

Let $M \in \calD^b \Mod^{\fg}(A)$.

\begin{enumerate}
\item 
\label{it:specialization-0} 
If $M_{\overline{p}}=0$ for any prime $p \in \Z$ not invertible in $R$, then $M=0$.
\item 
\label{it:specialization-degree-0}
If $M_{\overline{p}}$ is concentrated in degree $0$ for any prime $p$ not invertible in $R$, then $M$ is concentrated in degree $0$, and is flat over $R$.

\end{enumerate}

\end{lem}

\begin{proof}
First, one can clearly replace $M_{\overline{p}}$ by $M_p$ in these properties.

The ring $R$ has global dimension $1$, hence any object of $\calD^b \Mod(R)$ is isomorphic to the direct sum of its (shifted) cohomology objects. Hence it is enough to prove the following properties for a finitely generated $A$-module $M$:
\begin{itemize}

\item[$(*)$] 
\label{it:specialization-module}
if $M \otimes_R \F_p=0$ for any $p$, then $M=0$;
\item[$(**)$]
\label{it:specialization-module-2}
if $\Tor^R_{-1}(M,\F_p)=0$ for any $p$, then $M$ is flat over $R$.

\end{itemize}

Let us first prove $(*)$. It is enough to prove that for any maximal ideal $\frakM \subset A$, the localization $M_{\frakM}$ is zero. Then, by Nakayama's lemma, it is enough to prove that $M_{\frakM}/\frakM \cdot M_{\frakM}=M/\frakM \cdot M=0$. However, by the general form of the Nullstellensatz (\cite[Theorem 4.19]{Ei}), $A/\frakM$ is a field extension of $R/p\cdot R = \F_p$ for some prime $p$ not invertible in $R$. Hence $M/\frakM \cdot M$ is a quotient of $M \otimes_R \F_p$, hence is indeed zero.

Property $(**)$ implies that multiplication by $p$ is injective on $M$ for any $p$. Hence $M$ is torsion-free over $R$, hence flat.
\end{proof}

\begin{lem}
\label{lem:cohomology-Grothendieck-resolution}

Let $\bk$ be an algebraically closed field.

We have
\[
R^i\Gamma(\wfrakg_{\bk}, \calO_{\wfrakg_{\bk}})^{\Gm} \ = \left\{
\begin{array}{cl}
\bk & \text{if } \ i=0; \\
0 & \text{if } \ i \neq 0.
\end{array}
\right.
\]

\end{lem}

\begin{proof}
Let $r : \wfrakg_{\bk} \to \calB_{\bk}$ be the natural projection. Then we have $r_* \calO_{\wfrakg} \cong \rmS_{\calO_{\calB_{\bk}}}(\calM)$, where $\calM$ is the dual to the sheaf of sections of $\wfrakg_{\bk}$. The action of $\Gm$ on $\calB_{\bk}$ is trivial; the $\Gm$-equivariant structure on $\rmS_{\calO_{\calB_{\bk}}}(\calM)$ is given by the grading such that $\calM$ is in degree $-2$. Hence the claim follows from the similar result for $\calB_{\bk}$, which is well-known (see e.g.~\cite[Theorem 3.1.1]{BrK}).
\end{proof}

\begin{prop}
\label{prop:isom-Delta}

Let $R$ be any finite localization of $\Z$.

Let $\calM \in \calD^b \Coh^{G_R \times_R (\Gm)_R}(\wfrakg_R \times_R \wfrakg_R)$. Assume that for any prime $p$ not invertible in $R$ there is an isomorphism
\[
\calM \, \lotimes_R \, \overline{\F_p} \ \cong \ \calO_{\Delta \wfrakg_{\overline{\F_p}}}
\]
in $\calD^b \Coh^{G_{\overline{\F_p}} \times_{\overline{\F_p}} (\Gm)_{\overline{\F_p}}}(\wfrakg_{\overline{\F_p}} \times_{\overline{\F_p}} \wfrakg_{\overline{\F_p}})$. Then there exists an isomorphism
\[
\calM \ \cong \ \calO_{\Delta \wfrakg_R}
\]
in $\calD^b \Coh^{G_R \times_R (\Gm)_R}(\wfrakg_R \times_R \wfrakg_R)$.

\end{prop}

\begin{proof}
By Lemma \ref{lem:specialization}\eqref{it:specialization-degree-0}, $\calM$ is concentrated in degree $0$, i.e.~is an equivariant coherent sheaf, which is flat over $R$. Consider the object
\[
M \, := \ R\Gamma(\wfrakg_R \times_R \wfrakg_R, \, \calM)^{(\Gm)_R}
\]
in $\calD^b \Mod(R)$. Here, $(-)^{(\Gm)_R}$ is the functor of $(\Gm)_R$-fixed points. Note that this functor is exact by \cite[Lemma I.4.3.(b)]{JANAlg}, and commutes with specialization by \cite[Equation I.2.11.(10)]{JANAlg}. As $\wfrakg_R$ is proper over $\frakg^*_R$, the object $R\Gamma(\wfrakg_R \times_R \wfrakg_R, \, \calM)$ is an object of the category $\calD^b \Mod^{\fg}(\rmS_R(\frakg_R \oplus \frakg_R))$, hence $M$ is in fact in the subcategory $\calD^b \Mod^{\fg}(R)$. For any $p$ not invertible in $R$ we have
\begin{align*}
M_{\overline{p}} \ & \cong \ R\Gamma(\wfrakg_{\overline{\F_p}} \times_{\overline{\F_p}} \wfrakg_{\overline{\F_p}}, \, \calM_{\overline{\F_p}})^{(\Gm)_{\overline{\F_p}}} \\
& \cong \ R\Gamma(\wfrakg_{\overline{\F_p}}, \, \calO_{\wfrakg_{\overline{\F_p}}})^{(\Gm)_{\overline{\F_p}}},
\end{align*}
where the first isomorphism follows from base change. By Lemma \ref{lem:cohomology-Grothendieck-resolution}, this object is concentrated in degree $0$. Hence, by Lemma \ref{lem:specialization}\eqref{it:specialization-degree-0}, $M$ is concentrated in degree $0$, i.e.~is a finitely generated (graded) $R$-module.

For any prime $p$ not invertible in $R$ we have 
\[
M \otimes_R \overline{\F_p} \ \cong \ \overline{\F_p}.
\]
by Lemma \ref{lem:cohomology-Grothendieck-resolution}. Using the classical description of finitely generated modules over a principal ideal domain, it follows that $M \cong R$. Hence we also have
\[
\Gamma(\wfrakg_R \times_R \wfrakg_R, \, \calM)^{G_R \times_R (\Gm)_R} \, \cong \, R.
\]
Hence there exists a non-zero $(G_R \times_R (\Gm)_R)$-equivariant morphism
\[
\phi : \calO_{\wfrakg_R \times_R \wfrakg_R} \to \calM,
\]
which is unique up to an invertible scalar in $R$.

Let $\calJ$ be the cocone of $\phi$, so that we have a distinguished triangle
\[
\calJ \, \to \, \calO_{\wfrakg_R \times_R \wfrakg_R} \, \xrightarrow{\phi} \, \calM \, \xrightarrow{+1}.
\]
For any $p$ not invertible in $R$, $\calO_{\wfrakg_R \times_R \wfrakg_R} \, \lotimes_R \, \overline{\F_p}$ and $\calM \, \lotimes_R \, \overline{\F_p}$ are both concentrated in degree $0$, and $\phi \otimes_R \overline{\F_p}$ is surjective. It follows that $\calM \, \lotimes_R \, \overline{\F_p}$ is also concentrated in degree $0$. By Lemma \ref{lem:specialization}\eqref{it:specialization-degree-0}, we deduce that $\calJ$ is concentrated in degree $0$, i.e.~that $\phi$ is surjective. It follows that
\[
\calM \ \cong \ \calO_X
\]
for some $G_R \times_R (\Gm)_R$-stable closed subscheme $X \subset \wfrakg_R \times_R \wfrakg_R$.

For any $p$ not invertible in $R$ we have $X \times_{\mathrm{Spec}(R)} \mathrm{Spec}(\overline{\F_p}) \cong \Delta \wfrakg_{\overline{\F_p}}$. In particular, this fiber product is reduced. As moreover $X$ is flat over $R$, the arguments of \cite[Proof of Proposition 1.6.5]{BrK} imply that $X$ itself is reduced.

Let $Y$ be the restriction of $X$ to $(\wfrakg_R \times_R \wfrakg_R) \smallsetminus (\Delta \wfrakg_R)$. For any $p$ not invertible in $R$ we have $\calO_Y \, \lotimes_{R} \, \overline{\F_p}=0$. By Lemma \ref{lem:specialization}\eqref{it:specialization-0} we deduce that $Y$ is empty. As $X$ is reduced, it follows that $X$ is included in $\Delta \wfrakg_R$, i.e.~that $\phi$ factors through a $(G_R \times_R (\Gm)_R)$-equivariant morphism
\[
\psi : \calO_{\Delta \wfrakg_R} \to \calO_X.
\]
Let $\calI$ be the cone of $\psi$. Then for any $p$ not invertible in $R$ we have $\calI \, \lotimes_R \, \overline{\F_p}=0$. By Lemma \ref{lem:specialization}\eqref{it:specialization-0} again we deduce that $\calI=0$, i.e.~that $\psi$ is an isomorphism.
\end{proof}

\subsection{Reduction to the case of an algebraically closed field of positive characteristic}
\label{ss:reduction}

To prove Theorem \ref{thm:existenceaction-R}, it is easier to replace $R$ by an algebraically closed field of positive characteristic. In this subsection we explain how to justify this reduction.

Let $s \in \scS$. Then the subscheme $Z_{s,\Z} \subset \wfrakg_{\Z} \times_{\Z} \wfrakg_{\Z}$ can be described explicitly (see \cite[\S 1.4]{RAct}). Let $P_{s,\Z} \subset G_{\Z}$ be the minimal parabolic subgroup of $G_{\Z}$ containing $B_{\Z}$ associated with $s$, and let $\calP_{s,\Z}=G_{\Z}/P_{s,\Z}$ be the associated partial flag variety. Then $Z_{s,\Z}$ is a vector bundle over $\calB_{\Z} \times_{\calP_{s,\Z}} \calB_{\Z}$, of rank $\rk_{\Z}(\frakb_{\Z})-1$. (Note that in \cite{RAct} we work over an algebraically closed field; the case of $\Z$ is similar. Also, if $s=s_{\alpha}$ for $\alpha \in \Sigma$, the scheme $Z_s$ is denoted by $S_{\alpha}$ in \cite{RAct}.) In particular, $Z_{s,\Z}$ is a smooth scheme, flat over $\Z$. And, for any algebraically closed field $\bk$, the reduction $Z_{s,\bk}:=Z_s \times_{\mathrm{Spec}(R)} \mathrm{Spec}(\bk)$ is also smooth, and is also the closure of the inverse image in $\wfrakg_{\bk} \times_{\frakg^*_{\bk}} \wfrakg_{\bk}$ of the orbit of $(B_{\bk}/B_{\bk},sB_{\bk}/B_{\bk}) \in \calB_{\bk} \times_{\bk} \calB_{\bk}$.

\begin{lem}
\label{lem:inverse-R}

Let $R$ be as in Theorem {\rm \ref{thm:existenceaction-R}}.

There exist isomorphisms in $\calD^b \Coh^{G_R \times_R (\Gm)_R}(\wfrakg_R \times_R \wfrakg_R)$:
\begin{align*}
\calO_{Z_{s,R}}(-\rho, \rho-\alpha) \star \calO_{Z_{s,R}} \ & \cong \ \calO_{\Delta \wfrakg_{R}} \langle -2 \rangle, \\
\calO_{Z_{s,R}} \star \calO_{Z_{s,R}}(-\rho, \rho-\alpha) \ & \cong \ \calO_{\Delta \wfrakg_{R}} \langle -2 \rangle.
\end{align*}

\end{lem}

\begin{proof}
Let $\calM := \calO_{Z_{s,R}}(-\rho, \rho-\alpha) \star \calO_{Z_{s,R}} \langle 2 \rangle$. By Proposition \ref{prop:inverse} below, for any $p$ not invertible in $R$ we have a $(G_{\overline{\F_p}} \times_{\overline{\F_p}} (\Gm)_{\overline{\F_p}})$-equivariant isomorphism
\[
\calM \, \lotimes_R \, \overline{\F_p} \ \cong \ \calO_{\Delta \wfrakg_{\overline{\F_p}}}.
\]
By Proposition \ref{prop:isom-Delta}, we deduce that $\calM \cong \calO_{\Delta \wfrakg_R}$ in $\calD^b \Coh^{G_R \times_R (\Gm)_R}(\wfrakg_R \times_R \wfrakg_R)$. The second isomorphism can be proved similarly.
\end{proof}

\begin{proof}[Proof of Theorems {\rm \ref{thm:existenceaction-R}} and {\rm \ref{thm:existenceaction-R-equivariant}}]
Consider first the case of $\wfrakg_R$. What we have to prove is that the kernels $\calO_{Z_{s,R}}$ ($s \in \scS$) and $\calO_{\Delta \wfrakg_R}(x)$ ($x \in \bbX$) satisfy the relations of the affine braid group $\bB_{\aff}$ given in \S \ref{ss:notation} in the monoidal category $(\calD^b \Coh^{G_R \times_R (\Gm)_R}(\wfrakg_R \times_R \wfrakg_R),\star)$. Using Lemma \ref{lem:inverse-R}, one can rewrite these relations as stating that certain objects of $\calD^b \Coh^{G_R \times_R (\Gm)_R}(\wfrakg_R \times_R \wfrakg_R)$ are isomorphic to $\calO_{\Delta \wfrakg_R}$. By Proposition \ref{prop:isom-Delta}, it is enough to prove the isomorphisms after applying the derived specialization $(- \, \lotimes_R \, \overline{\F_p})$. In this setting, relation $\rmi$ is proved in Corollary \ref{cor:braidrelations}, and the other relations are proved in the end of \S \ref{ss:kernels-simple-reflections}.

The case of $\wcalN_R$ is similar. The compatibility of the actions with the functor $Ri_*$ easily follows from the definition and the projection formula. The compatibility with $Li^*$ follows by adjunction.
\end{proof}

\subsection{Statement for an algebraically closed field}
\label{ss:statement-k}

From now on and until the end of the section we fix an algebraically closed field $\bk$ of characteristic $p \geq 0$. All the schemes we will consider will be over $\bk$. In particular, we will consider the specialization of all the varieties defined over $\Z$ above. For simplicity we drop the index ``$\bk$.'' In particular we have the variety $Z_s$ defined in \S \ref{ss:reduction}. We set 
\[
Z_s' \ := \ Z_s \cap (\wcalN \times \wfrakg) \ \cong \ Z_{s,\Z}' \times_{\mathrm{Spec}(\Z)} \mathrm{Spec}(\bk).
\]
This is a closed subscheme of $\wcalN \times \wcalN$.

In the end of this section we will prove the following result, which is a version of Theorems \ref{thm:existenceaction-R} and \ref{thm:existenceaction-R-equivariant} over $\bk$. Note that the case $p=0$ is not excluded, though it is not needed to prove Theorems \ref{thm:existenceaction-R} and \ref{thm:existenceaction-R-equivariant}.

\begin{thm} \label{thm:existenceaction}

Assume $p \neq 2$ if $G$ has a component of type $\mbfF_4$, and $p \neq 3$ if $G$ has a component of type $\mbfG_2$. Then there exists an action $b \mapsto \mathbf{J}_b$, respectively $b \mapsto \mathbf{J}'_b$, of $\bB_{\aff}$ on the category $\calD^b \Coh(\wfrakg)$, respectively $\calD^b \Coh(\wcalN)$, such that
\begin{enumerate}
\item
For $s \in \scS$, $\mathbf{J}_{T_s}$, respectively $\mathbf{J}'_{T_s}$, is isomorphic to $F^{\calO_{Z_s}}_{\wfrakg}$, respectively $F^{\calO_{Z_s'}}_{\wcalN}$;
\item
For $x \in \bbX$, $\mathbf{J}_{\theta_x}$, respectively $\mathbf{J}'_{\theta_x}$, is isomorphic to $F^{\calO_{\Delta \wfrakg}(x)}_{\wfrakg}$, respectively $F^{\calO_{\Delta \wcalN}(x)}_{\wcalN}$, or equivalently to the functor of tensoring with the line bundle $\calO_{\wfrakg}(x)$, respectively $\calO_{\wcalN}(x)$.
\end{enumerate}

Moreover, these two actions are compatible (in the obvious sense) with the inverse and direct image functors \[ \xymatrix{ \calD^b \Coh(\wcalN) \ar@<0.5ex>[rr]^-{Ri_*} & & \calD^b \Coh(\wfrakg). \ar@<0.5ex>[ll]^-{Li^*} } \]

Similar results hold for the categories $\calD^b \Coh^{G \times \Gm}(\wfrakg)$, $\calD^b \Coh^{G \times \Gm}(\wcalN)$.

\end{thm}

This theorem was announced by the first author in \cite{BEZICM}. It was proved, in the case $G$ has no component of type $\mbfG_2$ (with a slightly stronger restriction on $p$), by the second author in \cite{RAct}. Here we do not make any assumption on the group, and we avoid case-by-case analysis.

The proof of the theorem is given in \S \ref{ss:kernels-simple-reflections} and Corollary \ref{cor:braidrelations}. It is based on a different interpretation of the scheme $Z_s$, see Remark \ref{rk:Z_s-graph}.

\begin{remark}
\label{rk:spherical}
Let $s \in \scS$, and consider the associated parabolic flag variety $\calP_s$ (see \S \ref{ss:reduction}). Set $\wcalN_s:=T^*(\calP_s)$, and let
$D_s:=\wcalN_s \times_{\calP_s} \calB$.
We have a closed embedding $i_s:D_s \hookrightarrow \wcalN$ (where $D_s$ is a divisor) and a projection
$p_s:D_s\to \wcalN_s$, which is a ${\mathbb P}^1$-bundle. It
follows from \cite[ Example 3]{An} that the functor
\[
R(i_s)_* L(p_s)^* : \calD^b \Coh(\wcalN_s) \to \calD^b \Coh(\wcalN)
\]
is a {\em spherical functor} of dimension 2. Hence this is also the case for the functor 
\[
\calM \mapsto \calO_{\wcalN}(-\rho) \otimes_{\calO_{\wcalN}} R(i_s)_* L(p_s)^* (\calM)
\]
(see \cite[Proposition 2.(2)]{An}). It is not hard to see, using the exact sequences of \cite[Lemma 6.1.1]{RAct}, that the
corresponding twist functor defined in \cite{An} (see also \cite{Ro}) is
isomorphic to the action of $(T_s)^{-1}$ constructed here. Since the element
$T_{s_{\alpha_0}}$ for an affine simple root $\alpha_0 \in \scS_{\aff} \smallsetminus \scS$ is conjugate to $T_s$ for some $s \in \scS$ (see \cite[Lemma 2.1.1]{BM}), we see
that all inverses of Coxeter generators of $\bB_{\aff}^{\mathrm{Cox}}$ act by reflection at a spherical
functor.
\end{remark}

\subsection{More notation}

For each positive root $\alpha$, there are subgroups $U_{\alpha}$, $U_{-\alpha}$ of $G$ naturally attached to $\alpha$ and $-\alpha$. We choose isomorphisms of algebraic groups
$u_{\alpha} : \bk \overset{\sim}{\to} U_{\alpha}$ and $u_{-\alpha}
: \bk \overset{\sim}{\to} U_{-\alpha}$ such that for all $t \in T$ we have $t
\cdot u_{\alpha}(x) \cdot t^{-1}=u_{\alpha}(\alpha(t) x)$ and $t
\cdot u_{-\alpha}(x) \cdot t^{-1}=u_{-\alpha}(\alpha(t)^{-1} x)$,
and such that these morphisms extend to a morphism of algebraic groups
$\psi_{\alpha} :
\SL(2,\bk) \to G$ such that
\[
\psi_{\alpha} \left(%
\begin{array}{cc}
  1 & x \\
  0 & 1 \\
\end{array}%
\right) = u_{\alpha}(x), \quad \psi_{\alpha} \left(%
\begin{array}{cc}
  1 & 0 \\
  x & 1 \\
\end{array}%
\right) = u_{-\alpha}(x), \quad x \in \bk,
\]
\[
\psi_{\alpha} \left(%
\begin{array}{cc}
  y & 0 \\
  0 & y^{-1} \\
\end{array}%
\right) = \alpha^{\vee}(y), \quad y \in \bk^{\times}.
\]
We define the elements
\[
e_{\alpha}:=d(u_{\alpha})_0(1), \quad
e_{-\alpha}:=d(u_{-\alpha})_0(1), \quad
h_{\alpha}:=[e_{\alpha},e_{-\alpha}]=d(\alpha^{\vee})_1(1).
\]

We let
\[
\frakn_+ \ := \ \bigoplus_{\alpha \in \Phi^+} \, \bk \cdot e_{\alpha}.
\]
Then $\frakn_+$ is the Lie algebra of the unipotent radical of the Borel subgroup of $G$ opposite to $B$ with respect to $T$.

If $P$ is a parabolic subgroup of $G$ and $V$ is any finite dimensional $P$-module, there exists a natural vector bundle $\calL_{G/P}(V)$ on $G/P$ associated to $V$ (see \cite[\S I.5.8]{JANAlg}). If $X \to G/P$ is a variety over $G/P$, we denote by $\calL_X(V)$ the inverse image of $\calL_{G/P}(V)$.

\subsection{Regular elements in $\frakg^*$}

In this subsection we prove some elementary facts about the coadjoint action of $G$ in $\frakg^*$. If $p$ is not too small, then $\frakg^*$ is isomorphic to $\frakg$ as a $G$-module, and all the results we need are well known (see e.g.~\cite{JANNil}). However, these facts have to be checked directly (without using $\frakg$) if we want to relax the assumptions on $p$. Some of these facts are proved in \cite[\S 5]{Xu}, using the same arguments as those of \cite{JANNil} and some results of \cite{KW}.

Let $\xi \in (\frakg/\frakn)^*$. The restriction morphism $(\frakg/\frakn)^* \to (\frakb/\frakn)^*$ is $B$-equiva\-riant, and the $B$-action on $(\frakb/\frakn)^*$ is trivial. It follows that the centralizer $Z_B(\xi)$ has dimension at least the rank $r$ of $G$. The same is true a fortiori for $Z_G(\xi)$. By \cite[Lemma 3.3]{KW} we have $\frakg^*=G \cdot (\frakg/\frakn)^*$. Hence for all $\xi \in \frakg^*$ we have $\dim(Z_G(\xi)) \geq r$. We denote by $\frakg^*_{\reg}$ the set of regular elements in $\frakg^*$, i.e.~the set of $\xi \in \frakg^*$ such that $\dim(Z_G(\xi))=\rk(G)$. By standard arguments, this set is open in $\frakg^*$ (see e.g.~\cite[\S 1.4]{HUMConj}).

\begin{lem} \label{lem:rs}

Let $\xi \in \frakg$, and assume that $\xi_{|\frakn \oplus \frakn^+}=0$. Then $\xi$ is regular if and only if for all $\alpha \in \Phi$, $\xi(h_{\alpha})\neq 0$.

\end{lem}

\begin{proof} This follows immediately from \cite[Lemma 3.1$\rmi$, $\rmii$, $\rmiv$]{KW}. \end{proof}

Let $\frakg^*_{\rs} \subset \frakg^*$ be the set of regular semi-simple elements, i.e.~of $\xi \in \frakg^*$ such that there exists $g \in G$ such that $g \cdot \xi$ satisfies the conditions of Lemma \ref{lem:rs}. By \cite[Theorem 4$\rmvi$]{KW}, under our hypotheses this set is non-empty and open in $\frakg^*$. By Lemma \ref{lem:rs} it is contained in $\frakg_{\reg}^*$.

\begin{remark} If $G=\mathrm{SL}(2,\bk)$ and $p=2$, then the set of regular semi-simple elements of $\frakg$ is empty, whereas $\frakg^*_{\rs} \neq \emptyset$.\end{remark}

Let us denote by $(\frakg/\frakn)^*_{\rs}$, respectively $(\frakg/(\frakn \oplus \frakn^+))^*_{\rs}$ the intersection $\frakg^*_{\rs} \cap (\frakg/\frakn)^*$, respectively $\frakg^*_{\rs} \cap (\frakg/(\frakn \oplus \frakn^+))^*$.

\begin{lem} \label{lem:rsb}

The action morphism induces an isomorphism of varieties \[ U \times (\frakg/(\frakn \oplus \frakn^+))^*_{\rs} \xrightarrow{\sim} (\frakg/\frakn)^*_{\rs}. \]

\end{lem}

\begin{proof} Let us first show that this morphism is surjective. Let $\xi \in (\frakg/\frakn)^*_{\rs}$. We have $\dim(Z_B(\xi)) \geq r$, hence $Z_G(\xi)^{\circ}=Z_B(\xi)^{\circ}$, hence in particular $Z_G(\xi)^{\circ} \subset B$. As $\xi$ is regular semi-simple, $Z_G(\xi)^{\circ}$ is a maximal torus. Hence there exists $b \in B$ such that $Z_G(b \cdot \xi)^{\circ}=b \cdot Z_G(\xi)^{\circ} \cdot b^{-1}=T$. Then $T$ stabilizes $b \cdot \xi$, hence $b \cdot \xi \in (\frakg/(\frakn \oplus \frakn^+))^*$. Writing $b=tu$ for some $t \in T$, $u \in U$, we have $u \cdot \xi \in (\frakg/(\frakn \oplus \frakn^+))^*_{\rs}$.

Now we prove that the morphism is injective. For $u \in U$ and $\xi \in (\frakg/(\frakn \oplus \frakn^+))^*$ we have $u \cdot \xi - \xi \in (\frakg / \frakb)^*$. Hence, using the decomposition $(\frakg/\frakn)^*=(\frakg/(\frakn \oplus \frakn^+))^* \oplus (\frakg / \frakb)^*$, it follows that $\xi$ is uniquely determined by $u \cdot \xi$. Then if $\xi$ is regular semi-simple and $u_1 \cdot \xi = u_2 \cdot \xi$ we have $(u_1)^{-1} u_2 \in Z_G(\xi) \subset N_G(Z_G(\xi)^{\circ})$. By Lemma \ref{lem:rs} we have $Z_G(\xi)^{\circ}=T$, hence $(u_1)^{-1} u_2 \in N_G(T) \cap U$, which implies $u_1=u_2$.

Then one can prove that the inverse bijection $(\frakg/\frakn)^*_{\rs} \to U \times (\frakg/(\frakn \oplus \frakn^+))^*_{\rs}$ is a morphism exactly as in \cite[\S 13.3]{JANNil}.\end{proof}

Let $\wfrakg_{\rs}$ be the inverse image of $\frakg^*_{\rs}$ under the natural morphism $\pi : \wfrakg \to \frakg^*$. Using \cite[Lemma 5.4, Lemma 5.5$\rmiii$]{Xu}, one obtains the following (see also \cite[\S 13.4]{JANNil}).

\begin{cor} \label{cor:actionW}

There exists a free action of $W$ on the variety $\wfrakg_{\rs}$, such that the restriction $\pi_{\rs}: \wfrakg_{\rs} \to \frakg^*_{\rs}$ of $\pi$ is a principal $W$-bundle.\qed

\end{cor}

\begin{remark}
\label{rk:Z_s-graph}
Let $s \in \scS$. Using the action of Corollary \ref{cor:actionW}, one can give a different interpretation of the variety $Z_s$: it is the closure of the graph of the action of $s$ on $\wfrakg_{\rs}$. Indeed, $Z_s$ clearly contains this closure, and both schemes are reduced, irreducible and of the same dimension.
\end{remark}

\subsection{Action of $W$ on $\wfrakg_{\reg}$}

Consider the morphism $\pi : \wfrakg \to \frakg^*$. Let now $\wfrakg_{\reg}$ be the inverse image of $\frakg^*_{\reg}$ under $\pi$. In this subsection we prove that the action of $W$ on $\wfrakg_{\rs}$ (see Corollary \ref{cor:actionW}) extends to $\wfrakg_{\reg}$.

\begin{lem} \label{lem:regsl2}

Let $\xi \in (\frakg/\frakn)^*_{\reg}$, and $\alpha \in \Sigma$. Then either $\xi(h_{\alpha})$ or $\xi(e_{\alpha})$ is non-zero.

\end{lem}

\begin{proof} Let $P_{\alpha}$ be the minimal standard parabolic subgroup of $G$ attached to $\alpha$, and let $\frakp_{\alpha}$ be its Lie algebra. Assume $\xi(h_{\alpha})=\xi(e_{\alpha})=0$. Consider the restriction morphism
\[
(\frakg/(\frakn \oplus \bk h_{\alpha} \oplus \bk e_{\alpha}))^* \xrightarrow{q} (\frakp_{\alpha}/(\frakn \oplus \bk h_{\alpha} \oplus \bk e_{\alpha}))^*.
\]
There are natural $P_{\alpha}$-actions on both spaces, and this morphism is $P_{\alpha}$-equivariant. Moreover, the $P_{\alpha}$-action on $(\frakp_{\alpha}/(\frakn \oplus \bk h_{\alpha} \oplus \bk e_{\alpha}))^*$ is trivial. Hence $q(P_{\alpha} \cdot \xi)=q(\xi)$, and $\dim(P_{\alpha} \cdot \xi) \leq \#(\Phi^+) - 1$. It follows that $\dim(G \cdot \xi) \leq \#\Phi - 2$, hence $\xi$ is not regular.\end{proof}

\begin{prop} \label{prop:actionreg}
 
There exists an action of $W$ on the variety $\wfrakg_{\reg}$, whose restriction to $\wfrakg_{\rs}$ is the action of Corollary {\rm \ref{cor:actionW}}.

\end{prop}

\begin{proof} Let $s \in \scS$, and let $Z_s^{\reg}$ be the restriction of $Z_s$ to $\wfrakg_{\reg}^2$. Consider the projection on the first component $p_1 : Z_s^{\reg} \to \wfrakg_{\reg}$. This morphism is proper, birational, with normal image. Moreover, it follows easily from the explicit description of $Z_s$ in \cite{RAct} and Lemma \ref{lem:regsl2} that it is bijective (see e.g.~Lemma \ref{lem:fibers} below). Hence it is an isomorphism of varieties. Similarly, the projection on the second component $p_2 : Z_s^{\reg} \to \wfrakg_{\reg}$ is an isomorphism. Let us denote by $f_s : \wfrakg_{\reg} \to \wfrakg_{\reg}$ the isomorphism given by the composition of $p_2$ with the inverse of $p_1$. We claim that the assignment \[ s \ \mapsto \ f_s \] extends to an action of the group $W$ on the variety $\wfrakg_{\reg}$.

First, for any $s \in \scS$, by symmetry of $Z_s$ under the exchange of the two copies of $\wfrakg$, $f_s$ is an involution. Hence we only have to check that these morphisms satisfy the braid relations. However, these morphisms stabilize the dense open subset $\wfrakg_{\rs} \subset \wfrakg_{\reg}$, and their restrictions satisfy the braid relations by Corollary \ref{cor:actionW}. Hence the braid relations are satisfied on the whole of $\wfrakg_{\reg}$. This finishes the proof.\end{proof}

The precise description of $\wfrakg_{\reg}$ will not be important for us. The only property that we need is the following. We prove it only under a very mild restriction on $p$. In this section, this assumption (which is probably not necessary) will only be used in the proof of the following proposition.

\begin{prop} \label{prop:codimreg}

Assume $p \neq 2$ if $G$ has a component of type $\mbfF_4$, and $p \neq 3$ if $G$ has a component of type $\mbfG_2$. Then the codimension of $\wfrakg \smallsetminus \wfrakg_{\reg}$ in $\wfrakg$ is at least $2$.

\end{prop}

\begin{proof} By $G$-equivariance, it is sufficient to prove that the codimension of $(\frakg/\frakn)^* \smallsetminus (\frakg/\frakn)^*_{\reg}$ in $(\frakg/\frakn)^*$ is at least $2$. We can assume that $G$ is quasi-simple.

First, we have
\begin{equation}
\label{eq:brs}
(\frakg/\frakn)^*_{\rs} = (\frakg/(\frakn \oplus \frakn^+))^*_{\rs} + (\frakg/\frakb)^*.
\end{equation}
Indeed, the inclusion $\subseteq$ follows from Lemma \ref{lem:rsb}. Now if $\xi$ is in the right-hand side, by the arguments of \cite[\S 3.8]{KW} there exists $u \in U$ such that $u \cdot \xi \in (\frakg/(\frakn \oplus \frakn^+))^*_{\rs}$. Hence $\xi \in (\frakg/\frakn)^*_{\rs}$.

By \eqref{eq:brs} and Lemma \ref{lem:rs}, the complement of $(\frakg/\frakn)^*_{\rs}$ in $(\frakg/\frakn)^*$ has components the \[ C_{\alpha}=\{ \, \xi \in (\frakg/\frakn)^* \mid \xi(h_{\alpha})=0 \, \} \] for $\alpha \in \Phi^+$. To prove the proposition, it is sufficient to prove that $(\frakg/\frakn)^*_{\reg}$ intersects any $C_{\alpha}$.

Assume that there is only a finite number of $G$-orbits in the ``dual nilpotent cone'' $\calN':=G \cdot (\frakg/\frakb)^*$. Then there exist regular nilpotent elements in $\frakg^*$, hence $(\frakg/\frakb)^*_{\reg} \neq \emptyset$. Moreover, elements of $(\frakg/\frakb)^*_{\reg}$ are obviously in $C_{\alpha} \cap (\frakg/\frakn)^*_{\reg}$ for any $\alpha$, which finishes the proof in this case. Unfortunately, we could not prove this property in general. Observe, however, that it is satisfied if $\frakg$ is a simple $G$-module (because in this case $\calN'$ identifies with the usual nilpotent cone in $\frakg$), for example if $G$ is of type $\mathbf{G}_2$ and $p \neq 3$ (see \cite[\S 0.13]{HUMConj}). This property is also proved if $G$ is of type $\mathbf{B}$, $\mathbf{C}$ or $\mathbf{D}$ and $p=2$ in \cite[\S 5.6]{Xu}.

The following alternative argument works if, in addition to the assumption of the proposition, we assume that $p \neq 2$ if $G$ is of type $\mbfB$, $\mbfC$ or $\mbfG$.

Fix some $\alpha \in \Phi^+$. We claim that, under our assumption, for $\beta \neq \alpha$ we have $h_{\alpha} \neq h_{\beta}$. Indeed for this fact we can assume that $\alpha$ is a simple root, using the action of $N_G(T)$. Then the claim follows from the facts that the $h_{\gamma}$'s, $\gamma \in \Sigma$, form a basis of $\frakt$ (because $G$ is simply connected), that $h_{\beta}=d(\beta^{\vee})_1(1)$, and from a look at the tables in \cite{BLie}. Hence one can consider some $\xi \in (\frakg/\frakn)^*$ such that $\xi(e_{\alpha})=1$, $\xi(h_{\alpha})=0$, $\xi(e_{\beta})=0$ for all $\beta \in \Phi^+ \smallsetminus \{\alpha\}$, and $\xi(h_{\beta}) \neq 0$ for all $\beta \in \Phi^+ \smallsetminus \{\alpha\}$. 

We claim that such a $\xi$ is regular. Indeed, the Jordan decomposition of $\xi$ (see \cite{KW}) is $\xi=\xi_s+\xi_n$, where $\xi_n$ is zero on $\frakt$ and on any $e_{\beta}$ for $\beta \neq \alpha$, and $\xi_n(e_{\alpha})=1$. By unicity of this decomposition, we have $Z_G(\xi)=Z_G(\xi_s) \cap Z_G(\xi_n)$. Now by construction and \cite[Lemma 3.1]{KW}, $L_{\alpha}:=Z_G(\xi_s)^{\circ}$ is a reductive subgroup of $G$ of semi-simple rank $1$, with Lie algebra $\frakg_{-\alpha} \oplus \frakt \oplus \frakg_{\alpha}$. And one easily checks that $\dim(Z_{L_{\alpha}}(\xi_n))=r$ (this is essentially a computation in $\mathrm{SL}(2,\bk)$), which finishes the proof of the proposition.\end{proof}

\begin{remark}
Let $\xi \in \frakg^*$ be such that $\xi_{|\frakg_{\alpha}}\neq 0$ for $\alpha \in \Sigma$, $\xi_{|\frakg_{\gamma}}=0$ for $\gamma \in \Phi \smallsetminus \Sigma$ and $\xi_{|\frakt}=0$. Then $\xi$ is nilpotent, and it follows from the proof of \cite[Lemma 5.13]{Xu} that $Z_G(\xi) \subset B$. Then it is easily checked that $Z_G(\xi) = Z(G) \times Z_U(\xi)$, where $Z(G)$ is the center of $G$. Hence, to prove that $\xi$ is regular, one only has to prove that $Z_U(\xi)$ has dimension $r$. It may be possible to use this observation to prove Proposition \ref{prop:codimreg} in the remaining two cases by direct computation, in the spirit of \cite{Ken}.
\end{remark}

\subsection{Kernels associated to simple reflections}
\label{ss:kernels-simple-reflections}

Fix $\alpha \in \Sigma$, and let $s=s_{\alpha}$. Consider the minimal standard parabolic subgroup $P_s$ associated to $\alpha$, and the corresponding partial flag variety $\calP_s=G/P_s$. Let $\frakp_s$ be the Lie algebra of $P_{\alpha}$, and $\frakp_s^{\mathrm{nil}}$ its nilpotent radical. Consider the variety
\[
\wfrakg_s \ := \ G \times^{P_s} (\frakg/\frakp_s^{\mathrm{nil}})^* \ = \ \{(X,gP_s) \in \frakg^* \times \calP_s \mid X_{|g \cdot \frakp_s^{\mathrm{nil}}}=0\}.
\]
There is a natural morphism $\widetilde{\pi}_s : \wfrakg \to \wfrakg_s$, and we consider the (scheme-theoretic) fiber product $\wfrakg \times_{\wfrakg_s} \wfrakg$. It is reduced, and has two irreducible components, $Z_s$ and $Z_1=\Delta \wfrakg$ (see \cite[\S 1.4]{RAct}). Recall that there exist exact sequences of $G \times \Gm$-equivariant coherent sheaves on $\wfrakg \times \wfrakg$
\begin{eqnarray}
\label{eq:exactsequence1}
& \calO_{\Delta \wfrakg} \langle -2 \rangle \ \hookrightarrow \ \calO_{\wfrakg \times_{\wfrakg_s} \wfrakg} \ \twoheadrightarrow \ \calO_{Z_s} ; & \\
\label{eq:exactsequence2}
& \calO_{Z_s}(-\rho, \rho-\alpha) \ \hookrightarrow \ \calO_{\wfrakg \times_{\wfrakg_s} \wfrakg} \ \twoheadrightarrow \ \calO_{\Delta \wfrakg}, &
\end{eqnarray}
see \cite[Corollary 5.3.2]{RAct}. Let us observe that for both sequences the surjection is induced by restriction of functions, and that the second sequence is induced by the natural exact sequence of sheaves on $\calB \times \calB$
\[
\calO_{\calB \times_{\calP_s} \calB}(-\rho,\rho-\alpha) \ \hookrightarrow \ \calO_{\calB \times_{\calP_s} \calB} \ \twoheadrightarrow \ \calO_{\Delta \calB}.
\]

These exact sequences allow to give a simpler proof of statement \eqref{it:inverse-2} of the following proposition, which was proved by explicit computation in \cite[Proposition 1.5.2]{RAct}.

\begin{prop} \label{prop:inverse}

\begin{enumerate}

\item
\label{it:inverse-1}
There exist isomorphisms
\[
\calO_{\wfrakg \times_{\wfrakg_{s}} \wfrakg} \star \calO_{Z_s} \, \cong \, \calO_{\wfrakg \times_{\wfrakg_{s}} \wfrakg}, \quad \calO_{Z_s} \star \calO_{\wfrakg \times_{\wfrakg_{s}} \wfrakg} \, \cong \, \calO_{\wfrakg \times_{\wfrakg_{s}} \wfrakg}
\]
in the category $\calD^b_{\mathrm{prop}} \Coh^{G \times \Gm}(\wfrakg \times \wfrakg)$.
\item 
\label{it:inverse-2}
The kernel $\calO_{Z_s}\langle 1 \rangle$ is invertible, with inverse $\calO_{Z_s}(-\rho, \rho-\alpha)\langle 1 \rangle$. In other words, there exist isomorphisms in $\calD^b_{\mathrm{prop}} \Coh^{G \times \Gm}(\wfrakg \times \wfrakg)$:
\[
\calO_{Z_s}(-\rho, \rho-\alpha) \star \calO_{Z_s} \, \cong \, \calO_{\Delta \wfrakg} \langle -2 \rangle, \quad \calO_{Z_s} \star \calO_{Z_s}(-\rho, \rho-\alpha) \, \cong \, \calO_{\Delta \wfrakg} \langle -2 \rangle.
\]

\end{enumerate}

\end{prop}

\begin{proof} 
\eqref{it:inverse-1} We only prove the first isomorphism, the second one can be treated similarly. Arguments similar to those of \cite[Proposition 5.2.2]{RAct} show that there is an isomorphism
\[
\calO_{\wfrakg \times_{\wfrakg_s} \wfrakg} \star \calO_{Z_s} \, \cong \, L(\Id \times \widetilde{\pi}_s)^* \circ R(\Id \times \widetilde{\pi}_s)_* \calO_{Z_s}.
\]
Now the morphism
\[
\Id \times \widetilde{\pi}_s : Z_s \ \to \ \wfrakg \times \wfrakg_s
\]
is proper, has normal image (indeed, this image coincides with the graph $\wfrakg \times_{\wfrakg_s} \wfrakg_s$ of $\widetilde{\pi}_s$, hence is isomorphic to $\wfrakg$), and is birational when considered as a morphism from $Z_s$ to this image. Hence, by Zariski's Main Theorem, $(\Id \times \widetilde{\pi}_s)_* \calO_{Z_s} \cong \calO_{\wfrakg \times_{\wfrakg_s} \wfrakg_s}$. Moreover, by the same arguments as in \cite[Proposition 2.4.1]{RAct}, $R^i (\Id \times \widetilde{\pi}_s)_* \calO_{Z_s} = 0$ for $i \geq 1$. To finish the proof, one observes that we have $L(\Id \times \widetilde{\pi}_s)^* \calO_{\wfrakg \times_{\wfrakg_s} \wfrakg_s} \cong \calO_{\wfrakg \times_{\wfrakg_s} \wfrakg}$. (The isomorphism in cohomological degree $0$ is the definition of the fiber product; the higher vanishing for inverse image follows from a dimension argument.)

\eqref{it:inverse-2} Again, we only prove the first isomorphism. Consider the exact sequence \eqref{eq:exactsequence2}, and convolve it with $\calO_{Z_s}$ on the right; we obtain a distinguished triangle in $\calD^b_{\mathrm{prop}} \Coh^{G \times \Gm}(\wfrakg \times \wfrakg)$:
\[
\calO_{Z_s}(-\rho, \rho-\alpha) \star \calO_{Z_s} \ \to \ \calO_{\wfrakg \times_{\wfrakg_s} \wfrakg} \star \calO_{Z_s} \ \to \ \calO_{\Delta \wfrakg} \star \calO_{Z_s} \xrightarrow{+1}.
\]
By \eqref{it:inverse-1}, the middle term is isomorphic to $\calO_{\wfrakg \times_{\wfrakg_s} \wfrakg}$, and the term on the right hand side is isomorphic to $\calO_{Z_s}$. Moreover, the morphism on the right identifies with the restriction of functions $\calO_{\wfrakg \times_{\wfrakg_s} \wfrakg} \to \calO_{Z_s}$. Hence this triangle can be identified with the one associated to the exact sequence \eqref{eq:exactsequence1}. In particular we obtain the expected isomorphism.\end{proof}

\begin{remark}
\label{rk:canonicalsheaf}
We have remarked above that $Z_s$ is a smooth variety. Statement \eqref{it:inverse-2} implies that its canonical sheaf is $\calO_{Z_s}(-\rho, \rho-\alpha)$ (see \cite{Huy}).
\end{remark}

To prove Theorem \ref{thm:existenceaction} for $\wfrakg$, we have to check that the kernels $\calO_{Z_s}$ ($s \in \scS$) and $\calO_{\Delta \wfrakg}(x)$ ($x \in \bbX$) satisfy the relations of the presentation of $\bB_{\aff}$ given in \S \ref{ss:notation} in the monoidal category $\calD^b \Coh^{G \times \Gm}(\wfrakg \times \wfrakg)$, up to isomorphism. It is explained in \cite[\S 1.6]{RAct} that relations $\rmii$ and $\rmiii$ are trivial, and that relations $\rmiv$ follow easily from Proposition \ref{prop:inverse}\eqref{it:inverse-2}. In the rest of this section we give a proof of relations $\rmi$ (finite braid relations), inspired by the methods of \cite{BEZCoh}. It is explained in \cite[\S 4]{RAct} how to deduce the theorem for $\wcalN$ from the case of $\wfrakg$. Note in particular that if $s=s_{\alpha}$ ($\alpha \in \Sigma$), the inverse of $\calO_{Z_s'} \langle 1 \rangle$ for the convolution product is $\calO_{Z_s'}(-\rho,\rho-\alpha) \langle 1 \rangle$. The compatibility of the actions with the functor $Ri_*$ easily follows from the definition and the projection formula. The compatibility with $Li^*$ follows by adjunction.

\subsection{Line bundles on $\wfrakg$}

Let us begin with some generalities on line bundles on $\wfrakg$. First, the following lemma immediately follows from \cite[II.8.5(1)]{JANAlg}.

\begin{lem}
\label{lem:lemample}

Let $\lambda \in \bbX$, such that $\lambda - \rho$ is dominant. Then $\calO_{\wfrakg}(\lambda)$ is an ample line bundle on $\wfrakg$.\qed

\end{lem}

Next, let $\bbX^+ \subset \bbX$ be the set of dominant weights. Consider the $\bbX$-graded $G \times \Gm$-equivariant algebra
\[
\mathbf{\Gamma}(\wfrakg) \ := \ \bigoplus_{\lambda \in \bbX^+} \, \Gamma(\wfrakg, \calO_{\wfrakg}(\lambda)).
\]
There is a natural functor
\[
\mathbf{\Gamma} : \left\{
\begin{array}{ccc}
\QCoh^{G \times \Gm}(\wfrakg) & \to & \Mod^{G \times \Gm}_{\bbX}(\mathbf{\Gamma}(\wfrakg)), \\
\calM & \mapsto & \bigoplus_{\lambda \in \bbX^+} \, \Gamma(\wfrakg, \calM \otimes_{\calO_{\wfrakg}} \calO_{\wfrakg}(\lambda))
\end{array}
\right. ,
\]
where $\Mod^{G \times \Gm}_{\bbX}(\mathbf{\Gamma}(\wfrakg))$ is the category of $\bbX$-graded $G \times \Gm$-equivariant modules over $\mathbf{\Gamma}(\wfrakg)$. We let $\mathrm{Tor}^{G \times \Gm}_{\bbX}(\mathbf{\Gamma}(\wfrakg)) \subset \Mod^{G \times \Gm}_{\bbX}(\mathbf{\Gamma}(\wfrakg))$ be the subcategory of objects which are direct limits of objects $M$ such that there exists $\mu \in \bbX$ such that the $\lambda$-component of $M$ is zero for any $\lambda \in \mu + \bbX^+$. As the morphism $\wfrakg \to \frakg^*$ is projective and using Lemma \ref{lem:lemample}, we have the following version of Serre's theorem. (To prove this result, one can e.g.~adapt the arguments of the proof of \cite[Theorem 1.3]{AV}.)

\begin{prop}
\label{prop:Serre-thm}

The composition
\[
\QCoh^{G \times \Gm}(\wfrakg) \xrightarrow{\mathbf{\Gamma}} \Mod^{G \times \Gm}_{\bbX}(\mathbf{\Gamma}(\wfrakg)) \to \Mod^{G \times \Gm}_{\bbX}(\mathbf{\Gamma}(\wfrakg)) / \mathrm{Tor}^{G \times \Gm}_{\bbX}(\mathbf{\Gamma}(\wfrakg))
\]
is an equivalence of abelian categories.\qed

\end{prop}

If $A$ is any subset of $\bbX$, we denote by $\calD_A$ the smallest strictly full triangulated subcategory of $\calD^b \Coh^{G \times \Gm}(\wfrakg)$ containing the line bundles $\calO_{\wfrakg}(\lambda)$ for $\lambda \in A$ and stable under the functors $\langle j \rangle$ for $j \in \Z$. We denote by ${\rm conv}(\lambda)$ the intersection of $\bbX$ with the convex hull of $W \cdot \lambda$, and by ${\rm conv}^0(\lambda)$ the complement of $W \cdot \lambda$ in ${\rm conv}(\lambda)$.

In the next lemma we will also use the following notation. If $\mathscr{B}$ is a triangulated category, and $\mathscr{A} \subset \mathscr{B}$ is a full triangulated subcategory, for $M, N \in \mathscr{B}$ we write $M \cong N \ {\rm mod} \ \mathscr{A}$ if the images of $M$ and $N$ in the quotient category $\mathscr{B} / \mathscr{A}$ (in the sense of Verdier) are isomorphic.

\begin{lem} \label{lem:lemreflections}

Let $\alpha \in \Sigma$, and $s=s_{\alpha}$.

\begin{enumerate}
\item
\label{it:lemreflections-1}
For any $\lambda \in \bbX$, the functors $F_{\wfrakg}^{\calO_{Z_s}}$, $F_{\wfrakg}^{\calO_{Z_s}(-\rho, \rho - \alpha)}$ preserve the subcategory $\calD_{{\rm conv(\lambda)}}$.
\item 
\label{it:lemreflections-2}
Let $\lambda \in \bbX$ such that $\langle \lambda, \alpha^{\vee} \rangle = 0$. Then we have
\[
F_{\wfrakg}^{\calO_{Z_s}}(\calO_{\wfrakg}(\lambda)) \ \cong \ \calO_{\wfrakg}(\lambda), \quad F_{\wfrakg}^{\calO_{Z_s(-\rho,\rho-\alpha)}}(\calO_{\wfrakg}(\lambda)) \ \cong \ \calO_{\wfrakg}(\lambda) \langle -2 \rangle.
\]
\item 
\label{it:lemreflections-3}
Let $\lambda \in \bbX$ such that $\langle \lambda, \alpha^{\vee} \rangle < 0$. Then
\[
F_{\wfrakg}^{\calO_{Z_s}}(\calO_{\wfrakg}(\lambda)) \ \cong \ \calO_{\wfrakg}(s \lambda) \ \langle -2 \rangle \quad {\rm mod} \ \calD_{{\rm conv}^0(\lambda)}.
\]
\item
\label{it:lemreflections-4}
Let $\lambda \in \bbX$ such that $\langle \lambda, \alpha^{\vee} \rangle > 0$. Then
\[
F_{\wfrakg}^{\calO_{Z_s}(-\rho, \rho - \alpha)}(\calO_{\wfrakg}(\lambda)) \ \cong \ \calO_{\wfrakg}(s \lambda) \quad {\rm mod} \ \calD_{{\rm conv}^0(\lambda)}.
\]

\end{enumerate}

\end{lem}

\begin{proof} Recall the notation for $P_s$, $\calP_s:=G/P_s$, $\wfrakg_s$ (see \S \ref{ss:kernels-simple-reflections}). The variety $\wfrakg_s$ is endowed with a natural $G \times \Gm$-action, such that the morphism $\widetilde{\pi}_s: \wfrakg \to \wfrakg_s$ is $G \times \Gm$-equivariant. Using exact sequences \eqref{eq:exactsequence1} and \eqref{eq:exactsequence2} and \cite[Proposition 5.2.2]{RAct}, for any $\calF$ in $\calD^b \Coh^{G \times \Gm}(\wfrakg)$ there exist distinguished triangles
\begin{eqnarray}
\label{eq:trianglefunctors1}
& \calF \langle -2 \rangle \ \to \ L(\widetilde{\pi}_s)^* \circ R(\widetilde{\pi}_s)_* \calF \ \to \ F_{\wfrakg}^{\calO_{Z_s}}(\calF) \ \xrightarrow{+1}; & \\
\label{eq:trianglefunctors2}
& F_{\wfrakg}^{\calO_{Z_s}(-\rho, \rho - \alpha)}(\calF) \ \to \ L(\widetilde{\pi}_s)^* \circ R(\widetilde{\pi}_s)_* \calF \ \to \ \calF \ \xrightarrow{+1}. &
\end{eqnarray}

Let $j: \wfrakg \hookrightarrow \wfrakg_s \times_{\calP_s} \calB$ be the natural inclusion. There exists an exact sequence
\begin{equation}
\label{eq:exactsequence}
\calO_{\wfrakg_s \times_{\calP_s} \calB}(-\alpha) \langle -2 \rangle \ \hookrightarrow \ \calO_{\wfrakg_s \times_{\calP_s} \calB} \ \twoheadrightarrow \ j_* \calO_{\wfrakg}.
\end{equation}
(Indeed, $(\frakg/\frakn)^* \subset (\frakg/\frakp_s^{\mathrm{nil}})^*$ is defined by one equation, of weight $(-\alpha, -2)$ for $B \times \Gm$.) Let also $p : \wfrakg_s \times_{\calP_s} \calB \to \wfrakg_s$ be the projection. Then $\widetilde{\pi}_s=p \circ j$.

Using triangles \eqref{eq:trianglefunctors1} and \eqref{eq:trianglefunctors2}, to prove \eqref{it:lemreflections-1} it is sufficient to prove that for any $\lambda \in \bbX$, $L(\widetilde{\pi}_s)^* \circ R(\widetilde{\pi}_s)_* \calO_{\wfrakg}(\lambda)$ is in $\calD_{{\rm conv}(\lambda)}$. The case $\langle \lambda, \alpha^{\vee} \rangle=0$ is trivial: in this case we have $L(\widetilde{\pi}_s)^* \circ R(\widetilde{\pi}_s)_* \calO_{\wfrakg}(\lambda) \cong \calO_{\wfrakg}(\lambda) \oplus \calO_{\wfrakg}(\lambda) \langle 2 \rangle$ by the projection formula. Here we have used the well-known isomorphism
\[
R(\widetilde{\pi}_s)_* \calO_{\wfrakg} \cong \calO_{\wfrakg_s} \oplus \calO_{\wfrakg_s} \langle -2 \rangle.
\]
The property \eqref{it:lemreflections-2} also follows, using triangles \eqref{eq:trianglefunctors1} and \eqref{eq:trianglefunctors2}.

Now, assume that $\langle \lambda, \alpha^{\vee} \rangle > 0$. Tensoring \eqref{eq:exactsequence} by $\calO_{\wfrakg_s \times_{\calP_s} \calB}(\lambda)$ we obtain an exact sequence
\[
\calO_{\wfrakg_s \times_{\calP_s} \calB}(\lambda - \alpha) \langle -2 \rangle \ \hookrightarrow \ \calO_{\wfrakg_s \times_{\calP_s} \calB}(\lambda) \ \twoheadrightarrow \ j_* \calO_{\wfrakg}(\lambda).
\]
Then, applying the functor $Rp_*$ and using \cite[\S I.5.19, Proposition II.5.2.(c)]{JANAlg} we obtain a distinguished triangle 
\[
\calL_{\wfrakg_s}(\Ind_B^{P_s}(\lambda - \alpha)) \langle -2 \rangle \ \to \ \calL_{\wfrakg_s}(\Ind_B^{P_s}(\lambda)) \ \to \ R(\widetilde{\pi}_s)_* \calO_{\wfrakg}(\lambda) \ \xrightarrow{+1}.
\]
(Observe that here $\langle \lambda - \alpha, \alpha^{\vee} \rangle \geq -1$.) Applying the functor $L(\widetilde{\pi}_s)^*$ we obtain a triangle
\begin{equation*}
\calL_{\wfrakg}(\Ind_B^{P_s}(\lambda - \alpha)) \langle -2 \rangle \ \to \ \calL_{\wfrakg}(\Ind_B^{P_s}(\lambda)) \ \to \ L(\widetilde{\pi}_s)^* \circ R(\widetilde{\pi}_s)_* \calO_{\wfrakg}(\lambda) \ \xrightarrow{+1}.
\end{equation*}
Now it is well known (see again \cite[Proposition II.5.2.(c)]{JANAlg}) that the $P_s$-module $\Ind_B^{P_s}(\lambda)$ has weights $\lambda, \ \lambda - \alpha, \ \cdots, \ s \lambda$. Hence $\calL_{\wfrakg}(\Ind_B^{P_s}(\lambda))$ has a filtration with subquotients $\calO_{\wfrakg}(\lambda)$, $\calO_{\wfrakg}(\lambda-\alpha)$, $\cdots$, $\calO_{\wfrakg}(s \lambda)$. Similarly, $\calL_{\wfrakg}(\Ind_B^{P_s}(\lambda - \alpha))$ has a filtration with subquotients $\calO_{\wfrakg}(\lambda-\alpha)$, $\cdots$, $\calO_{\wfrakg}(s \lambda + \alpha)$. This proves \eqref{it:lemreflections-1} in this case, and also \eqref{it:lemreflections-4}.

Now assume $\langle \lambda, \alpha^{\vee} \rangle < 0$. Using similar arguments, there exists a distinguished triangle
\[
L(\widetilde{\pi}_s)^* \circ R(\widetilde{\pi}_s)_* \calO_{\wfrakg}(\lambda) \ \to \ \calL_{\wfrakg}(R^1 \Ind_B^{P_s}(\lambda - \alpha)) \langle -2 \rangle \ \to \ \calL_{\wfrakg}(R^1 \Ind_B^{P_s}(\lambda)).
\]
Moreover, $\calL_{\wfrakg}(R^1 \Ind_B^{P_s}(\lambda))$ has a filtration with subquotients $\calO_{\wfrakg}(s \lambda - \alpha)$, $\cdots$, $\calO_{\wfrakg}(\lambda + \alpha)$, and $\calL_{\wfrakg}(R^1 \Ind_B^{P_s}(\lambda - \alpha))$ has a filtration with subquotients $\calO_{\wfrakg}(s \lambda)$, $\cdots$, $\calO_{\wfrakg}(\lambda)$. As above, this proves \eqref{it:lemreflections-1} in this case, and \eqref{it:lemreflections-3}. \end{proof}

\begin{remark} 
\begin{enumerate}
\item The case \eqref{it:lemreflections-2} of the proposition is not needed for our arguments. We only include it for completeness.
\item 
\label{it:remark-2}
It follows from the proof of the proposition that if $\langle \lambda,\alpha^{\vee} \rangle=1$, then the isomorphism of \eqref{it:lemreflections-4} can be lifted to an isomorphism in $\calD^b \Coh^{G \times \Gm}(\wfrakg)$. Similarly, if $\langle \lambda,\alpha^{\vee} \rangle=-1$, the isomorphism of \eqref{it:lemreflections-3} can be lifted to an isomorphism in $\calD^b \Coh^{G \times \Gm}(\wfrakg)$.
\item If $\langle \lambda,\alpha^{\vee} \rangle \notin \{-1,0,1 \}$, we do not have an explicit description of the objects $F_{\wfrakg}^{\calO_{Z_s}}(\calO_{\wfrakg}(\lambda))$ or $F_{\wfrakg}^{\calO_{Z_s}(-\rho,\rho-\alpha)}(\calO_{\wfrakg}(\lambda))$ as in \eqref{it:remark-2}. The proof of the proposition gives a recipe for computing their class in equivariant $K$-theory, however. The answer can be given in terms of Demazure--Lusztig operators as in \cite[Theorem 7.2.16]{CG}. (See \S \ref{ss:homology} below for more details in this direction.)

\end{enumerate}
\end{remark}

The following lemma is a generalization of \cite[Lemma 5]{BEZCoh} (where it is assumed that $p=0$). The proof is similar.

\begin{lem} \label{lem:lemorthogonality}

Let $\lambda, \mu \in \bbX$.

We have $\Ext^{\bullet}_{\calD^b \Coh^{G}(\wfrakg)}(\calO_{\wfrakg}(\lambda), \calO_{\wfrakg}(\mu))=0$ unless $\lambda - \mu \in \mathbb{Z}_{\geq 0} R^+$.

Similarly, for any $i \in \mathbb{Z}$ we have $\Ext^{\bullet}_{\calD^b \Coh^{G \times \Gm}(\wfrakg)}(\calO_{\wfrakg}(\lambda), \calO_{\wfrakg}(\mu) \langle i \rangle)=0$ unless $\lambda - \mu \in \mathbb{Z}_{\geq 0} R^+$.

\end{lem}

\begin{proof} We give a proof only in the first case. Recall that $\calD^b \Coh^{G}(\wfrakg)$ is equivalent to the full subcategory of $\calD^b \QCoh^{G}(\wfrakg)$ whose objects have coherent cohomology (see \cite[Corollary 2.11]{BEZPer}). Hence we can replace $\calD^b \Coh^{G}(\wfrakg)$ by $\calD^b \QCoh^{G}(\wfrakg)$ in the statement. Moreover, for any $i \in \mathbb{Z}$ there is a natural isomorphism \begin{equation} \label{eq:isomext} \Ext^i_{\calD^b \QCoh^{G}(\wfrakg)}(\calO_{\wfrakg}(\lambda), \calO_{\wfrakg}(\mu)) \cong H^i \bigl( R(\Gamma^G)(\calO_{\wfrakg}(\mu - \lambda)) \bigr), \end{equation} where $\Gamma^G$ denotes the functor which sends a $G$-equivariant quasi-coherent sheaf $\calF$ to the $G$-invariants in its global sections, and $R(\Gamma^G)$ is its derived functor.

Recall also that, by definition, we have $\wfrakg = G \times^B (\frakg / \frakn)^*$. Hence the restriction functor $\calF \mapsto \calF|_{\{1\} \times (\frakg/\frakn)^*}$ induces an equivalence of categories \[ \QCoh^G(\wfrakg) \xrightarrow{\sim} \QCoh^B\bigl((\frakg/\frakn)^* \bigr) \] (see e.g.~\cite[\S 2]{BRIMul}). Moreover, the following diagram commutes, where $\Gamma^B$ is defined as $\Gamma^G$ above: \[ \xymatrix{\QCoh^G(\wfrakg) \ar[d]^-{\wr} \ar[rrd]^-{\Gamma^G} & & \\ \QCoh^B((\frakg / \frakn)^*) \ar[rr]^{\Gamma^B} & & {\rm Vect}(\bk). } \] It follows, using isomorphism \eqref{eq:isomext}, that for any $i \in \mathbb{Z}$ we have \begin{equation} \label{eq:isomext2} \Ext^i_{\calD^b \QCoh^{G}(\wfrakg)}(\calO_{\wfrakg}(\lambda), \calO_{\wfrakg}(\mu)) \cong H^i \bigl( R(\Gamma^B)(\calO_{(\frakg / \frakn)^*} \otimes_{\bk} \bk_B(\mu - \lambda)) \bigr). \end{equation}

The functor $\Gamma^B$ is the composition of the functor $$\Gamma((\frakg / \frakn)^*, -) : \QCoh^B((\frakg/\frakn)^*) \xrightarrow{\sim} \Mod^B({\rm S}(\frakg / \frakn)),$$ which is an equivalence of categories because $(\frakg / \frakn)^*$ is affine, and the $B$-fixed points functor $I^B : \Mod^B({\rm S}(\frakg / \frakn)) \to {\rm Vect(\bk)}$. (Here, ${\rm S}(\frakg / \frakn)$ is the symmetric algebra of the vector space $\frakg / \frakn$.) Hence, using isomorphism \eqref{eq:isomext2} we deduce that for any $i \in \mathbb{Z}$ we have \begin{equation} \label{eq:isomext3} \Ext^i_{\calD^b \QCoh^{G}(\wfrakg)}(\calO_{\wfrakg}(\lambda), \calO_{\wfrakg}(\mu)) \cong H^i \bigl( R(I^B)({\rm S}(\frakg / \frakn) \otimes_{\bk} \bk_B(\mu - \lambda)) \bigr). \end{equation}

Now $I^B$ is the composition of the forgetful functor $\For : \Mod^B({\rm S}(\frakg / \frakn)) \to {\rm Rep}(B)$ and the $B$-fixed points functor $J^B : {\rm Rep}(B) \to {\rm Vect}(\bk)$. Of course the functor $\For$ is exact, and in the category $\Mod^B({\rm S}(\frakg / \frakn))$ there are enough objects of the form $\Ind_{\{1\}}^B(M) \cong M \otimes_{\bk} \bk[B]$, for $M$ a ${\rm S}(\frakg / \frakn)$-module, whose images under $\For$ are acyclic for the functor $J^B$. Hence for any $i \in \mathbb{Z}$ we have \begin{equation} \label{eq:isomext4} \Ext^i_{\calD^b \QCoh^{G}(\wfrakg)}(\calO_{\wfrakg}(\lambda), \calO_{\wfrakg}(\mu)) \cong H^i \bigl( R(J^B)({\rm S}(\frakg / \frakn) \otimes_{\bk} \bk_B(\mu - \lambda)) \bigr), \end{equation} where for simplicity we have omitted the functor $\For$.

Finally, as $B \cong T \ltimes U$, the functor $J^B$ is the composition of the $U$-fixed points functor $J^U$, followed by the $T$-fixed points functor $J^T$ (which is exact). Hence $RJ^B \cong J^T \circ RJ^U$, and we only have to prove that \begin{equation} \label{eq:vanishing} J^T \bigl( R(J^U)({\rm S}(\frakg / \frakn) \otimes_{\bk} \bk_B(\nu)) \bigr) = 0 \end{equation} unless $\nu$ is a sum of negative roots. But $R(J^U)({\rm S}(\frakg / \frakn) \otimes_{\bk} \bk_B(\nu))$ can be computed by the Hochschild complex $C(U,\, {\rm S}(\frakg / \frakn) \otimes_{\bk} \bk_B(\nu))$ (see \cite[I.4.16]{JANAlg}). And the $T$-weights of this complex are all in $\mathbb{Z}_{\geq 0} R^+$ (because all weights of ${\rm S}(\frakg / \frakn)$ and of $\bk[U]$ are in $\mathbb{Z}_{\geq 0} R^+$). Then \eqref{eq:vanishing} easily follows.\end{proof}

\subsection{Braid relations}
\label{ss:braid-relations}

\begin{prop} \label{prop:propfunctorsbraidrel}

Let $\alpha, \beta \in \Sigma$, and $s=s_{\alpha}$, $t=s_{\beta}$. For any dominant weight $\lambda$ we have an isomorphism
\begin{multline*}
F_{\wfrakg}^{\calO_{Z_s}(-\rho, \rho - \alpha)} \circ F_{\wfrakg}^{\calO_{Z_t}(-\rho, \rho-\beta)} \circ \cdots \bigl( \calO_{\wfrakg}(\lambda) \bigr) \\
\cong \ F_{\wfrakg}^{\calO_{Z_t}(-\rho, \rho - \beta)} \circ F_{\wfrakg}^{\calO_{Z_s}(-\rho, \rho - \alpha)} \circ \cdots \bigl( \calO_{\wfrakg}(\lambda) \bigr),
\end{multline*}
in $\calD^b \Coh^{G \times \Gm}(\wfrakg)$, where the number of functors appearing on each side is $n_{s,t}$.

\end{prop}

\begin{proof} To fix notation, let us assume that $\alpha$ and $\beta$ generate a sub-system of type $\mbfA_2$. (The proof is similar in the other cases.) By Proposition \ref{prop:inverse}\eqref{it:inverse-2} we have an isomorphism of functors 
\[
\bigl( F_{\wfrakg}^{\calO_{Z_s}} \bigr)^{-1} \ \cong \ F_{\wfrakg}^{\calO_{Z_s}(-\rho, \rho - \alpha)} \langle 2 \rangle,
\]
and similarly for $\beta$. Hence proving the proposition is equivalent to proving that
\begin{multline}
\label{eq:braidrelfunctor}
E_{\lambda} := F_{\wfrakg}^{\calO_{Z_s}} \circ F_{\wfrakg}^{\calO_{Z_t}} \circ F_{\wfrakg}^{\calO_{Z_s}} \circ F_{\wfrakg}^{\calO_{Z_t}(-\rho, \rho - \beta)} \\ \circ F_{\wfrakg}^{\calO_{Z_s}(-\rho, \rho - \alpha)} \circ F_{\wfrakg}^{\calO_{Z_t}(-\rho, \rho - \beta)} (\calO_{\wfrakg}(\lambda)) 
\end{multline}
is isomorphic to $\calO_{\wfrakg}(\lambda) \langle -6 \rangle$. First, it follows from Lemma \ref{lem:lemreflections} that $E_{\lambda} \in \calD_{{\rm conv}(\lambda)}$ and that
\begin{equation}
\label{eq:isomquotient}
E_{\lambda} \ \cong \ \calO_{\wfrakg}(\lambda) \langle -6 \rangle \ {\rm mod} \ \calD_{{\rm conv}^0(\lambda)}.
\end{equation}
Hence $E_{\lambda}$ is even in $\calD_{{\rm conv}^0 (\lambda) \cup \{\lambda\}}$.

For any full subcategory $\mathscr{A}$ of a category $\mathscr{B}$, we denote by $(\mathscr{A}^{\bot})_{\mathscr{B}}$ the full subcategory of $\mathscr{B}$ with objects the $M$ such that $\Hom_{\mathscr{B}}(A,M)=0$ for any $A$ in $\mathscr{A}$. By Lemma \ref{lem:lemorthogonality}, $\calO_{\wfrakg}(\lambda)$ is in $(\calD_{{\rm conv}^0(\lambda)}^{\bot})_{\calD_{{\rm conv}^0 (\lambda) \cup \{\lambda\}}}$. Hence, as all the functors involved preserve the subcategory $(\calD_{{\rm conv}^0(\lambda)}^{\bot})_{\calD_{{\rm conv}(\lambda)}}$ (because their inverse preserves $\calD_{{\rm conv}^0(\lambda)}$ by Lemma \ref{lem:lemreflections}), also $E_{\lambda}$ is in the subcategory $(\calD_{{\rm conv}^0(\lambda)}^{\bot})_{\calD_{{\rm conv}^0 (\lambda) \cup \{\lambda\}}}$. Now it follows easily from \cite[Propositions 1.5 and 1.6]{BKRep} that the projection
\[
(\calD_{{\rm conv}^0(\lambda)}^{\bot})_{\calD_{{\rm conv}^0 (\lambda) \cup \{\lambda\}}} \to \calD_{{\rm conv}^0 (\lambda) \cup \{\lambda\}} / \calD_{{\rm conv}^0 (\lambda)}
\]
is an equivalence of categories. Using again \eqref{eq:isomquotient}, we deduce that $E_{\lambda} \cong \calO_{\wfrakg}(\lambda) \langle -6 \rangle$ in $\calD^b \Coh^{G \times \Gm}(\wfrakg)$, as claimed. \end{proof}

Before the next corollary we introduce some notation. If $\lambda$ is a dominant weight, we write that a property is true for $\lambda \gg 0$ if there exists a positive integer $N$ such that the property is true for any weight $\lambda$ such that $\langle \lambda, \alpha^{\vee} \rangle \geq N$ for any positive root $\alpha$.

\begin{cor} \label{cor:braidrelations}

The kernels $\calO_{Z_s}$, $s \in \scS$, satisfy the finite braid relations in the monoidal category $\calD^b \Coh_{\pro}(\wfrakg \times \wfrakg)$. More precisely, for $s,t \in \scS$ there exists a canonical isomorphism 
\[
\calO_{Z_s} \star \calO_{Z_t} \star \cdots \ \cong \ \calO_{Z_t} \star \calO_{Z_s} \star \cdots,
\]
in $\calD^b \Coh^{G \times \Gm}(\wfrakg \times \wfrakg)$, where the number of terms on each side is $n_{s,t}$.

\end{cor}

\begin{proof} To fix notation, let us assume that $s$ and $t$ generate a subgroup of $W$ of type $\mbfA_2$. (The other cases are similar.) The kernel $\calO_{Z_s}$ is invertible (see Proposition \ref{prop:inverse}\eqref{it:inverse-2}), with inverse
\[
(\calO_{Z_s})^{-1} \ := \ R\sheafHom_{\calO_{\wfrakg \times \wfrakg}}(\calO_{Z_s}, \calO_{\wfrakg \times \wfrakg}) \otimes_{\calO_{\wfrakg \times \wfrakg}} p_2^* \omega_{\wfrakg} [\dim(\frakg)].
\]
The same is true for $t$ instead of $s$. Hence we only have to prove that 
\[
\calO_{Z_s} \star \calO_{Z_t} \star \calO_{Z_s} \star (\calO_{Z_t})^{-1} \star (\calO_{Z_s})^{-1} \star (\calO_{Z_t})^{-1} \ \cong \ \calO_{\Delta \wfrakg}.
\]
For simplicity, let us denote by $\calK_{s,t}$ the object on the left hand side of this equation.

First, let $j : \wfrakg_{\reg} \times \wfrakg_{\reg} \hookrightarrow \wfrakg \times \wfrakg$ be the inclusion. We claim that there exists a canonical isomorphism
\begin{equation}
\label{eq:reg-part}
j^* \calK_{s,t} \ \cong \ \calO_{\Delta \wfrakg_{\reg}}.
\end{equation}
Indeed, the functor $j^*$ is monoidal. And there exist isomorphism
\begin{align*}
 j^* \calO_{Z_s} \ & \cong \ \calO_{Z_s^{\reg}} & \\
 j^* (\calO_{Z_s})^{-1} \cong R\sheafHom_{\calO_{\wfrakg_{\reg} \times \wfrakg_{\reg}}}(\calO_{Z_s^{\reg}}, & \calO_{\wfrakg_{\reg} \times \wfrakg_{\reg}}) \otimes_{\calO_{\wfrakg_{\reg} \times \wfrakg_{\reg}}} p_2^* \omega_{\wfrakg_{\reg}} [\dim(\frakg)], &
\end{align*}
where as above $Z_s^{\reg}=Z_s \cap (\wfrakg_{\reg} \times \wfrakg_{\reg})$. For simplicity, we denote the right hand side of the second line by $(\calO_{Z_s^{\reg}})^{-1}$. The same formulas also hold for $t$ instead of $s$. Now there exists a canonical isomorphism
\[
\calO_{Z_s^{\reg}} \star \calO_{Z_t^{\reg}} \star \calO_{Z_s^{\reg}} \ \cong \ \calO_{Z_t^{\reg}} \star \calO_{Z_s^{\reg}} \star \calO_{Z_t^{\reg}},
\]
because each side is canonically isomorphism to functions on the graph of the action of $sts=tst$. (Here we use Proposition \ref{prop:actionreg}.) By usual adjunction properties for Fourier-Mukai kernels, this isomorphism induces a canonical morphism
\[
j^* \calK_{s,t} \to \calO_{\wfrakg_{\reg}},
\]
which also is an isomorphism.

To prove the isomorphism of the corollary it is sufficient, using Lemma \ref{lem:lemample} and Proposition \ref{prop:Serre-thm}, to prove that for $\lambda,\mu \gg 0$ we have $R^{\neq 0}\Gamma(\wfrakg \times \wfrakg, \ \calK_{s,t}(\lambda, \mu))=0$ (this implies that $\calK_{s,t}$ is concentrated in degree $0$, i.e.~is a sheaf), and that there exist canonical isomorphisms
\[
\Gamma(\wfrakg \times \wfrakg, \, \calK_{s,t}(\lambda, \mu)) \ \cong \ \Gamma(\wfrakg \times \wfrakg, \, \calO_{\Delta \wfrakg}(\lambda,\mu)), 
\]
compatible with the natural action of $\mathbf{\Gamma}(\wfrakg \times \wfrakg)$.

The object $\calK_{s,t}$ is the kernel associated with the functor
\[
F_{s,t}:=F_{\wfrakg}^{\calO_{Z_s}} \circ F_{\wfrakg}^{\calO_{Z_t}} \circ F_{\wfrakg}^{\calO_{Z_s}} \circ (F_{\wfrakg}^{\calO_{Z_t}})^{-1} \circ (F_{\wfrakg}^{\calO_{Z_s}})^{-1} \circ (F_{\wfrakg}^{\calO_{Z_t}})^{-1}.
\]
We have seen in Proposition \ref{prop:propfunctorsbraidrel} that $F_{s,t}$ fixes any line bundle $\calO_{\wfrakg}(\lambda)$ with $\lambda \in \bbX^+$. Moreover, for any $\lambda,\mu$ we have, by the projection formula, $R\Gamma(\wfrakg, \, F^{\calK_{s,t}}_{\wfrakg}(\calO_{\wfrakg}(\lambda)) \otimes_{\calO_{\wfrakg}} \calO_{\wfrakg}(\mu)) \cong R\Gamma(\wfrakg \times \wfrakg, \, \calK_{s,t}(\lambda, \mu))$. It follows, using \cite[Theorem III.5.2]{HARAG}, that for $\lambda,\mu \gg 0$ we have $R^{\neq 0}(\wfrakg \times \wfrakg, \, \calK_{s,t}(\lambda, \mu))=0$ and, moreover, there is an isomorphism 
\[
\Gamma(\wfrakg \times \wfrakg, \, \calK_{s,t}(\lambda, \mu)) \ \cong \ \Gamma(\wfrakg, \, \calO_{\wfrakg}(\lambda + \mu)) \ \cong \ \Gamma(\wfrakg \times \wfrakg, \, \calO_{\Delta \wfrakg}(\lambda,\mu)).
\]

It remains to show that these isomorphisms can be chosen in a canonical way, and are compatible with the action of $\mathbf{\Gamma}(\wfrakg \times \wfrakg)$. We claim that the restriction morphisms induced by $j^*$:
\begin{align*}
\Gamma(\wfrakg \times \wfrakg, \, \calK_{s,t}(\lambda, \mu)) \ & \to \ \Gamma(\wfrakg_{\reg} \times \wfrakg_{\reg}, \, j^* \calK_{s,t}(\lambda, \mu)), \\
\Gamma(\wfrakg \times \wfrakg, \, \calO_{\Delta \wfrakg}(\lambda,\mu)) \ & \to \ \Gamma(\wfrakg_{\reg} \times \wfrakg_{\reg}, \, \calO_{\Delta \wfrakg_{\reg}}(\lambda,\mu))
\end{align*}
are isomorphisms. Indeed, the first morphism coincides (via the projection formula) with the restriction morphism
\[
\Gamma(\wfrakg, \, F^{\calK_{s,t}}_{\wfrakg}(\calO_{\wfrakg}(\lambda)) \otimes_{\calO_{\wfrakg}} \calO_{\wfrakg}(\mu)) \to \Gamma(\wfrakg_{\reg}, \, F^{j^* \calK_{s,t}}_{\wfrakg_{\reg}}(\calO_{\wfrakg_{\reg}}(\lambda)) \otimes_{\calO_{\wfrakg_{\reg}}} \calO_{\wfrakg_{\reg}}(\mu))
\]
induced by the inverse image under the inclusion $\wfrakg_{\reg} \hookrightarrow \wfrakg$. As the sheaf $F^{\calK_{s,t}}_{\wfrakg}(\calO_{\wfrakg}(\lambda)) \otimes_{\calO_{\wfrakg}} \calO_{\wfrakg}(\mu)$ is a line bundle, and the complement of $\wfrakg_{\reg}$ in $\wfrakg$ has codimension at least $2$ (see Proposition \ref{prop:codimreg}), the latter morphism is an isomorphism. The arguments for the second morphism are similar. 

It follows that there exists a unique isomorphism
\[
\Gamma(\wfrakg \times \wfrakg, \, \calK_{s,t}(\lambda, \mu)) \ \cong \ \Gamma(\wfrakg \times \wfrakg, \, \calO_{\Delta \wfrakg}(\lambda,\mu))
\]
which is compatible with our canonical isomorphism \eqref{eq:reg-part}. With this choice, the compatibility with the action of the algebra
\[
\mathbf{\Gamma}(\wfrakg \times \wfrakg) \ \cong \ \bigoplus_{\lambda,\mu \in \bbX^+} \ \Gamma(\wfrakg_{\reg} \times \wfrakg_{\reg}, \, \calO_{\wfrakg_{\reg} \times \wfrakg_{\reg}}(\lambda,\mu))
\]
is clear. \end{proof}

This corollary finishes the proof of Theorem \ref{thm:existenceaction}.

\section{Description of the kernels}
\label{sec:kernels}

\subsection{Statement}
\label{ss:statement-kernels-R}

Let $R=\Z[\frac{1}{h!}]$, where $h$ is the Coxeter number of $G_{\Z}$. For any $w \in W$, we define
\[
Z_{w,R}
\]
as the closure in $\wfrakg_R \times_R \wfrakg_R$ of the inverse image under the morphism $\wfrakg_R \times_{\frakg^*_R} \wfrakg_R \hookrightarrow \wfrakg_R \times_R \wfrakg_R \to \calB_R \times_R \calB_R$ of the $G_R$-orbit of $(B_R/B_R, w^{-1}B_R/B_R) \in \calB_R \times_R \calB_R$ (for the diagonal action). It is a reduced closed subscheme of $\wfrakg_R \times_{\frakg^*_R} \wfrakg_R$.

We also set
\[
Z_{w,R}' \ := \ Z_{w,R} \cap (\wcalN_R \times_R \wfrakg_R).
\]
(Note that here we take the \emph{scheme-theoretic} intersection.) It is easy to see that $Z_{w,R}'$ is in fact a closed subscheme of $\wcalN_R \times_R \wcalN_R$. (See Lemma \ref{lem:Z_w'} below in the case of a field.)

The main result of this section is the following.

\begin{thm}
\label{thm:descriptionkernels-R}

Let $R=\Z[\frac{1}{h!}]$.

Let $w \in W$, and let $w=s_1 \cdots s_n$ be a reduced expression (where $s_i \in \scS$). There exists an isomorphism in $\calD^b \Coh^{G_R \times_R (\Gm)_R}(\wfrakg_R \times_R \wfrakg_R)$, respectively in $\calD^b \Coh^{G_R \times_R (\Gm)_R}(\wcalN_R \times_R \wcalN_R)$:
\begin{eqnarray*}
& \calO_{Z_{s_1,R}} \star \cdots \star \calO_{Z_{s_n,R}} \ \cong \ \calO_{Z_{w,R}}, & \\
& \text{respectively } \ \calO_{Z_{s_1,R}'} \star \cdots \star \calO_{Z_{s_n,R}'} \ \cong \ \calO_{Z_{w,R}'}. &
\end{eqnarray*}

\end{thm}

Let again $w \in W$, and let $w=s_1 \cdots s_n$ be a reduced expression (where $s_i \in \scS$). Then by definition $T_w=T_{s_1} \cdots T_{s_n}$. Theorem \ref{thm:descriptionkernels-R} allows to give an explicit description of the action of $T_w$ on $\calD^b \Coh(\wfrakg_R)$ or $\calD^b \Coh(\wcalN_R)$ constructed in Theorem \ref{thm:existenceaction-R}. Namely, there exist isomorphisms of functors
\[
\mathbf{J}_{T_w} \ \cong \ F_{\wfrakg_R}^{\calO_{Z_{w,R}}}, \quad \mathbf{J}_{T_w}' \ \cong \ F_{\wcalN_R}^{\calO_{Z_{w,R}'}}.
\]

As for Theorem \ref{thm:existenceaction-R}, the proof of Theorem \ref{thm:descriptionkernels-R} is based on the reduction to the case of an algebraically closed field of positive characteristic, which will be treated using Representation Theory. We first treat this case. The proof of Theorem \ref{thm:descriptionkernels-R} for $\wfrakg_R$ is given in \S \ref{ss:proof-thm-description-kernels}. The case of $\wcalN_R$ is treated in \S \ref{ss:action-wcalN}.

\subsection{Kernels for the finite braid group: case of a field} \label{ss:kernels-statement}

From now on and until \S \ref{ss:proof-thm-description-kernels}, we fix an algebraically closed field of characteristic $p>h$. We use the same notation as in Section \ref{sec:braidrelations}, and we drop the index ``$\bk$'' for simplicity.

For $w \in W$, we denote by $\frakX_w^0 \subset \calB \times \calB$ the $G$-orbit of $(B/B,wB/B)$ for the diagonal action (a \emph{Schubert cell}), and by $\frakX_w$ its closure (a \emph{Schubert variety}). We denote by
\[
Z_w
\]
the closure of the inverse image of $\frakX_{w^{-1}}^0$ under the morphism $\wfrakg \times_{\frakg^*} \wfrakg \hookrightarrow \wfrakg \times \wfrakg \twoheadrightarrow \calB \times \calB$. This is a reduced closed subscheme of $\wfrakg \times_{\frakg^*} \wfrakg$. Note that it is not clear at this point that $Z_w$ is isomorphic to $Z_{w,R} \times_{\mathrm{Spec}(R)} \mathrm{Spec}(\bk)$ (because the latter scheme is a priori not reduced). This will follow from our results (see Remark \ref{rk:reduction-Z_w} below). The fiber of $Z_w$ over $(B/B,w^{-1}B/B)$ is 
\[
\bigl( \frakg / (\frakn + w^{-1} \cdot \frakn) \bigr)^*.
\]
In particular, by $G$-equivariance the restriction of $Z_w$ to the inverse image of $\frakX_{w^{-1}}^0$ is a vector bundle over $\frakX_{w^{-1}}^0$, of rank $\dim(\frakb) - \ell(w)$.

For $w \in W$, we define
\[
Z_w^{\reg} \subset \wfrakg_{\reg} \times \wfrakg_{\reg}
\]
to be the graph of the action of $w$ provided by Proposition \ref{prop:actionreg}. Then one easily checks that $Z_w^{\reg}=Z_w \cap (\wfrakg_{\reg} \times \wfrakg_{\reg})$ and that $Z_w$ is the closure of $Z_w^{\reg}$.

We also set 
\[
Z_w' \ := \ Z_w \cap (\wcalN \times \wfrakg).
\]
It is easy to check that $Z_w'$ is in fact a closed subscheme of $\wcalN \times \wcalN$ (see Lemma \ref{lem:Z_w'} below).

The version over $\bk$ of Theorem \ref{thm:descriptionkernels-R} is the following.

\begin{thm} \label{thm:descriptionkernels}

Assume $p>h$.

Let $w \in W$, and let $w=s_1 \cdots s_n$ be a reduced expression ($s_i \in \scS$). There is an isomorphism in $\calD^b \Coh^{G \times \Gm}(\wfrakg \times \wfrakg)$, respectively in $\calD^b \Coh^{G \times \Gm}(\wcalN \times \wcalN)$:
\begin{eqnarray*}
& \calO_{Z_{s_1}} \star \cdots \star \calO_{Z_{s_n}} \ \cong \ \calO_{Z_w}, & \\
& \text{respectively } \ \calO_{Z_{s_1}'} \star \cdots \star \calO_{Z_{s_n}'} \ \cong \ \calO_{Z_w'}.& 
\end{eqnarray*}
Moreover, $Z_w$, respectively $Z_w'$, is Cohen--Macaulay with dualizing sheaf 
\begin{eqnarray*}
& \calO_{Z_{s_n}}(-\rho, \rho-\alpha_n) \star \cdots \star \calO_{Z_{s_1}}(-\rho, \rho-\alpha_1), & \\
& \text{respectively } \
\calO_{Z_{s_n}'}(-\rho, \rho-\alpha_n) \star \cdots \star \calO_{Z_{s_1}'}(-\rho, \rho-\alpha_1), & 
\end{eqnarray*}
where $s_i=s_{\alpha_i}$ ($\alpha_i \in \Sigma$). In particular, these objects are also concentrated in degree $0$, i.e.~are coherent sheaves.

\end{thm}

\begin{remark} 
\label{rk:thm-description-kernels}
\begin{enumerate}
\item 
Probably, the restriction on the characteristic is not necessary. For instance, it follows from \cite{RAct} that the theorem is true for all $p$ if $w$ is an element of a parabolic subgroup of $W$ of type $\mbfA_1 \times \mbfA_1$ or $\mbfA_2$, and for $p \neq 2$ if $w$ is an element of a parabolic subgroup of $W$ of type $\mbfB_2$. To obtain a weaker restriction on $p$ in the general case, it would certainly be necessary to understand better the geometric properties of the varieties $Z_w$. For instance, it is proved in \cite{RAct} that $Z_w$ is normal if $p$ and $w$ are as above. We don't know if this property is true in general.
\item 
\label{it:rk-thm-description-kernels-char-0}
The case $p=0$ is excluded from Theorem \ref{thm:descriptionkernels}. However, it follows from Theorem \ref{thm:descriptionkernels-R} that Theorem \ref{thm:descriptionkernels} is also true in this case (by compatibility of convolution with change of scalars). One can also prove Theorem \ref{thm:descriptionkernels} for $\bk$ of characteristic zero directly, using Saito's theory of mixed Hodge modules. This will be the subject of a future publication.
\item As in the proof of Corollary \ref{cor:braidrelations}, once we have proved the theorem, one can check using restriction to $\wfrakg_{\reg}$ that the isomorphisms
\[
\calO_{Z_{s_1}} \star \cdots \star \calO_{Z_{s_n}} \ \cong \ \calO_{Z_w}, \quad \calO_{Z_{s_1}'} \star \cdots \star \calO_{Z_{s_n}'} \ \cong \ \calO_{Z_w'}
\]
can be chosen in a \emph{canonical} way.
\item
It follows from the theorem that for any $v,w \in W$ such that $\ell(vw)=\ell(v)+\ell(w)$, there exists an isomorphism
\[
\calO_{Z_v} \star \calO_{Z_w} \ \cong \ \calO_{Z_{vw}}.
\]
And, again by the same arguments as in the proof of Corollary \ref{cor:braidrelations}, this isomorphism can be chosen canonically. With this choice, one easily checks that the condition of \cite[Theorem 1.5]{De} holds. Hence the restriction of our action to the Artin braid group $\bB$ can be ``lifted'' to an action in the strong sense of \cite{De}. We don't know if this property holds for the whole of $\bB_{\aff}$.
\end{enumerate}
\end{remark}

We first concentrate on the case of $\wfrakg$. The proof of Theorem \ref{thm:descriptionkernels} in this case is given in \S \ref{ss:endproof}. The case of $\wcalN$ is treated in \S \ref{ss:action-wcalN}.

Let again $w \in W$, and let $w=s_1 \cdots s_n$ be a reduced expression ($s_i \in \scS$). For $i=1, \cdots, n$, let $\alpha_i \in \Sigma$ be the simple root attached to the simple reflection $s_i$. Then we define 
\begin{align*}
\calK_w \ & := \ \calO_{Z_{s_1}} \star \cdots \star \calO_{Z_{s_n}}, \\
\calK_w^{\dag} \ & := \ \calO_{Z_{s_n}}(-\rho, \rho-\alpha_n) \star \cdots \star \calO_{Z_{s_1}}(-\rho, \rho-\alpha_1),
\end{align*}
considered as objects in $\calD^b \Coh(\wfrakg \times \wfrakg)$. By Corollary \ref{cor:braidrelations}, these objects do not depend on the choice of the reduced expression (up to isomorphism). We will sometimes use the fact that $\calK_w$ has a canonical lift to an object of $\calD^b \Coh^{G \times \Gm}(\wfrakg \times \wfrakg)$. (This is also the case for $\calK_w^{\dag}$, but we will not use it.)

By definition we have isomorphisms of functors $\mathbf{J}_{T_w} \cong F^{\calK_w}_{\wfrakg}$ and $\mathbf{J}_{T_{w}^{-1}} \cong F^{\calK_w^{\dag}}_{\wfrakg}$. Moreover, we have
\begin{equation}
\label{eq:duality}
\calK_w^{\dag} \ \cong \ R\sheafHom_{\calO_{\wfrakg^2}}(\calK_w, \calO_{\wfrakg^2})[\dim(\frakg)]
\end{equation}
(see \cite[Proposition 5.9]{Huy} for details). Let $\delta$ be the automorphism of $\wfrakg \times \wfrakg$ which exchanges the two factors. Then we have
\begin{equation}
\label{eq:exchange}
\delta^* \calK_w \cong \calK_{w^{-1}}, \quad \delta^* \calK_w^{\dag} \cong \calK_{w^{-1}}^{\dag}.
\end{equation}

The proof of Theorem \ref{thm:descriptionkernels} is based on Representation Theory, and more precisely on localization theory for Lie algebras in positive characteristic, as studied in \cite{BMR, BMR2, BM}. Hence sometimes we will rather consider $\calK_w$ and $\calK_w^{\dag}$ as sheaves on $\wfrakg^{(1)} \times \wfrakg^{(1)}$, where ${}^{(1)}$ denotes the Frobenius twist, i.e.~we will twist the structures as sheaves of $\bk$-vector spaces. For simplicity we do not indicate this in the notation.

\subsection{More notation}
\label{ss:morenotation}

Recall that if $X$ is a scheme, and $Y \subset X$ a closed subscheme, one says that an
$\calO_X$-module $\calF$ is \emph{supported on} $Y$ if $\calF_{|X - Y}=0$. We write
$\Coh_Y(X)$ for the full subcategory of $\Coh(X)$ whose objects are supported on $Y$.

We denote by $\frakh$ the universal Cartan subalgebra of $\frakg$. It is isomorphic to $\frakb_0/[\frakb_0,\frakb_0]$ for any Lie algebra $\frakb_0$ of a Borel subgroup of $G$. In particular, $\frakh$ identifies naturally with $\frakt$ via the morphism $\frakt \xrightarrow{\sim} \frakb/\frakn \cong \frakh$.

For $s=s_{\alpha}$ ($\alpha \in \Sigma$), we denote by $\fraksl(2,s)$ the subalgebra of $\frakg$ generated by $h_{\alpha}$, $e_{\alpha}$ and $e_{-\alpha}$. It is isomorphic to $\fraksl(2,\bk)$.

The extended affine Weyl group $W_{\aff}$ naturally acts on $\bbX$. We will also consider the ``twisted'' action defined by
\[
w \bullet \lambda \ := \ w(\lambda+\rho) - \rho
\]
for $w \in W_{\aff}$, $\lambda \in \bbX$, where $\rho$ is the half sum of the positive roots.

Recall that an element $\chi \in \frakg^*$ is said to be \emph{nilpotent} if it is conjugate to an element of $(\frakg/\frakb)^*$.

Let $\frakZ$ be the center of $\calU \frakg$, the enveloping
algebra of $\frakg$. The subalgebra of $G$-invariants
$\frakZ_{{\rm HC}}:=(\calU \frakg)^G$ is central in $\calU \frakg$. This
is the ``Harish-Chandra part'' of the center, isomorphic
to $S(\frakt)^{(W,\bullet)}$, the algebra of $W$-invariants in the
symmetric algebra of $\frakt$, for the dot-action. The center
$\frakZ$ also has an other part, the ``Frobenius part''
$\frakZ_{{\rm Fr}}$, which is generated as an algebra by the elements
$X^p - X^{[p]}$ for $X \in \frakg$. It is isomorphic to
$S(\frakg^{(1)})$, the functions on the Frobenius twist
$\frakg^* {}^{(1)}$ of $\frakg^*$. Under our assumption $p>h$, there is an isomorphism (see \cite{KW} or \cite{MRCen}):
\[ \frakZ_{{\rm HC}} \otimes_{\frakZ_{{\rm Fr}} \cap \frakZ_{{\rm HC}}}
\frakZ_{{\rm Fr}}
\xrightarrow{\sim} \frakZ.\] Hence a character of $\frakZ$ is
given by a ``compatible pair'' $(\nu,\chi) \in \frakt^* \times
\frakg^{*}{}^{(1)}$. In this paper we only consider the case where
$\chi$ is nilpotent, and $\nu$ is \emph{integral}, i.e.~in the image
of the natural map $\bbX \to \frakt^*$. (Such a pair is always ``compatible.'') If $\lambda \in
\bbX$, we still denote by $\lambda$ its image in $\frakt^*$. If $\lambda \in \bbX$, we denote by
$\Mod^{\fg}_{(\lambda,\chi)}(\calU
\frakg)$ the category of finitely generated $\calU \frakg$-modules on
which $\frakZ$ acts with generalized character $(\lambda,\chi)$. Similarly, we denote by $\Mod^{\fg}_{\lambda}(\calU \frakg)$ the category of finitely generated $\calU \frakg$-modules on which $\frakZ_{\mathrm{HC}}$ acts with generalized character $\lambda$.

The translation functors for $\calU \frakg$-modules are defined in \cite[\S 6.1]{BMR}. More precisely, for $\lambda,\mu \in \bbX$ and $\chi \in \frakg^*{}^{(1)}$, the functor
\[
T_{\lambda}^{\mu} : \Mod^{\fg}_{(\lambda,\chi)}(\calU
\frakg) \to \Mod^{\fg}_{(\mu,\chi)}(\calU
\frakg)
\]
sends the module $M$ to
\[
\mathrm{pr}_{\mu} \bigl( M \otimes_{\bk} V(\mu - \lambda) \bigr),
\]
where $V(\mu - \lambda)$ is the standard (induced) $G$-module with extremal weight $\mu-\lambda$ (i.e.~with highest weight the dominant $W$-conjugate of $\mu-\lambda$), and $\mathrm{pr}_{\mu}$ is the functor which sends a locally finite $\frakZ_{\mathrm{HC}}$-module to its generalized eigenspace associated with $\mu$.

We will also consider the baby Verma modules, as defined e.g.~in \cite[\S 3.1.4]{BMR}. Namely, let $\frakb_0$ be the Lie algebra of a Borel subgroup of $G$. Let $\chi \in \frakg^*{}^{(1)}$ be such that $\chi_{|\frakb_0^{(1)}}=0$, and let $\lambda \in \bbX$. Then the associated baby Verma module $M_{\frakb_0,\chi;\lambda}$ is by definition the $\calU \frakg$-module
\[
M_{\frakb_0,\chi;\lambda} \ := \ (\calU \frakg)_{\chi} \otimes_{(\calU \frakb_0)_{\chi}} \bk_{\lambda}.
\]
Here $(\calU \frakg)_{\chi}:=(\calU \frakg) \otimes_{\frakZ_{\mathrm{Fr}}} \bk_{\chi}$, $(\calU \frakb_0)_{\chi}$ is the subalgebra generated by the image of $\frakb_0$, and $\lambda$ defines a character of $(\calU \frakb_0)_{\chi}$ via the morphism \[ \frakb_0 \to \frakb_0/[\frakb_0,\frakb_0] \cong \frakh \cong \frakt \xrightarrow{\lambda} \bk.\]

Finally, we will need the functors $\mathrm{Tens}_{\lambda}^{\mu}$ of \cite[\S 2.1.3]{BMR2}. Consider the sheaf of algebras $\wcalD$ on $\calB$ defined as follows (see \cite[\S 3.1.3]{BMR}). Let $p: G/U \to G/B$ be the natural quotient; it is a torsor for the abelian group $B/U$ (acting via multiplication on the right). Then
\[
\wcalD \ := \ p_*(\calD_{G/U})^{B/U},
\]
where $\calD_{G/U}$ is the sheaf of crystalline differential operators on $G/U$ (see \cite[\S 1.2]{BMR}). The sheaf of algebras $\wcalD$ can be naturally considered as a sheaf of algebras on the scheme $\wfrakg^{(1)} \times_{\frakh^*{}^{(1)}} \frakh^*$, where the morphism $\wfrakg \to \frakh^*$ sends $(X,gB)$ to $X_{|g \cdot \frakb} \in (g \cdot \frakb / [g \cdot \frakb, g \cdot \frakb])^* \cong \frakh^*$, and the morphism $\frakh^* \to \frakh^*{}^{(1)}$ is the Artin--Schreier map defined by the algebra morphism \[ \left\{ \begin{array}{ccc} \mathrm{S}(\frakh^{(1)}) & \to & \mathrm{S}(\frakh) \\ X \in \frakh^{(1)} & \mapsto & X^p - X^{[p]} \end{array} . \right. \] Moreover, it is an Azumaya algebra on this scheme. 

For $\chi \in \frakg^*{}^{(1)}$, we denote by $\calB_{\chi}^{(1)}$ the inverse image of $\chi$ under the natural morphism $\wfrakg^{(1)} \to \frakg^*{}^{(1)}$, endowed with the reduced subscheme structure. We denote by
$\Mod^{{\rm c}}(\wcalD)$ the
category of quasi-coherent, locally finitely generated
$\wcalD$-modules (equivalently, either on $\calB$ or on $\wfrakg^{(1)}
\times_{\frakh^*{}^{(1)}} \frakh^*$). For $\nu \in \bbX$ (considered as a linear form on $\frakt \cong \frakh$) and $\chi \in \frakg^*{}^{(1)}$ we denote by $\Mod^{{\rm
    c}}_{\nu}(\wcalD)$, respectively $\Mod^{{\rm c}}_{(\nu,\chi)}(\wcalD)$, the
full subcategory of $\Mod^{{\rm c}}(\wcalD)$ whose objects are
supported on $\wcalN^{(1)} \times \{ \nu \} \subset
\wfrakg^{(1)} \times_{\frakh^*{}^{(1)}} \frakh^*$, respectively on
$\calB_{\chi}^{(1)} \times \{ \nu \} \subset \wfrakg^{(1)}
\times_{\frakh^*{}^{(1)}} \frakh^*$. Then for $\lambda,\mu \in \bbX$ we consider the equivalence of categories \[ \mathrm{Tens}_{\lambda}^{\mu} : \left\{ \begin{array}{ccc} \Mod^{{\rm c}}_{(\lambda,\chi)}(\wcalD) & \to & \Mod^{{\rm c}}_{(\mu,\chi)}(\wcalD) \\ \calF & \mapsto & \calO_{\calB}(\mu - \lambda) \otimes_{\calO_{\calB}} \calF
\end{array} . \right. \]

\subsection{Reminder on localization in positive characteristic} \label{ss:reminderlocalization} 

Recall that a weight $\lambda \in \bbX$ is called
\emph{regular} if, for any root $\alpha$, $\langle \lambda +
\rho, \alpha^{\vee} \rangle \notin p \mathbb{Z}$, i.e.~if $\lambda$ is not on any reflection hyperplane of $W_{\aff}^{\mathrm{Cox}}$ (for the dot-action). Under our assumption $p>h$, $0 \in \bbX$ is regular; in particular, regular weights exist. By \cite[Theorem 5.3.1]{BMR} we have:

\begin{thm}\label{thm:thmBMR}

Let $\lambda \in \bbX$ be regular, and let $\chi \in \frakg^* {}^{(1)}$ be nilpotent. There exists an equivalence 
\[
\calD^b
  \Coh_{\calB_{\chi}^{(1)}}(\wfrakg^{(1)})
  \ \cong \ \calD^b\Mod^{\fg}_{(\lambda,\chi)}(\calU \frakg)
\]
of triangulated categories.\qed

\end{thm}

Let us recall briefly how this equivalence can be
constructed. We use
the notation of \cite{BMR}. Consider the sheaf of
algebras $\wcalD$ (see \S \ref{ss:morenotation}). If $\lambda \in \bbX$ is
regular, the global sections functor $R\Gamma_{\lambda} : \calD^b \Mod^{{\rm c}}_{\lambda}(\wcalD)
\to \calD^b \Mod^{\fg}_{\lambda}(\calU \frakg)$ is an equivalence of
categories. Its inverse is the localization functor
$\calL^{\hat{\lambda}}$. These functors restrict to equivalences
between the categories $\calD^b \Mod^c_{(\lambda,\chi)}(\wcalD)$ and
$\calD^b\Mod^{\fg}_{(\lambda,\chi)}(\calU \frakg)$.

Next, the Azumaya algebra $\wcalD$ splits on the formal neighborhood
of $\calB_{\chi}^{(1)} \times \{ \lambda \}$ in $\wfrakg^{(1)}
\times_{\frakh^*{}^{(1)}} \frakh^*$. Hence, the choice of a splitting
bundle on this formal neighborhood yields an equivalence of
categories $\Coh_{\calB_{\chi}^{(1)} \times \{ \lambda \} }(\wfrakg^{(1)}
\times_{\frakh^*{}^{(1)}} \frakh^*) \cong
\Mod^c_{(\lambda,\chi)}(\wcalD)$. Finally, the projection $\wfrakg^{(1)}
\times_{\frakh^{*}{}^{(1)}} \frakh^* \to
\wfrakg^{(1)}$ induces an isomorphism between the formal
neighborhood of $\calB_{\chi}^{(1)} \times \{ \lambda \}$
and the formal neighborhood of $\calB_{\chi}^{(1)}$ (see \cite[\S 1.5.3.c]{BMR2}). This isomorphism induces
an equivalence $\Coh_{\calB_{\chi}^{(1)} \times
  \{ \lambda \} }(\wfrakg^{(1)} \times_{\frakh^*{}^{(1)}} \frakh^*) \cong
\Coh_{\calB^{(1)}_{\chi}}(\wfrakg^{(1)})$, which finally gives the equivalence of Theorem \ref{thm:thmBMR}.

We choose the normalization of the splitting bundles as in
\cite[\S 1.3.5]{BMR2}. We denote by $\calM_{(\lambda,\chi)}$
the splitting bundle associated to $\lambda$, and by \[\gamma_{(\lambda,\chi)} :
\calD^b\Coh_{\calB_{\chi}^{(1)}}(\wfrakg^{(1)}) \xrightarrow{\sim}
\calD^b\Mod^{\fg}_{(\lambda,\chi)}(\calU \frakg)\] the associated equivalence.

\subsection{Reminder and corrections of \cite{BMR2} and \cite{RAct}}

First, let us recall the representation-theoretic interpretation of the braid group action in positive characteristic.

Let us fix a character $\lambda \in \bbX$ in the fundamental alcove (i.e.~such that for all $\alpha \in \Phi^+$ we have $0 < \langle \lambda + \rho, \alpha^{\vee} \rangle < p$), and some $\chi \in \frakg^*{}^{(1)}$ nilpotent. There is a natural \emph{right} action of the group $W_{\aff}$ on the set $W_{\aff} \bullet \lambda$, defined by 
\[
(w \bullet \lambda) * v := wv \bullet \lambda
\]
for $v,w \in W_{\aff}$. 

For $\alpha \in \Sigma_{\aff}$, let $\mu_{\alpha} \in\bbX$ be a weight on the $\alpha$-wall of the fundamental alcove, and on no other wall. Then, for $s=s_{\alpha}$ we define the functor
\[
\mathbf{I}_{T_{s}} := \mathrm{Cone} \bigl( \Id \to T_{\mu_{\alpha}}^{\lambda} \circ T_{\lambda}^{\mu_{\alpha}} \bigr) : \ \calD^b \Mod^{\fg}_{(\lambda,\chi)}(\calU \frakg) \to \calD^b \Mod^{\fg}_{(\lambda,\chi)}(\calU \frakg).
\]
(Here the morphism of functors is induced by adjunction.) This functor is well defined,\footnote{Recall that if $\mathscr{A}$, $\mathscr{B}$ are abelian categories, $F,G: \mathscr{A} \to \mathscr{B}$ are exact functors and $\phi: F \to G$ is a morphism of functors, the functor
\[
\left\{
\begin{array}{ccc}
\calC^b(\mathscr{A}) & \to & \calC^b(\mathscr{B}) \\
X & \mapsto & \mathrm{Cone}(FX \xrightarrow{\phi(X)} GX)
\end{array}
\right.
\]
descends to derived categories. The resulting functors is denoted by $\mathrm{Cone}(F \to G)$.} and does not depend on the choice of $\mu_{\alpha}$, see \cite[Corollary 2.2.7]{BMR2}. Also, for $\omega \in \Omega$, we put
\[
\mathbf{I}_{T_{\omega}} \ := \ T_{\lambda}^{\omega \bullet \lambda} : \ \calD^b \Mod^{\fg}_{(\lambda,\chi)}(\calU \frakg) \to \calD^b \Mod^{\fg}_{(\lambda,\chi)}(\calU \frakg).
\]
(Observe that $\lambda$ and $\omega \bullet \lambda$ have the same image in $\frakt^*/(W,\bullet)$.)

Let $w \in W_{\aff}$, and $\mu \in W_{\aff} \bullet \lambda$. Write $w=s_1 \cdots s_n \omega$, where $s_i \in \scS_{\aff}$, $\omega \in \Omega$, and $\ell(w)=n$. Recall (\cite[\S 2.1.3]{BMR2}) that one says that $w$ \emph{increases} $\mu$ if for each $i=1, \cdots, n$, $\mu * (s_1 \cdots s_{i-1}) < \mu * (s_1 \cdots s_i)$ for the standard order on $\bbX$.

The following theorem is a corrected version of \cite[Theorem 2.1.4]{BMR2}. The same proof applies.

\begin{thm} \label{thm:intertwiner}

\begin{enumerate}

\item 
The assignment $T_{s} \mapsto \mathbf{I}_{T_s}$ ($s \in \scS_{\aff}$), $T_{\omega} \mapsto \mathbf{I}_{T_{\omega}}$ extends to a \emph{right} action of the group $\bB_{\aff}$ on $\calD^b \Mod^{\fg}_{(\lambda,\chi)}(\calU \frakg)$.
\item 
\label{it:thm-intertwiner-2}
For $b \in \bB_{\aff}$, we denote by
\[
\mathbf{I}_{b} : \ \calD^b \Mod^{\fg}_{(\lambda,\chi)}(\calU \frakg) \to \calD^b \Mod^{\fg}_{(\lambda,\chi)}(\calU \frakg)
\]
the action of $b$. Then for $\mu \in W_{\aff} \bullet \lambda$ and $w \in W_{\aff}$ such that $w$ increases $\mu$, there is a natural isomorphism of functors
\[
\calL^{\widehat{\mu * w}} \circ \mathbf{I}_{T_w} \ \cong \ \mathrm{Tens}_{\mu}^{\mu * w} \circ \calL^{\widehat{\mu}}.
\]

\end{enumerate}

\end{thm}

Let us recall that property \eqref{it:thm-intertwiner-2} implies the following description of the action of (lifts of) dominants weights (see \cite[Proposition 2.3.3]{BMR2}).

\begin{cor} \label{cor:dominant}

For $\nu \in \bbX$ dominant and for $\calF$ in $\calD^b \Coh_{\calB_{\chi}^{(1)}}(\wfrakg^{(1)})$, there is a isomorphism \[ \mathbf{I}_{T_{\nu}} \circ \gamma_{(\lambda,\chi)}(\calF) \ \cong \ \gamma_{(\lambda,\chi)}(\calF \otimes_{\calO_{\wfrakg^{(1)}}} \calO_{\wfrakg^{(1)}}(\nu)) \]
which is functorial in $\calF$.\qed

\end{cor}

Now we consider the Bernstein presentation of $\bB_{\aff}$ (see \S \ref{ss:notation}). The relations in this presentation are symmetric. Hence there is a natural anti-automorphism $\iota : \bB_{\aff} \xrightarrow{\sim} \bB_{\aff}$ which is the identity on the generators $T_s$ ($s \in \scS$) and $\theta_x$ ($x \in \bbX$). The following theorem is a corrected version of \cite[Theorem 5.4.1]{RAct}. (The same proof as in \cite{RAct} works, taking into account the corrections to \cite{BMR2} given above.) It provides a link between the \emph{geometric} braid group action of Theorem \ref{thm:existenceaction} and the \emph{representation-theoretic} braid group action constructed in \cite{BMR2} and recalled in Theorem \ref{thm:intertwiner}.

\begin{thm} \label{thm:RTGeom}

For any $b \in \bB_{\aff}$ the following diagram:
\[
\xymatrix{ \calD^b \Coh_{\calB_{\chi}^{(1)}}(\wfrakg^{(1)}) \ar[d]_{\gamma_{(\lambda,\chi)}}^{\wr} \ar[rr]^{\mathbf{J}_b} & & \calD^b \Coh_{\calB_{\chi}^{(1)}}(\wfrakg^{(1)}) \ar[d]^{\gamma_{(\lambda,\chi)}}_{\wr} \\ \calD^b \Mod^{\fg}_{(\lambda,\chi)}(\calU \frakg) \ar[rr]^{\mathbf{I}_{\iota(b)}} & & \calD^b \Mod^{\fg}_{(\lambda,\chi)}(\calU \frakg) }
\]
is commutative up to an isomorphism of functors.\qed

\end{thm}

\begin{remark} \label{rk:formal} The same formalism also applies to formal completions. For simplicity, here we assume that $\chi=0$. For $\lambda \in \bbX$ regular, by \cite[Theorem 5.4.1]{BMR} there is an equivalence of categories \[ \widehat{\gamma}_{\lambda} : \calD^b \Coh(\widehat{\calB^{(1)}}) \xrightarrow{\sim} \calD^b \Mod^{\fg}((\calU \frakg)^{\widehat{\lambda}}) \] between the derived category $\calD^b \Mod^{\fg}((\calU \frakg)^{\widehat{\lambda}})$ of finitely generated modules over the completion $(\calU \frakg)^{\widehat{\lambda}}$ of $\calU \frakg$ at the central character $(\lambda,0)$ and the derived category $\calD^b \Coh(\widehat{\calB^{(1)}})$ of coherent sheaves on the formal neighborhood $\widehat{\calB^{(1)}}$ of the zero section $\calB^{(1)}$ in $\wfrakg^{(1)}$. There exist a geometric action of $\bB_{\aff}$ on $\calD^b \Coh(\widehat{\calB^{(1)}})$ and a representation-theoretic action on $\calD^b \Mod^{\fg}((\calU \frakg)^{\widehat{\lambda}})$, and a statement similar to Theorem \ref{thm:RTGeom} also holds in this context.\end{remark}

\subsection{The kernels are sheaves}

From now on and until \S \ref{ss:irreducibility}, we fix some $w \in W$.

For $\lambda \in \bbX$ regular, recall the splitting bundle $\calM_{(\lambda,0)}$ for the Azumaya algebra $\wcalD$ on the formal neighborhood of $\calB^{(1)} \times \{\lambda\}$ in $\wfrakg^{(1)} \times_{\frakh^* {}^{(1)}} \frakh^*$ (see \S \ref{ss:reminderlocalization}). Here, for simplicity, we write $\calM_{\lambda}$ for $\calM_{(\lambda,0)}$. We will rather consider $\calM_{\lambda}$ as a vector bundle on the formal neighborhood $\widehat{\calB^{(1)}}$ of $\calB^{(1)}$ in $\wfrakg^{(1)}$, see \S \ref{ss:reminderlocalization}. Recall that there exists a $\Gm$-equivariant vector bundle $\mathbf{M}_{\lambda}$ on $\wfrakg^{(1)}$ such that $\calM_{\lambda}$ is the restriction of $\mathbf{M}_{\lambda}$ to $\widehat{\calB^{(1)}}$, i.e.~the completion of $\mathbf{M}_{\lambda}$ for the $\calI$-adic topology, where $\calI$ is the sheaf of ideals of the closed subscheme $\calB^{(1)} \subset \wfrakg^{(1)}$ (see \cite[\S 9.4]{R2}).

Recall also that the restriction $\wcalD^{\widehat{\lambda}}$ of the sheaf of algebras $\wcalD$ to the formal neighborhood of $\calB^{(1)} \times \{\lambda\}$ in $\wfrakg^{(1)} \times_{\frakh^* {}^{(1)}} \frakh^*$, identified with $\widehat{\calB^{(1)}}$, is isomorphic to $\sheafEnd_{\calO_{\widehat{\calB^{(1)}}}}(\calM_{\lambda})$. Hence  $(\wcalD^{\widehat{\lambda}})^{\mathrm{opp}}$ is isomorphic to $\sheafEnd_{\calO_{\widehat{\calB^{(1)}}}}(\calM_{\lambda}^{\vee})$, where $\calM_{\lambda}^{\vee}$ is the vector bundle dual to $\calM_{\lambda}$. Observe finally that $\calM_{\lambda}^{\vee}$ is the restriction of $\mathbf{M}_{\lambda}^{\vee}:=\sheafHom_{\calO_{\wfrakg^{(1)}}}(\mathbf{M}_{\lambda}, \calO_{\wfrakg^{(1)}})$ to $\widehat{\calB^{(1)}}$.

Let $p_1,p_2 : \wfrakg^{(1)} \times \wfrakg^{(1)} \to \wfrakg^{(1)}$ be the projections. Then \begin{equation} \label{eq:bimodule} p_1^*(\calM_0) \otimes_{\calO_{\widehat{\calB^{(1)}}^2}} \calK_w^{\dag} \otimes_{\calO_{\widehat{\calB^{(1)}}^2}} p_2^*(\calM_{w^{-1} \bullet 0}^{\vee}) \end{equation} is naturally a $\wcalD^{\widehat{0}} \boxtimes (\wcalD^{\widehat{w^{-1} \bullet 0}})^{\mathrm{opp}}$-module. (Here, the tensor products can be replaced by \emph{derived} tensor products without any change.) As a sheaf, it is the restriction of \[ p_1^*(\mathbf{M}_0) \otimes_{\calO_{(\wfrakg^{(1)})^2}} \calK_w^{\dag} \otimes_{\calO_{(\wfrakg^{(1)})^2}} p_2^*(\mathbf{M}_{w^{-1} \bullet 0}^{\vee}) \] to the formal neighborhood of the zero section in $\wfrakg^{(1)} \times \wfrakg^{(1)}$. Taking the derived global sections of \eqref{eq:bimodule}, we obtain an object in the derived category of $(\calU \frakg^{\widehat{0}}, \ \calU \frakg^{\widehat{w^{-1} \bullet 0}})$-bimodules. Moreover, here we have $\calU \frakg^{\widehat{w^{-1} \bullet 0}} = \calU \frakg^{\widehat{0}}$. Our first observation is the following.

\begin{lem} \label{lem:K_w^0}
 
There is an isomorphism in the derived category of $\calU \frakg^{\widehat{0}}$-bimodules \[ R\Gamma \bigl( p_1^*(\calM_0) \otimes_{\calO_{\widehat{\calB^{(1)}}^2}} \calK_w^{\dag} \otimes_{\calO_{\widehat{\calB^{(1)}}^2}} p_2^*(\calM_{w^{-1} \bullet 0}^{\vee}) \bigr) \ \cong \ (\calU \frakg)^{\widehat{0}}. \]

\end{lem}

\begin{proof} As above, let $\calI$ be the sheaf of ideals of the closed subscheme $\calB^{(1)} \subset \wfrakg^{(1)}$, and let $n \in \mathbb{Z}_{>0}$. By the projection formula, and then the definition of Fourier--Mukai transform we have
\begin{multline*}
R\Gamma(p_1^*(\calM_0) \, \lotimes_{\calO_{\widehat{\calB^{(1)}}^2}} \, \calK_w^{\dag} \, \lotimes_{\calO_{\widehat{\calB^{(1)}}^2}} \, p_2^*(\calM_{w^{-1} \bullet 0}^{\vee} / \calI^n \calM_{w^{-1} \bullet 0}^{\vee})) \\
\cong \ R\Gamma_0( \calM_0 \lotimes_{\calO_{\wfrakg^{(1)}}} \bigl( R(p_1)_* (\calK_w^{\dag} \, \lotimes_{\calO_{(\wfrakg^{(1)})^2}} \, p_2^*(\calM_{w^{-1} \bullet 0}^{\vee} / \calI^n \calM_{w^{-1} \bullet 0}^{\vee})) \bigr) ) \\
\cong R\Gamma_0 \bigl( \calM_0 \lotimes_{\calO_{\wfrakg^{(1)}}} F^{\delta^* \calK^{\dag}_w}_{\wfrakg^{(1)}}(\calM_{w^{-1} \bullet 0}^{\vee} / \calI^n \calM_{w^{-1} \bullet 0}^{\vee}) \bigr).
\end{multline*}
Now by isomorphisms \eqref{eq:exchange} and Theorem \ref{thm:RTGeom} the latter object is isomorphic to
\[
\mathbf{I}_{T_{w}^{-1}} \circ R\Gamma_0(\calM^0 \otimes_{\calO_{\wfrakg^{(1)}}} (\calM_{w^{-1} \bullet 0}^{\vee} / \calI^n \calM_{w^{-1} \bullet 0}^{\vee}) ).
\]
(Observe that $\iota(T_{w^{-1}}^{-1})=T_w^{-1}$.) The element $w$ increases the weight $w^{-1} \bullet 0$, and $(w^{-1} \bullet 0) * w = 0$. Hence by Theorem \ref{thm:intertwiner}\eqref{it:thm-intertwiner-2} we have an isomorphism of functors
\begin{equation}
\label{eq:intertwiner}
\calL^{\widehat{0}} \circ \mathbf{I}_{T_{w}} \ \cong \ \mathrm{Tens}_{w^{-1} \bullet 0}^0 \circ \calL^{\widehat{w^{-1} \bullet 0}}.
\end{equation}
Taking inverses on both sides, we obtain
\[
\mathbf{I}_{T_{w}^{-1}} \circ R\Gamma_0 \ \cong \ R\Gamma_{w^{-1} \bullet 0} \circ \mathrm{Tens}^{w^{-1} \bullet 0}_0.
\]
Moreover, by definition $\calO_{\calB}(w^{-1} \bullet 0) \otimes_{\calO_{\calB}} \calM_0 \cong \calM_{w^{-1} \bullet 0}$ (see \cite[\S 1.3.5]{BMR2}). Hence we obtain finally
\begin{multline*}
R\Gamma(p_1^*(\calM_0) \, \lotimes_{\calO_{\widehat{\calB^{(1)}}^2}} \, \calK_w^{\dag} \, \lotimes_{\calO_{\widehat{\calB^{(1)}}^2}} \, p_2^*(\calM_{w^{-1} \bullet 0}^{\vee} / \calI^n \calM_{w^{-1} \bullet 0}^{\vee})) \\
\cong \ R\Gamma_{w^{-1} \bullet 0} (\calM_{w^{-1} \bullet 0} \, \lotimes_{\calO_{\wfrakg^{(1)}}} \, (\calM_{w^{-1} \bullet 0}^{\vee} / \calI^n \calM_{w^{-1} \bullet 0}^{\vee})).
\end{multline*}
Then, taking the inverse limit over $n$ and using \cite[Theorem 4.1.5]{EGA3_1} applied to the projective morphisms $\wfrakg^{(1)} \to \frakg^* {}^{(1)}$ and $(\wfrakg^{(1)})^2 \to (\frakg^* {}^{(1)})^2$ we obtain
\begin{multline*}
R\Gamma(p_1^*(\calM_0) \otimes_{\calO_{\widehat{\calB^{(1)}}^2}} \calK_w^{\dag} \otimes_{\calO_{\widehat{\calB^{(1)}}^2}} p_2^*(\calM_{w^{-1} \bullet 0}^{\vee})) \\
\cong \ R\Gamma_{w^{-1} \bullet 0} (\calM_{w^{-1} \bullet 0} \otimes_{\calO_{\wfrakg^{(1)}}} \calM_{w^{-1} \bullet 0}^{\vee}).
\end{multline*}
The latter object is isomorphic to
\[
R\Gamma(\wcalD^{\widehat{w^{-1} \bullet 0}}) \ \cong \ \calU \frakg^{\widehat{w^{-1} \bullet 0}} \ = \ \calU \frakg^{\widehat{0}}.
\]
(For the first isomorphism we have used \cite[Proposition 3.4.1 and \S 5.4]{BMR}.) One easily checks that all our isomorphisms are compatible with the $\calU \frakg^{\widehat{0}}$-bimodule structures.
\end{proof}

\begin{cor} \label{cor:K_w^0}

The object $\calK_w^{\dag}$ has cohomology only in non-positive degrees.

\end{cor}

\begin{proof} It is sufficient to prove the same property for the object \[ p_1^*(\mathbf{M}_0) \otimes_{\calO_{(\wfrakg^{(1)})^2}} \calK_w^{\dag} \otimes_{\calO_{(\wfrakg^{(1)})^2}} p_2^*(\mathbf{M}_{w^{-1} \bullet 0}^{\vee}).\] Moreover, this object is in the essential image of the forgetful functor
\[
\calD^b \Coh^{\Gm}(\wfrakg^{(1)} \times \wfrakg^{(1)}) \to \calD^b \Coh(\wfrakg^{(1)} \times \wfrakg^{(1)}).
\]
Hence it is sufficient to show that its restriction to the formal neighborhood of the zero-section, namely
\[
p_1^*(\calM_0) \otimes_{\calO_{\widehat{\calB^{(1)}}^2}} \calK_w^{\dag} \otimes_{\calO_{\widehat{\calB^{(1)}}^2}} p_2^*(\calM_{w^{-1} \bullet 0}^{\vee}),
\]
has the same property. Now this fact follows directly from Lemma \ref{lem:K_w^0}, using the fact that the derived global sections functor $R\Gamma$ is an equivalence of categories between the derived categories of coherent $\wcalD^{\widehat{0}} \boxtimes (\wcalD^{\widehat{w^{-1} \bullet 0}})^{\mathrm{opp}}$-modules and of coherent $(\calU \frakg)^{\widehat{0}}$-bimodules, with inverse the localization functor (see \cite[Theorem 5.4.1]{BMR} and Remark \ref{rk:formal}), which is right-exact.\end{proof}

Now we consider the kernel $\calK_w$. We consider the couple $(-2\rho,0)$ as a Harish-Chandra character for the Lie algebra $\frakg \times \frakg$.

\begin{lem} \label{lem:K_w}

The object
\[
R\Gamma_{(-2\rho,0)}((\calM_{-2\rho} \boxtimes \calM_0) \otimes_{\wfrakg^{(1)} \times \wfrakg^{(1)}} \calK_w)
\]
is concentrated in degree $0$.

\end{lem}

\begin{proof} By the projection formula, and Theorem \ref{thm:RTGeom} (see also Remark \ref{rk:formal}) we have
\begin{align*} R\Gamma_{(-2\rho,0)}((\calM_{-2\rho} \boxtimes \calM_0) \otimes_{\wfrakg^{(1)} \times \wfrakg^{(1)}} \calK_w) \ & \cong \ R\Gamma_0(\calM_0 \otimes_{\wfrakg^{(1)}} F^{\calK_w}_{\wfrakg^{(1)}}(\calM_{-2\rho})) \\
& \cong \ \mathbf{I}_{T_{w^{-1}}} \circ R\Gamma_0(\calM_0 \otimes_{\wfrakg^{(1)}} \calM_{-2\rho}).
\end{align*}
(Here we use $\iota(T_w)=T_{w^{-1}}$.) Now by \cite[Lemma 3.0.6 and its proof]{BMR2}, we have an isomorphism $\calM_{-2\rho} \cong \calM_0^{\vee}$. Hence we obtain
\[
R\Gamma_{(-2\rho,0)}((\calM_{-2\rho} \boxtimes \calM_0) \otimes_{\wfrakg^{(1)} \times \wfrakg^{(1)}} \calK_w) \ \cong \ \mathbf{I}_{T_{w^{-1}}} \circ R\Gamma(\wcalD^{\widehat{0}}) \ \cong \ \mathbf{I}_{T_{w^{-1}}}(\calU \frakg^{\widehat{0}}).
\]
The object on the left hand side of these isomorphisms has, by definition, only cohomology in non-negative degrees. On the other hand, again by definition, the functor $\mathbf{I}_{T_{w^{-1}}}$ stabilizes the subcategory $\calD^{b,\leq 0} \Mod^{\fg}(\calU \frakg^{\widehat{0}})$. (It suffices to prove this statement when $\ell(w)=1$, in which case it is trivial.) Hence the object on the right hand side has cohomology only in non positive degrees. As these objects are isomorphic, they have to be concentrated in degree $0$.\end{proof}

A proof similar to that of Corollary \ref{cor:K_w^0} gives the following result.

\begin{cor} \label{cor:K_w}

The object $\calK_w$ has only cohomology in non-positive degrees.\qed

\end{cor}

We can can finally prove:

\begin{prop} \label{prop:K_wsheaf}

The objects $\calK_w$ and $\calK_w^{\dag}$ are concentrated in degree $0$.

\end{prop}

\begin{proof} Recall the following well-known fact, if $X$ is a smooth $\bk$-scheme and $\calF$ a coherent sheaf on $X$:
\[
\sheafExt^i_{\calO_X}(\calF,\calO_X)=0 \quad \text{for } \ i<\mathrm{codim}(\mathrm{Supp}(\calF)),
\leqno{(\#)}
\]
where $\mathrm{Supp}(\calF)$ is the support of $\calF$.

By Corollary \ref{cor:K_w}, $\calK_w$ has cohomology only in non-positive degrees. Moreover, by definition $\mathrm{Supp}(\calK_w)=Z_w$, hence $\mathrm{codim}(\mathrm{Supp}(\calK_w))=\dim(\frakg)$. Hence, using $(\#)$, it follows that \[ R\sheafHom_{\calO_{\wfrakg^{(1)} \times \wfrakg^{(1)}}}(\calK_w,\calO_{\wfrakg^{(1)} \times \wfrakg^{(1)}}) \] has cohomology only in degrees $\geq \dim(\frakg)$. Using \eqref{eq:duality} we deduce that $\calK_w^{\dag}$ has cohomology only in non-negative degrees. Then, by Corollary \ref{cor:K_w^0}, $\calK_w^{\dag}$ is concentrated in degree $0$.

The same arguments apply to $\calK_w$.\end{proof}

\subsection{Two preliminary lemmas}

Now we want to prove that the sheaf $\calK_w$ is generated by its global sections. For this we need two lemmas: one on the representation-theoretic side, and one on the geometric side.

Let us fix $\lambda \in \bbX$ in the fundamental alcove. Its differential induces a linear form on $\frakh$, via the natural isomorphism $\frakt \xrightarrow{\sim} \frakb/\frakn \cong \frakh$. Recall the definition of baby Verma modules, \S \ref{ss:morenotation}. As usual, we let $w_0 \in W$ be the longest element. We will use only the special case $\lambda=0$, $\chi=0$ of the next lemma (which we state in full generality for completeness).

\begin{lem} \label{lem:translationbabyVerma}

Let $\frakb_0 \subset \frakg$ be the Lie algebra of a Borel subgroup, and $\chi \in \frakg^*{}^{(1)}$ be such that $\chi_{|\frakb_0^{(1)}}=0$. For any nonzero submodule $N \subseteq M_{\frakb_0,\chi;w_0(\lambda)}$ we have $T_{\lambda}^{-\rho} N \neq 0$.

\end{lem}

\begin{proof} This claim is obviously independent of the choice of $\frakb_0$ (because the action of $G$ on Borel subgroups is transitive), hence we can assume $\frakb_0=\frakb$. It follows from \cite[B.3]{JANSub} that $M_{\frakb,\chi;w_0(\lambda)}$ is a submodule of $T_{-\rho}^{\lambda}(M_{\frakb,\chi;\rho})$. (Beware that Jantzen uses \emph{positive} Borel subalgebras to define baby Verma modules, while here $\frakb$ is the \emph{negative} Borel subalgebra.) Hence there is a nonzero morphism $N \to T_{-\rho}^{\lambda}(M_{\frakb,\chi;\rho})$. By adjunction, we have an isomorphism \[ \Hom_{\frakg}(N,T_{-\rho}^{\lambda}(M_{\frakb,\chi;\rho})) \ \cong \ \Hom_{\frakg}(T^{-\rho}_{\lambda}(N),M_{\frakb,\chi;\rho}). \] Hence the right hand side is nonzero, which implies $T_{\lambda}^{-\rho} N \neq 0$. \end{proof}

In the next lemma, we consider $\calK_w$ as a sheaf on $\wfrakg \times \wfrakg$. We let $p_2 : \wfrakg \times \wfrakg \to \wfrakg$ be the projection on the second component.

\begin{lem}
\label{lem:globalsections}
 
For any $w \in W$ we have an isomorphism
\[
R(p_2)_* \calK_w \ \cong \ \calO_{\wfrakg}
\]
in $\calD^b \Coh^{G \times \Gm}(\wfrakg)$. In particular, 
\[
R\Gamma(\wfrakg \times \wfrakg, \, \calK_w) \ \cong \ \Gamma(\wfrakg, \, \calO_{\wfrakg}) \ \cong \ \mathrm{S}(\frakg) \otimes_{\mathrm{S}(\frakh)^W} \mathrm{S}(\frakh)
\]
as $(G \times \Gm)$-equivariant algebras.

\end{lem}

\begin{proof} Let us begin with the first statement. By definition of convolution functors, $R(p_2)_* \calK_w \cong F_{\wfrakg}^{\calK_w}(\calO_{\wfrakg})$. By definition again, the functor $F_{\wfrakg}^{\calK_w}$ is the composition of functors of the form $F_{\wfrakg}^{\calO_{Z_s}}$ for $s \in \scS$. Hence it is sufficient to prove that $F_{\wfrakg}^{\calO_{Z_s}}(\calO_{\wfrakg}) \cong R(p_2)_* \calO_{Z_s}$ is isomorphic to $\calO_{\wfrakg}$. However, the morphism $p_2 : Z_s \to \wfrakg$ is proper, birational, and has normal image. Hence $(p_2)_* \calO_{Z_s} \cong \calO_{\wfrakg}$ by Zariski's Main Theorem. And the vanishing of $R^i (p_2)_* \calO_{Z_s}$ for $i>0$ can be proved exactly as in \cite[Proposition 3.4.1]{RAct}.

The second statement follows from the obvious isomorphism $R\Gamma(\wfrakg \times \wfrakg, -) \cong R\Gamma(\wfrakg, -) \circ R(p_2)_*$, and the isomorphisms
\[
R^i\Gamma(\wfrakg, \calO_{\wfrakg}) \cong \left\{
\begin{array}{cl}
0 & \text{if } i \neq 0 \\
\rmS(\frakg) \otimes_{\rmS(\frakh)^W} \rmS(\frakh) & \text{if } i=0
\end{array} \right.
\]
(see \cite[Proof of Proposition 3.4.1]{BMR} and references therein).
\end{proof}

\begin{remark} In other words, the first statement in Lemma \ref{lem:globalsections} is that we have an isomorphism \[ \mathbf{J}_{T_w}(\calO_{\wfrakg}) \ \cong \ \calO_{\wfrakg} \] for any $w \in W$. See also \cite[Lemma 1.3.4(c)]{BM}.
\end{remark}

\subsection{The sheaf $\calK_w$ is globally generated}

In the proof of the following proposition, we will consider localization functors for both groups $G$ and $G \times G$, and always for the case $\chi=0$. For simplicity and to avoid confusion, we write $\gamma^G_{\lambda}$, respectively $\gamma^{G \times G}_{(\lambda,\mu)}$, for $\gamma_{(\lambda,0)}$, respectively $\gamma_{((\lambda,0),(\mu,0))}$. In the second functor, $\lambda,\mu \in \bbX$, so that the pair $(\lambda,\mu)$ defines a character of the maximal torus $T \times T \subset G \times G$. We use similar notation for the functors $\widehat{\gamma}^G_{\lambda}$, $\widehat{\gamma}^{G \times G}_{(\lambda,\mu)}$ (see Remark \ref{rk:formal}).

\begin{prop} \label{prop:globalgeneration}
 
The sheaf $\calK_w$ is generated by its global sections.

\end{prop}

\begin{proof}
It is enough to prove the following property: $(\ddag)$ For any Borel subalgebras $\frakb_1,\frakb_2$ of $\frakg$, and any nonzero morphism of $\calO_{\wfrakg^{(1)} \times \wfrakg^{(1)}}$-modules $\calK_w \to \calO_{\{(0,\frakb_1), (0,\frakb_2)\}}$, the morphism
\[
\Gamma(\wfrakg^{(1)} \times \wfrakg^{(1)}, \, \calK_w) \ \to \ \bk
\]
obtained by taking global sections is also nonzero. 

Indeed, assume that property $(\ddag)$ is satisfied, and that $\calK_w$ is not globally generated. Let $\calF$ be the cokernel of the natural morphism 
\[
\Gamma(\wfrakg^{(1)} \times \wfrakg^{(1)}, \, \calK_w) \otimes_{\bk} \calO_{\wfrakg^{(1)} \times \wfrakg^{(1)}} \ \to \ \calK_w.
\]
By assumption, $\calF \neq 0$. The sheaf $\calK_w$ is $\Gm$-equivariant, hence the same is true for $\calF$. It follows that there exist Borel subalgebras $\frakb_1,\frakb_2$ of $\frakg$ such that the fiber of $\calF$ at the point $\bigl( (0,\frakb_1), (0,\frakb_2) \bigr) \in \wfrakg^{(1)} \times \wfrakg^{(1)}$ is nonzero. By Nakayama's lemma (and adjuction), it follows that there exists a nonzero morphism $\calF \to \calO_{\{(0,\frakb_1), (0,\frakb_2)\}}$. Consider the composition
\begin{equation}
\label{eq:morphism0}
\calK_w \ \twoheadrightarrow \ \calF \ \to \ \calO_{\{(0,\frakb_1), (0,\frakb_2)\}},
\end{equation}
which is also nonzero. By property $(\ddag)$, the morphism obtained by taking global sections is nonzero. In other words, there exists a morphism $\calO_{\wfrakg^{(1)} \times \wfrakg^{(1)}} \to \calK_w$ whose composition with \eqref{eq:morphism0} is nonzero. This is absurd since this morphism can be factored through a morphism
\[
\calO_{\wfrakg^{(1)} \times \wfrakg^{(1)}} \ \to \ \Gamma(\wfrakg^{(1)} \times \wfrakg^{(1)}, \, \calK_w) \otimes_{\bk} \calO_{\wfrakg^{(1)} \times \wfrakg^{(1)}} \ \to \ \calK_w \ \to \ \calF,
\]
which is zero by definition of $\calF$.

So, let us now prove $(\ddag)$. The idea of the proof is to translate this property in representation-theoretic terms; then it easily follows from Lemma \ref{lem:translationbabyVerma}. Consider a nonzero morphism \begin{equation} \label{eq:morphism1} \calK_w \to \calO_{\{(0,\frakb_1), (0,\frakb_2)\}}. \end{equation} Restricting this morphism to the formal neighborhood of the zero section of $\wfrakg^{(1)} \times \wfrakg^{(1)}$, and applying the equivalence $\widehat{\gamma}^{G \times G}_{(-2\rho,0)}$, we obtain a nonzero morphism
\begin{equation}
\label{eq:morphism2}
R\Gamma_{(-2\rho,0)}((\calM_{-2\rho} \boxtimes \calM_0) \otimes_{\wfrakg^{(1)} \times \wfrakg^{(1)}} \calK_w) \ \to \ \gamma_{(-2\rho,0)}^{G \times G}(\calO_{\{(0,\frakb_1), (0,\frakb_2)\}}).
\end{equation}
By definition (and using the K\"{u}nneth formula), $\gamma_{(-2\rho,0)}^{G \times G}(\calO_{\{(0,\frakb_1), (0,\frakb_2)\}}) \cong \gamma_{-2\rho}^G(\calO_{\{(0,\frakb_1)\}}) \otimes_{\bk} \gamma_0^G(\calO_{\{(0,\frakb_2)\}})$. And, by \cite[Proposition 3.1.4]{BMR}, we have isomorphisms
\[
\gamma_{-2\rho}^G(\calO_{\{(0,\frakb_1)\}}) \cong M_{(\frakb_1,0;0)}, \quad  \gamma_0^G(\calO_{\{(0,\frakb_2)\}}) \cong M_{(\frakb_2,0;2\rho)}.
\]
(Here we consider the baby Verma modules for the Lie algebra $\frakg$.)

Let $q : \wfrakg^{(1)} \times \wfrakg^{(1)} \to \frakg^* {}^{(1)} \times \wfrakg^{(1)}$ be the natural morphism, and let $\widehat{q}$ be the induced morphism on the formal neighborhoods of the zero sections. Let also $\eta$, respectively $\epsilon$, be the inclusion of the formal neighborhood of the zero section in $\wfrakg^{(1)} \times \wfrakg^{(1)}$, respectively $\frakg^* {}^{(1)} \times \wfrakg^{(1)}$, so that $q \circ \eta = \epsilon \circ \widehat{q}$.

By Lemma \ref{lem:K_w}, the object on the left hand side of morphism \eqref{eq:morphism2} is concentrated in degree $0$. Then, by Lemma \ref{lem:translationbabyVerma}, the morphism obtained by applying the functor $T_{(-2\rho,0)}^{(-\rho,0)}$ to this morphism is nonzero. (Observe that $T_{(-2\rho,0)}^{(-\rho,0)}=T_{(w_0 \bullet (-2\rho),0)}^{(w_0 \bullet (-\rho),0)}=T_{(0,0)}^{(-\rho,0)} $.) On the other hand, by \cite[Lemma 2.2.5]{BMR2}, the morphism obtained by applying $T_{(-2\rho,0)}^{(-\rho,0)}$ is a morphism \begin{equation*} \widehat{\gamma}^{G \times G}_{(-\rho,0)} (R(\widehat{q})_* \circ \eta^* \calK_w) \ \to \ \gamma_{(-\rho,0)}^{G \times G}(Rq_* \calO_{\{(0,\frakb_1),(0,\frakb_2)\}}). \end{equation*} There are isomorphisms $R(\widehat{q})_* \circ \eta^* \calK_w \cong \epsilon^* \circ Rq_* \calK_w$ (see \cite[Theorem 4.1.5]{EGA3_1}) and $Rq_* \calO_{\{(0,\frakb_1),(0,\frakb_2)\}} \cong \calO_{\{0,(0,\frakb_2)\}}$. It follows that the morphism \begin{equation} \label{eq:morphism3} R(p_2)_* \calK_w \ \to \ \calO_{\{(\frakb_2,0)\}} \end{equation} obtained by applying the functor $R(p_2)_*$ to \eqref{eq:morphism1} is nonzero. Here, $p_2 : \wfrakg^{(1)} \times \wfrakg^{(1)} \to \wfrakg^{(1)}$ is the projection on the second factor. 

By Lemma \ref{lem:globalsections}, we have an isomorphism $R(p_2)_* \calK_w \cong \calO_{\wfrakg^{(1)}}$. Hence the morphism obtained by applying $R\Gamma(\wfrakg^{(1)},-)$ to the morphism \eqref{eq:morphism3} is nonzero. But, as $R\Gamma(\wfrakg^{(1)} \times \wfrakg^{(1)},-) \cong R\Gamma(\wfrakg^{(1)},-) \circ R(p_2)_*$, the latter morphism is obtained by taking global sections of \eqref{eq:morphism1}. This finishes the proof.\end{proof}

Let us remark that, using similar arguments, one can prove the following fact (which will not be used in this paper).

\begin{prop}

Let $\lambda \in \bbX$ be in the fundamental alcove. Let $\calF$ be an object of $\Coh^{\Gm}(\wfrakg^{(1)})$, and let $\widehat{\calF}$ be its restriction to the formal neighborhood of the zero section $\calB^{(1)} \subset \wfrakg^{(1)}$. 

Assume that the object $\widehat{\gamma}_{w_0 \bullet \lambda}(\widehat{\calF})$ is concentrated in degree $0$. Then $\calF$ is generated by its global sections.\qed

\end{prop}

\subsection{Irreducibility of a fibre product}
\label{ss:irreducibility}

The main result of this subsection is Proposition \ref{prop:irreducible}. Before we can prove it, we need several easy preliminary results.

%
%
%

The following fact follows from the explicit description of $Z_s$ in
\cite{RAct}. We denote as above by $p_1$ and $p_2$ the projections $\wfrakg^2 \to \wfrakg$.

\begin{lem} \label{lem:fibers}

Let $s \in \scS$. The fiber of the projection $p_1 : Z_{s} \to \wfrakg$, respectively $p_2 : Z_{s} \to \wfrakg$, over $(X,gB)$ is one point if $X_{|g \cdot \fraksl(2,s)} \neq 0$, and
is isomorphic to $\mathbb{P}^1$ otherwise.\qed

\end{lem}

We denote by $F_s$ the closed subvariety of codimension $2$ of
$\wfrakg$ defined by $F_s:=\{(X,gB) \in \wfrakg \mid X_{|g \cdot
  \fraksl(2,s)} = 0\}$.

Recall that there is a natural morphism
\[
\nu: \wfrakg \to
\frakh^*,
\]
sending $(X,gB)$ to the restriction $X_{|g \cdot \frakb}$, considered
as a linear form on $(g \cdot \frakb) / (g \cdot \frakn) \cong
\frakh$. Then for any $w \in W$ the image of $Z_w$ under $\nu \times
\nu : \wfrakg \times \wfrakg \to \frakh^* \times \frakh^*$ is the
graph ${\rm Graph}(w,\frakh^*)$ of the action of $w$
on $\frakh^*$. Indeed the inverse image of ${\rm
  Graph}(w,\frakh^*)$ under $(\nu \times \nu)$ contains $Z_w^{\reg}$,
hence also $Z_w$; it follows that $(\nu \times \nu)(Z_w) \subseteq {\rm
  Graph}(w,\frakh^*)$; we conclude using the fact the morphism
$p_1 : Z_w \to \wfrakg$ is proper and birational, hence surjective.

\begin{lem} \label{lem:nu}

Let $w \in W$, and let $v \in W$ such that $v < w$. Let $V$ be the fiber of $Z_w$ over $(B/B,v^{-1}B/B) \in \calB \times \calB$ (a closed subvariety of $(\frakg/\frakn)^*$). Then $\nu(V \times \{B/B\})$ is included in the space of fixed points of $v^{-1}w$ on $\frakh^*$.

\end{lem}

\begin{proof} Let $X \in V$. By definition, $\nu(X,v^{-1}B/B)=v \cdot \nu(X,B/B)$. On the other hand, as $(X,B/B,v^{-1}B/B) \in Z_w$, by the remarks above we have $\nu(X,v^{-1}B/B)=w \cdot \nu(X,B/B)$. We deduce that $v \cdot \nu(X,B/B)=w \cdot \nu(X,B/B)$, which gives the result. \end{proof}

Recall the remarks of \S \ref{ss:kernels-statement}. We set $d:=\dim(\frakg)$.

\begin{cor} \label{cor:dimension}

Let $w \in W$ and $s \in \scS$ such that $ws>w$. Then
\[
\dim\bigl( (F_s \times \wfrakg) \cap Z_w \bigr) \leq d - 2.
\]

\end{cor}

\begin{proof} For simplicity, let us denote the variety $(F_s \times \wfrakg) \cap Z_w$ by $Y_w^s$. (We will use this notation only in this proof.) First, $Y_w^s$ is included strictly in $Z_w$, hence has dimension lower than $d - 1$. Assume that it has dimension $d - 1$. By $G$-equivariance, and as $G$ has only finitely many orbits on $\calB \times \calB$, there exists $u \leq w$ such that the restriction $(Y_w^s)_{|\frakX^0_{u^{-1}}}$ has dimension $d-1$. Moreover, $u \neq w$ as the restriction of $Y^s_w$ to $\frakX_{w^{-1}}^0$ has dimension $\dim(\frakg)-2$. (This follows from our assumption $ws>w$ or, equivalently, $w(\alpha) > 0$ for $s=s_{\alpha}$, $\alpha \in \Sigma$.) By $G$-equivariance again, the fiber of $Y_w^s$ over $(B/B,u^{-1}B/B)$ has
dimension $d-1-\dim(\mathfrak{X}_{u^{-1}}^0)=d-1-\ell(u)-\dim(\calB)$. On the other hand, this fiber is
included in the fiber of $(\wfrakg \times_{\frakg^*} \wfrakg) \cap
(F_s \times \wfrakg)$ over $(B/B,u^{-1}B/B)$, which is itself included in \[ V:=\{X \in \frakg^* \mid X_{|\frakn + u^{-1} \cdot \frakn}=0 \ \text{and} \
X(h_s)=0\}.\] The subspace $V$ has dimension $d-1-\ell(u)-\dim(\calB)$. Hence the fiber of $Y_w^s$ over $(B/B,u^{-1}B/B)$
equals $V$. In particular, the fiber of $Z_w$ over $(B/B,u^{-1}B/B)$
contains $V$.

Now we derive a contradiction. By definition we have $\nu(V \times \{B/B\})=\{X \in (\frakb/\frakn)^* \mid
X(h_s)=0\}$, i.e.~$\nu(V \times \{B/B\}) \subset \frakh^*$ is the reflection hyperplane of $s$. By Lemma \ref{lem:nu}, this subspace is included in the fixed points of $u^{-1}w$. Hence either $u^{-1}w=s$, or $u^{-1}w=1$. This is absurd since $u < w < ws$. \end{proof}

Let $w \in W$ and $s \in \scS$. Consider the scheme $Z_s \times_{\wfrakg} Z_w$, where the morphism $Z_s \to \wfrakg$ (respectively $Z_w \to \wfrakg$) is induced by the second (respectively the first) projection.

\begin{prop} \label{prop:irreducible}

Let $w \in W$ and $s \in \scS$ such that $ws>w$. The scheme $Z_s
\times_{\wfrakg} Z_w$ is irreducible, of dimension $\dim(\frakg)$.

\end{prop}

\begin{proof} The scheme $Z_s \times_{\wfrakg} Z_w$ is the
scheme-theoretic intersection of the subvarieties $(Z_s \times
\wfrakg)$ and $(\wfrakg \times Z_w)$ of $\wfrakg^3$. Each of these
subvarieties has dimension $2\dim(\frakg)$. As $\wfrakg$ has a finite covering
by open subsets isomorphic to $\mathbb{A}^{\dim(\frakg)}$, the dimension of each
irreducible component of $Z_s \times_{\wfrakg} Z_w$ is at least $\dim(\frakg)$
(see \cite[Proposition I.7.1]{HARAG}). 


We have $(Z_s \times_{\wfrakg} Z_w) \cap (\wfrakg_{\reg})^3 = Z_s^{\reg} \times_{\wfrakg_{\reg}} Z_w^{\reg}$, hence this intersection is irreducible (because it is isomorphic to $\wfrakg_{\reg}$). Hence any irreducible component of $Z_s \times_{\wfrakg} Z_w$ is either the closure of $Z_s^{\reg} \times_{\wfrakg_{\reg}} Z_w^{\reg}$, or is included in $(\wfrakg \smallsetminus \wfrakg_{\reg})^3$. Assume that there is a component $Y$ included in $(\wfrakg \smallsetminus \wfrakg_{\reg})^3$. By the arguments above, $\dim(Y) \geq \dim(\frakg)$. Consider the image $Y'$ of $Y$ under the projection $p_{2,3} : Z_s \times_{\wfrakg} Z_w \to Z_w$. Then $Y'$ is strictly included in $Z_w$, hence has dimension lower than $\dim(\frakg) - 1$. As the fibers of $p_{2,3}$ have dimension at most $1$ (see Lemma \ref{lem:fibers}), we have $\dim(Y')=\dim(\frakg)-1$, $\dim(Y)=\dim(\frakg)$, and all the fibers of the restriction $(p_{2,3})_{|Y}$ have dimension exactly $1$. It follows that $Y' \subset (F_s \times \wfrakg) \cap Z_w$. This is absurd, since $\dim((F_s \times \wfrakg) \cap Z_w) \leq \dim(\frakg)-2$ by Corollary \ref{cor:dimension}.\end{proof}

\subsection{End of the proof of Theorem \ref{thm:descriptionkernels} (case of $\wfrakg$)} \label{ss:endproof}

We can finally finish the proof of Theorem \ref{thm:descriptionkernels} in the case of $\wfrakg$.

\begin{proof}[Proof of Theorem {\rm \ref{thm:descriptionkernels}}]
We prove the theorem by induction on the Bruhat order. First, the statement is clear by definition if $\ell(w)$ is $0$ or $1$. Now assume it is known for $w$, and let $s \in \scS$ be a simple reflection such that $ws>w$. We only have to prove that $\calK_{ws} \cong \calO_{Z_{ws}}$ as $(G \times \Gm)$-equivariant sheaves. Indeed, once we know this isomorphism, the fact that $\calK_{ws}^{\dag}$ is a sheaf implies that $Z_{ws}$ is Cohen--Macaulay, and that $\calK_{ws}^{\dag}$ is its dualizing sheaf (see Equation \eqref{eq:duality}).

By definition of convolution, and the induction hypothesis,
\[
\calK_{ws} \ \cong \ \calK_w \star \calO_{Z_s} \ \cong \ R(p_{1,3})_*(\calO_{Z_s \times \wfrakg} \, \lotimes_{\wfrakg^3} \, \calO_{\wfrakg \times Z_w}).
\]
By induction hypothesis, $\wfrakg \times Z_w$ is a Cohen--Macaulay scheme. Moreover, $Z_s \times \wfrakg$ is defined locally in $\wfrakg^3$ be a regular sequence of length $\dim(\frakg)$. (This is a general fact for smooth subvarieties of smooth varieties, and it is checked explicitly in this case in \cite{RAct}.) We have proved in Proposition \ref{prop:irreducible} that $(Z_s \times \wfrakg) \cap (\wfrakg \times Z_w)=Z_s \times_{\wfrakg} Z_w$ has dimension $\dim(\frakg)$. Using a Koszul complex (\cite[(18.D) Theorem 43]{MATCA}) and \cite[(16.B) Theorem 31]{MATCA}, the derived tensor product $\calO_{Z_s \times \wfrakg} \, \lotimes_{\calO_{\wfrakg^3}} \, \calO_{\wfrakg \times Z_w}$ is concentrated in degree $0$. Hence it equals $\calO_{Z_s \times_{\wfrakg} Z_w}$. Moreover, $Z_s \times_{\wfrakg} Z_w$ is Cohen--Macaulay (see \cite[(16.A) Theorem 30]{MATCA}). By Proposition \ref{prop:irreducible} it is also irreducible, and it is smooth on an open subscheme (e.g.~the intersection with $\wfrakg_{\reg}^3$). Hence it is reduced (see \cite[p.~125]{MATCA}). It follows that $p_{1,3}$ induces a birational and proper (hence surjective) morphism $p_{1,3} : Z_s \times_{\wfrakg} Z_w \to Z_{ws}$, and that we have
\begin{equation}
\label{eq:K_wsdirectimage}
\calK_{ws} \ \cong \ (p_{1,3})_* \calO_{Z_s \times_{\wfrakg} Z_w}
\end{equation}
(Recall that $\calK_{ws}$ is a sheaf by Proposition \ref{prop:K_wsheaf}.)

On the other hand, let
\[
f : \wfrakg \times \wfrakg \ \xrightarrow{p_2} \ \wfrakg \ \to \ (\frakg^* \times_{\frakh^*/W} \frakh^*)
\]
be the natural morphism. By Proposition \ref{prop:globalgeneration}, $\calK_{ws}$ is globally generated. Hence the adjunction morphism $f^* f_* \calK_w \to \calK_w$ is surjective. Using also Lemma \ref{lem:globalsections}, it follows that there exists a natural surjective $(G \times \Gm)$-equivariant morphism
\[
\calO_{\wfrakg \times \wfrakg} \ \twoheadrightarrow \ \calK_{ws}.
\]
In other words, $\calK_{ws}$ is the structure sheaf of a closed subscheme of $\wfrakg \times \wfrakg$. By equation \eqref{eq:K_wsdirectimage}, this subscheme is reduced, and coincides with $Z_{ws}$. This finishes the proof.
\end{proof}

\subsection{Proof of Theorem \ref{thm:descriptionkernels-R} for $\wfrakg_R$}
\label{ss:proof-thm-description-kernels}

We come back to the setting of Theorem \ref{thm:descriptionkernels-R}. In particular, let $R=\Z[\frac{1}{h!}]$, let $w \in W$, and let $w=s_1 \cdots s_n$ be a reduced expression.

Consider the object
\[
\calK_{w,R} \ := \ \calO_{Z_{s_1,R}} \star \cdots \star \calO_{Z_{s_n,R}}
\]
of $\calD^b \Coh^{G_R \times_R (\Gm)_R}(\wfrakg_R \times_R \wfrakg_R)$. Note that for every prime $p>h$, $\calK_{w,R} \, \lotimes_R \, \overline{\F_p}$ is the object ``$\calK_w$'' of \S \ref{ss:kernels-statement} for the field $\bk=\overline{\F_p}$. Hence, by Lemma \ref{lem:specialization}\eqref{it:specialization-degree-0} and Proposition \ref{prop:K_wsheaf}, $\calK_{w,R}$ is a coherent sheaf, which is flat over $R$.

Let $M:=R\Gamma(\wfrakg_R \times_R \wfrakg_R, \, \calK_{w,R})$. By Lemma \ref{lem:specialization}\eqref{it:specialization-degree-0} again and Lemma \ref{lem:globalsections}, $M$ is concentrated in degree $0$. By the same arguments as in the proof of Proposition \ref{prop:isom-Delta}, we have $M^{(\Gm)_R} \cong R$. Hence also
\[
\Gamma(\wfrakg_R \times_R \wfrakg_R, \, \calK_{w,R})^{G_R \times_R (\Gm)_R} \ \cong \ R.
\]
It follows that there exists a non-zero $(G_R \times_R (\Gm)_R)$-equivariant morphism
\[
\phi : \calO_{\wfrakg_R \times_R \wfrakg_R} \to \calK_w
\]
(uniquely defined up to a scalar which is invertible in $R$.) By the same arguments as in the proof of Proposition \ref{prop:isom-Delta}, $\phi$ is surjective. In other words, we have an isomorphism $\calO_{X_w} \xrightarrow{\sim} \calK_w$ for some $(G_R \times_R (\Gm)_R)$-stable closed subscheme $X_w \subset \wfrakg_R \times_R \wfrakg_R$. Again by the same arguments as in the proof of Proposition \ref{prop:isom-Delta}, $X_w$ is reduced, and its restriction to $(\wfrakg_R \times_R \wfrakg_R) \smallsetminus Z_{w,R}$ is empty. (Here we use that $Z_{w,R} \times_{\mathrm{Spec}(R)} \mathrm{Spec(\overline{\F_p})}$ and the scheme $Z_w$ defined in \S \ref{ss:kernels-statement} for $\bk=\overline{\F_p}$ have the same underlying topological space.) Hence $X_w$ is a reduced closed subscheme of $Z_{w,R}$.

On the other hand, let $U_w$ be the inverse image under the projection $\wfrakg_R \times_R \wfrakg_R \to \calB_R \times_R \calB_R$ of the $G_R$-orbit of $(B_R/B_R, w^{-1}B_R/B_R) \in \calB_R \times_R \calB_R$, and let $Z_{w,R}^*$ be the intersection of $U_w$ with $\wfrakg_R \times_{\frakg^*_R} \wfrakg_R$ (so that, by definition, $Z_{w,R}$ is the closure of $Z_{w,R}^*$). It is easy to check (using the fact that the Bott--Samelson resolution associated with the reduced expression $w^{-1}=s_n \cdots s_1$ is an isomorphism over the inverse image of the orbit of $w^{-1}$) that we have
\[
\calK_{w,R}|_{U_w} \ \cong \ \calO_{Z_{w,R}^*}.
\]
Hence $X_w \cap U_w = Z_{w,R}^*$, which implies that $X_w$ contains $Z_{w,R}$. We deduce that $X_w=Z_{w,R}$, which finishes the proof of Theorem \ref{thm:descriptionkernels-R} in the case of $\wfrakg_R$.

\begin{remark}
\label{rk:reduction-Z_w}
Let $\bk$ be an algebraically closed field of characteristic $p>h$. As we have observed in the proof above, $\calK_{w,R} \, \lotimes_R \, \bk$ is the object ``$\calK_w$'' of \S \ref{ss:kernels-statement} for the field $\bk$ (by compatibility of convolution with change of scalars). On the other hand, the isomorphism $\calK_{w,R} \cong \calO_{Z_{w,R}}$ implies that $\calK_{w,R} \otimes_R \bk \cong \calO_{Z_{w,R} \times_{\mathrm{Spec}(R)} \mathrm{Spec}(\bk)}$. We deduce that
\[
Z_{w,R} \times_{\mathrm{Spec}(R)} \mathrm{Spec}(\bk) \ \cong \ Z_w,
\]
where ``$Z_w$'' is the scheme defined in \S \ref{ss:kernels-statement} for the field $\bk$. In other words, $Z_{w,R} \times_{\mathrm{Spec}(R)} \mathrm{Spec}(\bk)$ is reduced.
\end{remark}

\subsection{Geometric action for $\wcalN$ and $\wcalN_R$}
\label{ss:action-wcalN}

In this subsection we prove Theorems \ref{thm:descriptionkernels-R} and \ref{thm:descriptionkernels} in the case of $\wcalN$. For simplicity, we only treat the case of a field; the case of $R$ is similar. Hence we fix again an algebraically closed field $\bk$ of characteristic $p > h$.

Let $w \in W$. Recall the scheme $Z_w'$ defined in \S \ref{ss:kernels-statement}. As $Z_w$ is Cohen--Macaulay (by Theorem \ref{thm:descriptionkernels}), and $\dim(Z_w') \leq \dim(\wcalN \times_{\frakg^*} \wcalN)=\dim(\frakg) - \dim(\frakt)$, one easily checks that $Z_w'$ is Cohen--Macaulay, and that
\[
\calO_{Z_w'} \cong \, \calO_{Z_w} \, \lotimes_{\wfrakg^2} \, \calO_{\wcalN \times \wfrakg}
\]
(see \S \ref{ss:endproof} for a similar argument).

\begin{lem}
\label{lem:Z_w'}
 
$Z_w'$ is a closed subscheme of $\wcalN \times \wcalN$.

\end{lem}

\begin{proof} It is sufficient to prove that the morphism $\mathrm{S}(\frakh) \to \calO_{Z_w'}$ induced by the morphism \[ Z_w' \hookrightarrow \wfrakg \times \wfrakg \xrightarrow{p_2} \wfrakg \to \frakh^* \] is zero. But this morphisms coincides with \[ Z_w' \hookrightarrow \wfrakg \times \wfrakg \xrightarrow{p_1} \wfrakg \to \frakh^* \xrightarrow{w} \frakh^* \] (see the remarks before Lemma \ref{lem:nu}). Hence indeed it is zero, by definition of $Z_w'$.\end{proof}

It follows in particular from this lemma that $\calO_{Z_w'}$ can be considered as a coherent sheaf on $\wcalN \times \wcalN$. Now an easy argument, similar to that of \cite[Corollary 4.3]{RAct}, proves the isomorphism of Theorem \ref{thm:descriptionkernels}. The description of the dualizing sheaf for $Z_w'$ can be proved using an analogue of formula \eqref{eq:duality} for $\wcalN$.

\subsection{Application to homology and $K$-theory} \label{ss:homology}

In \cite[\S\S 6.1, 6.2]{RAct}, we have explained how the braid group action of Theorem \ref{thm:existenceaction}, for $\bk=\mathbb{C}$, is related to the $K$-theoretic description of the affine Hecke algebra, and Springer's geometric construction of the representations of $W$. Here we explain some consequences of Theorem \ref{thm:descriptionkernels} in this context. (Recall that, by Remark \ref{rk:thm-description-kernels}\eqref{it:rk-thm-description-kernels-char-0}, Theorem \ref{thm:descriptionkernels} is also true for $\bk=\mathbb{C}$.)

For the moment, consider an arbitrary algebraically closed field $\bk$, and consider the fiber product
\[
Z \ := \ \wfrakg \times_{\frakg^*} \wfrakg.
\]

\begin{lem}

The scheme $Z$ is reduced.

\end{lem}

\begin{proof}
It is well known that $Z$ has dimension $d:=\dim(\frakg)$ (see e.g.~\cite[\S 10]{JANNil}), and is defined by $d$ equations in the smooth variety $\wfrakg \times \wfrakg$. Hence it is Cohen--Macaulay (see \cite[(16.A) Theorem 30 and (16.B) Theorem 31]{MATCA}). Moreover, it is smooth on the inverse image of $\frakg^*_{\rs} \subset \frakg^*$, which is dense. One concludes using \cite[p.~125]{MATCA}.
\end{proof}

\begin{remark}
One can interpret the schemes $Z_w$ from yet another point of view using this definition: they are the irreducible components of the variety $Z$ (see e.g.~\cite[\S 10]{JANNil}).
\end{remark}

On the other hand, the scheme-theoretic fiber product $\wcalN \times_{\frakg^*} \wcalN$ is not reduced in general. (For example, it is not reduced for $G=\SL(2,\bk)$ if $\mathrm{char}(\bk)=2$.) We set
\[
Z' \ := \ (\wcalN \times_{\frakg^*} \wcalN)_{\mathrm{red}}.
\]

Now for simplicity we specialize to the case $\bk=\mathbb{C}$. First, recall the algebra isomorphism
\begin{equation}
\label{eq:hecke}
\calH_{\aff} \ \xrightarrow{\sim} \ K^{G \times \mathbb{C}^{\times}}(Z')
\end{equation}
due to Kazhdan--Lusztig and Ginzburg. Here $G \times \mathbb{C}^{\times}$ acts on $Z'$ via $(g,t) \cdot (X,g_1B, g_2B) := (t^{-2} g \cdot X, gg_1B, gg_2B)$,
\[
\calH_{\aff}\ := \ \mathbb{Z}[v,v^{-1}][\bB_{\aff}] \, / \, \langle (T_s+v^{-1})(T_s-v), \ s \in \scS \rangle
\]
is the (extended) affine Hecke algebra, and we consider Lusztig's isomorphism defined in \cite[Theorem 8.6]{LUSBas} (and not Ginzburg's isomorphism, defined in \cite{CG}, which is slightly different). In order to follow Lusztig's conventions, we consider the algebra structure on the right-hand side of \eqref{eq:hecke} which is induced by the analogue of the convolution product on the category $\calD^b \Coh_{Z'}(\wcalN \times \wcalN)$ defined in \S \ref{ss:convolution1}, but where the role of the two copies of $\wcalN$ is exchanged. By \cite[\S 6.1]{RAct}, this isomorphism sends $T_s$ to $(-v^{-1}) \cdot [\calO_{Z_s'}]$.

Similarly, there is an algebra isomorphism $K^{G \times \mathbb{C}^{\times}}(Z') \cong K^{G \times \mathbb{C}^{\times}}(Z)$ (see e.g.~\cite[\S 6.1]{MRHec}). Composing it with \eqref{eq:hecke} we obtain an algebra isomorphism
\begin{equation}
\label{eq:hecke2}
\calH_{\aff} \ \xrightarrow{\sim} \ K^{G \times \mathbb{C}^{\times}}(Z),
\end{equation}
which sends $T_s$ to $(-v^{-1}) \cdot [\calO_{Z_s}]$.

Theorem \ref{thm:descriptionkernels} implies the following.

\begin{prop}

Let $w \in W$. Under isomorphism \eqref{eq:hecke2}, respectively isomorphism \eqref{eq:hecke}, $T_w \in \calH_{\aff}$ is sent to $(-v^{-1})^{\ell(w)} \cdot [\calO_{Z_{w^{-1}}}]$, respectively to $(-v^{-1})^{\ell(w)} \cdot [\calO_{Z_{w^{-1}}'}]$.\qed

\end{prop}

In particular, this result gives a geometric description of two standard bases of the affine Hecke algebra, given by $\{T_w \theta_x \mid w \in W, \ x \in \bbX\}$ and $\{\theta_x T_w \mid w \in W, \ x \in \bbX\}$. Namely, these bases are given by the classes of some shifts of the following $(G \times \mathbb{C}^{\times})$-equivariant coherent sheaves on $Z$:
\[
\calO_{Z_{w^{-1}}}(x,0), \quad \text{respectively} \quad \calO_{Z_{w^{-1}}}(0,x). \]

Now, consider the rational top Borel--Moore homology $\calH^{\mathrm{BM}}_{\mathrm{top}}(Z')$ of the variety $Z'$. By \cite[Proposition 5.3]{G1}, \cite[Theorem 3.4.1]{CG} or \cite{KT}, there exists an algebra isomorphism 
\begin{equation}
\label{eq:BorelMoore}
\mathbb{Q}[W] \ \xrightarrow{\sim} \ \calH^{\mathrm{BM}}_{\mathrm{top}}(Z'),
\end{equation}
where the algebra structure on the right-hand side is again given by convolution. (Here again, to follow the conventions of \cite{G1, KT, CG}, the role of the two copies of $\wcalN$ is exchanged in the definition of the convolution.) This isomorphism can be defined as follows. For $w \in W$, consider the regular holonomic system $\frakM_w$ on $\calB$ corresponding to the Verma module with highest weight $-w(\rho) - \rho$ under Beilinson--Bernstein's equivalence (see \cite[\S 2]{KT} for details). Then $\frakM_w$ is $B$-equivariant, hence induces a $G$-equivariant regular holonomic system $\widetilde{\frakM}_w$ on $\calB \times \calB$. Then isomorphism \eqref{eq:BorelMoore} sends $w$ to the characteristic cycle $\mathbf{Ch}(\widetilde{\frakM}_w)$ of $\widetilde{\frakM}_w$. By \cite[\S 6.2]{RAct}, the image of $s \in \scS$ is also the characteristic class $[Z'_s]$ of $Z'_s$.

Now Theorem \ref{thm:descriptionkernels} implies the following.

\begin{prop} \label{prop:BorelMoore}

Let $w \in W$. Under isomorphism \eqref{eq:BorelMoore}, $w \in W$ is sent to the characteristic class $[Z_{w^{-1}}'] \in \calH^{\mathrm{BM}}_{\mathrm{top}}(Z')$. In particular we have an equality $\mathbf{Ch}(\widetilde{\frakM}_w)=[Z_{w^{-1}}']$.

\end{prop}

\begin{remark}
The equality $\mathbf{Ch}(\widetilde{\frakM}_w)=[Z_{w^{-1}}']$ can also be derived from \cite[Equation (6.2.3)]{G1}.
\end{remark}

Proposition \ref{prop:BorelMoore} can be used to check that $Z_w'$ is not reduced in general. Recall that the irreducible components of the Steinberg variety $Z'$ are the closures of the conormal bundles to the $G$-orbits on $\calB \times \calB$ (see \cite{CG} or \cite[\S 4.1]{KT}). For $y \in W$, let us denote by 
\[
Y_y:=\overline{T^*_{\frakX^0_y}(\calB \times \calB)}
\]
the corresponding component, endowed with the reduced subscheme structure. It can be easily checked that the reduced subscheme associated to $Z_w'$ is 
\[
(Z_w')_{\mathrm{red}} = \bigcup_{y \leq w^{-1}} Y_y.
\]
(Here, $\leq$ is the Bruhat order.) However, the multiplicity of $Y_y$ in $Z_w'$ can be more than $1$, which will prove that $Z_w'$ is not reduced in these cases. 

Indeed, for $w \in W$, in addition to $\frakM_w$ defined above, consider the regular holonomic system $\frakL_w$ on $\calB$ which corresponds to the simple $\calU \frakg$-module with highest weight $-w(\rho) - \rho$ under Beilinson--Bernstein's equivalence (see again \cite[\S 2]{KT} for details). By Proposition \ref{prop:BorelMoore}, the multiplicity of $Y_y$ in $Z_w'$ is the coefficient of the cycle $[\overline{T^*_{ByB/B}(\calB)}]$ in the decomposition of the characteristic cycle $\mathbf{Ch}(\frakM_w)$ as a sum of elements $[\overline{T^*_{BzB/B}(\calB)}]$, $z \in W$. To prove that one of these coefficients is greater than $1$, it is sufficient to prove that one of the coefficients of the decomposition of $\mathbf{Ch}(\frakM_w)$ as a sum of elements $\mathbf{Ch}(\frakL_y)$ is greater than $1$. However, the latter coefficients are given by values at $1$ of Kazhdan--Lusztig polynomials, which are related to singularities of Schubert varieties. 

For example, consider the group $G=\mathrm{SL}(4)$. Let $s_1$, $s_2$, $s_3$ be the standard generators of the Weyl group $W$, and let $w=s_2 s_1 s_3 s_2$, $y=s_2$. Then $P_{y,w}(q)=1+q$, hence the coefficient of $\mathbf{Ch}(\frakL_{w_0 w})$ in the decomposition of $\mathbf{Ch}(\frakM_{w_0 y})$ is $2$. It follows that the multiplicity of $Y_{w_0 w}$ in $Z_{y^{-1} w_0}'$ is at least\footnote{In fact, it is checked in \cite{KS} that if $n \leq 7$, for the group $G=\mathrm{SL}(n)$ we have $\mathbf{Ch}(\frakL_y)=[\overline{T^*_{ByB/B}(\calB)}]$ for $y \in W$. Hence here the multiplicity is exactly $2$.} $2$, hence that $Z_{y^{-1} w_0}'$ is not reduced.

\section{Generalities on dg-schemes} \label{sec:dgschemes}

In the next two sections, we develop a general framework to define group actions on derived categories of coherent (dg-)sheaves on (dg-)schemes. We will use this framework to extend the action of Theorem \ref{thm:existenceaction} to other related categories.

First, in this section we extend the formalism of dg-schemes of \cite[\S 1]{R2} to a setting more adapted to \emph{quasi-coherent} dg-sheaves. Then we extend the base change theorem and the projection formula to this setting.

In this section and the next one, a \emph{scheme} is always assumed to be separated and noetherian of
finite Krull dimension. It follows that the morphisms of schemes are
always quasi-compact and separated.

\subsection{Definitions} \label{ss:definitionsdg}

Recall the following definitions (\cite[\S 2.2]{CK}, \cite[\S 1.8]{R2}).

\begin{defin}
\label{def:dg-schemes}

\begin{enumerate}
\item A \emph{dg-scheme} is a pair $(X,\calA)$ where $(X,\calO_{X})$ is a scheme
(with the conventions stated above), and $\calA$ is a quasi-coherent,
non-positively graded, graded-commutative sheaf of $\calO_X$-dg-algebras on $X$.
\item A \emph{morphism of dg-schemes} $f : (X,\calA_X) \to (Y,\calA_Y)$ is
the data of a morphism of schemes $f_0 : X \to Y$, and a morphism of sheaves of dg-algebras $(f_0)^* \calA_Y \to \calA_X$.

\end{enumerate}

\end{defin}

Let us fix a dg-scheme $(X,\calA)$. Remark that, as the image of
$\calO_X$ in $\calA$ is in degree $0$, it is killed by the
differential $d_{\calA}$ of $\calA$. It follows that $d_{\calA}$ is
$\calO_X$-linear. The same applies to all sheaves of
$\calA$-dg-modules.

We denote by $\calC(X,\calA)$ the category of all sheaves of
$\calA$-dg-modules, and by $\calC^{\qc}(X,\calA)$ the full subcategory
of dg-modules whose cohomology sheaves are quasi-coherent over $\calO_X$. We denote by
$\calD(X,\calA)$ and $\calD^{\qc}(X,\calA)$ the corresponding derived
categories. Note that $\calD^{\qc}(X,\calA)$ is equivalent to the full
subcategory of $\calD(X,\calA)$ whose objects are in
$\calC^{\qc}(X,\calA)$. The category $\calD(X,\calA)$ is triangulated
in a natural way. It is not clear from the definition that
$\calD^{\qc}(X,\calA)$ is a triangulated subcategory; we will eventually prove that this is the case under our hypotheses (see Proposition \ref{prop:equivdgalg}).

We denote by $\calC \QCoh(X,\calA)$ the category of
sheaves of $\calA$-dg-modules which are quasi-coherent over $\calO_X$, and by $\calD
\QCoh(X,\calA)$ the corresponding derived category.

Recall the following definition (see \cite[Definition 1.1]{Sp}, \cite[Definition 1.3.1]{R2}).

\begin{defin}

An object $\calF$ of $\calC(X,\calA)$ (respectively $\calC \QCoh(X,\calA)$)
is said to be \emph{K-injective} if for any acyclic object $\calG$ in
$\calC(X,\calA)$ (respectively $\calC \QCoh(X,\calA)$), the complex of
abelian groups $\Hom_{\calA}(\calG,\calF)$ is acyclic.

\end{defin}

Let us consider $\calO_X$ as a sheaf of dg-algebras
concentrated in degree $0$, with trivial differential. We have defined
the categories $\calC(X,\calO_X)$,
$\calC^{\qc}(X,\calO_X)$, $\calC \QCoh(X,\calO_X)$ and the
corresponding derived categories. Recall that the forgetful functor
\[ \For^{\calO_X} : \calC \QCoh(X,\calO_X) \to \calC(X,\calO_X) \] has a
right adjoint \[ Q^{\calO_X} : \calC(X,\calO_X) \to \calC \QCoh(X,\calO_X), \]
called the
\emph{quasi-coherator} (see \cite[p.~187, Lemme 3.2]{SGA6}). As $\calC(X,\calO_X)$
has enough K-injectives (see \cite[Theorem 4.5]{Sp}), $Q^{\calO_X}$ admits a right derived
functor \[ RQ^{\calO_X} : \calD(X,\calO_X) \to \calD\QCoh(X,\calO_X). \] The
functor $Q^{\calO_X}$ sends K-injective objects of $\calC(X,\calO_X)$ to
K-injective objects of $\calC \QCoh(X,\calO_X)$ (because it has an
exact left adjoint functor). One easily deduces that $RQ^{\calO_X}$ is right adjoint to
the forgetful functor $\For^{\calO_X} :
\calD \QCoh(X,\calO_X) \to \calD(X,\calO_X)$. The
functor $Q^{\calO_X}$ also induces a functor \[ Q^{\calA} : \calC(X,\calA) \to \calC
\QCoh(X,\calA), \] which is right adjoint to \[ \For^{\calA} : \calC
\QCoh(X,\calA) \to \calC(X,\calA).\] Under our hypotheses
$\calC(X,\calA)$ has enough K-injectives (see \cite[Theorem 1.3.6]{R2}), hence $Q^{\calA}$ has a right-derived functor $RQ^{\calA}$,
which is right adjoint to $\For^{\calA} : \calD \QCoh(X,\calA) \to
\calD(X,\calA)$ (for the same reasons as above).

\begin{prop}
\label{prop:thmSGA}

\begin{enumerate}
\item 
\label{it:prop-thmSGA-1}
The functors
\begin{align*} 
\For^{\calO_X} : & \ \calD \QCoh(X,\calO_X) \to
\calD^{\qc}(X,\calO_X), \\ 
RQ^{\calO_X} : & \ \calD^{\qc}(X,\calO_X) \to \calD
\QCoh(X,\calO_X) 
\end{align*}
are quasi-inverse equivalences of categories.
\item 
\label{it:prop-thmSGA-2}
$\calD^{\qc}(X,\calO_X)$ is a triangulated subcategory
of $\calD(X,\calO_X)$.
\end{enumerate}

\end{prop}

\begin{proof}
Statement \eqref{it:prop-thmSGA-1} is proved in \cite[p.~191, Proposition 3.7]{SGA6} (see also \cite[Proposition 1.3]{AJL} for a more general version).

Let us prove \eqref{it:prop-thmSGA-2}. Let
\[
\calF \xrightarrow{f} \calG \to \calH \xrightarrow{+1}
\]
be a distinguished triangle in $\calD(X,\calO_X)$, and assume that $\calF, \calG$ have quasi-coherent cohomology. By \eqref{it:prop-thmSGA-1} there exists $\calF'$, $\calG'$ in $\calD \QCoh(X,\calO_X)$, and a morphism $f' :
\calF' \to \calG'$ such that $\calF=\For^{\calO_X}(\calF')$,
$\calG=\For^{\calO_X}(\calG')$, $f=\For^{\calO_X}(f')$. By usual properties of triangulated categories, one can complete the morphism $f'$ to a distinguished triangle
\[
\calF' \xrightarrow{f'} \calG' \to \calH' \xrightarrow{+1}
\]
in $\calD \QCoh(X,\calO_X)$. Then, again by usual properties of
triangulated categories, there exists an isomorphism
$\For^{\calO_X}(\calH') \cong \calH$ in $\calD(X,\calO_X)$. It follows
that $\calH$ has quasi-coherent cohomology.
\end{proof}

\subsection{K-flats and inverse image}

Let $(X,\calA_X)$ be a dg-scheme. As $\calA_X$ is graded-commutative, we have an equivalence between left and right $\calA_X$-dg-modules (see \cite{BL}). In particular we can take tensor products of two left $\calA_X$-dg-modules, and we still obtain a left $\calA_X$-dg-module. Also, the tensor product is commutative.

Recall the definition of a K-flat object (see \cite[Definition 5.1]{Sp}, \cite[Definition 1.3.1]{R2}).

\begin{defin}

An object $\calF$ of $\calC(X,\calA_X)$ is said to be \emph{K-flat} if for any acyclic object $\calG$ in $\calC(X,\calA_X)$, the dg-module $\calF \otimes_{\calA_X} \calG$ is acyclic.

\end{defin}

Under our hypotheses, every quasi-coherent $\calO_X$-module is the quotient of a quasi-coherent flat $\calO_X$-module. (This follows e.g.~from \cite[Proof of Proposition 1.1]{AJL}.) Statement \eqref{it:prop-existenceK-flat-1} of the following proposition is proved in \cite[Theorem 1.3.3]{R2}; statement \eqref{it:prop-existenceK-flat-2} can be proved similarly using this remark.

\begin{prop} \label{prop:existenceK-flat}

\begin{enumerate}

\item
\label{it:prop-existenceK-flat-1}
For any $\calF$ in $\calC(X,\calA_X)$ there exists a K-flat $\calA_X$-dg-module $\calP$ and a quasi-isomorphism $\calP \xrightarrow{\qis} \calF$.
\item 
\label{it:prop-existenceK-flat-2}
For any $\calF$ in $\calC \QCoh(X,\calA_X)$ there exists a
quasi-coherent $\calA_X$-dg-module $\calP$ such that $\For^{\calA_X}(\calP)$ is
K-flat in the category $\calC(X,\calA_X)$, and a quasi-isomorphism
$\calP \xrightarrow{\qis} \calF$.\qed

\end{enumerate}

\end{prop}

It follows from this proposition that if $f: (X,\calA_X) \to (Y,\calA_Y)$ is a morphism of dg-schemes, the left derived functors
\begin{align*} Lf^* : \calD(Y,\calA_Y) \, & \to \, \calD(X,\calA_X), \\ Lf_{\qc}^* : \calD \QCoh(Y,\calA_Y) \, & \to \, \calD \QCoh(X,\calA_X) \end{align*}
are well defined and can be computed by means of left K-flat resolutions, and that the following diagram commutes:
\begin{equation}
\label{eq:diaginverseimage}
\vcenter{
\xymatrix{ \calD
  \QCoh(Y,\calA_Y) \ar[rr]^{Lf_{\qc}^*} \ar[d]_{\For^{\calA_Y}} & & \calD
  \QCoh(X,\calA_X) \ar[d]^{\For^{\calA_X}} \\ \calD(Y,\calA_Y) \ar[rr]^{Lf^*} & &
  \calD(X,\calA_X).
}
}
\end{equation} Because of this compatibility, we will not write the subscript ``$\qc$'' if no confusion can arise.

\begin{cor} \label{cor:compositioninverseimage}

Let $f : (X,\calA_X) \to (Y,\calA_Y)$ and $g : (Y,\calA_Y) \to (Z,\calA_Z)$ be morphisms of dg-schemes. Then we have isomorphisms \[ L(g \circ f)^* \cong Lf^* \circ Lg^* \quad \text{and} \quad L(g \circ f)_{\qc}^* \cong Lf_{\qc}^* \circ Lg_{\qc}^*. \]

\end{cor}

\begin{proof} The first isomorphism is proved in \cite[Proposition 1.5.4]{R2}. (It follows from the fact that $g^*$ sends K-flats to K-flats.) The second one follows using diagram \eqref{eq:diaginverseimage}, or can be proved similarly.\end{proof}

Note that the natural diagram \[ \xymatrix{ \calD
  \QCoh(Y,\calA_Y) \ar[rr]^{Lf^*} \ar[d]_{\For} & & \calD
  \QCoh(X,\calA_X) \ar[d]^{\For} \\ \calD \QCoh(Y,\calO_Y)
  \ar[rr]^{L(f_0)^*} & & \calD \QCoh(X,\calO_X) } \] is not commutative
in general. It is commutative in certain cases, however, e.g.~if $(f_0)^* \calA_Y \cong
\calA_X$ and $\calA_Y$ is K-flat as an $\calO_Y$-dg-module.

\subsection{K-injectives and direct image}

Let again $(X,\calA)$ be a dg-scheme. First, let us prove a
compatibility result for quasi-coherator functors attached to
$\calO_X$ and to $\calA$.

\begin{lem} \label{lem:compatibilityQ}

The following diagram commutes, where the vertical arrows are
natural forgetful functors: \[ \xymatrix{ \calD(X,\calA)
  \ar[rr]^-{RQ^{\calA}} \ar[d]_-{\For} & & \calD \QCoh(X,\calA)
  \ar[d]^-{\For} \\ \calD(X,\calO_X) \ar[rr]^-{RQ^{\calO_X}} & & \calD
  \QCoh(X,\calO_X). } \]

\end{lem}

\begin{proof} The corresponding diagram for non-derived categories
  (and non-deri\-ved functors) clearly commutes. Hence it is sufficient
  to prove that if $\calI$ is K-injective
  in $\calC(X,\calA)$, then $\For(\calI)$ is split on the right for
  the functor $Q^{\calO_X}$. Recall that an object $\calF$ of $\calC(X,\calO_X)$ is said to be \emph{weakly K-injective} if for any acyclic, K-flat object $\calG$ of $\calC(X,\calO_X)$, the complex of abelian groups $\Hom_{\calO_X}(\calG,\calF)$ is acyclic (see \cite[Definition 5.11]{Sp}). In particular, if $\calI$ is K-injective in $\calC(X,\calA)$, then $\For(\calI)$ is weakly K-injective in $\calC(X,\calO_X)$. (This follows from \cite[Proposition 5.15(b)]{Sp} applied to the natural morphism of dg-ringed spaces $(X,\calA) \to (X,\calO_X)$.) Hence it is enough to prove that weakly K-injective objects of $\calC(X,\calO_X)$ are split on the right for the functor $Q^{\calO_X}$. And for this, using the existence of right K-injective (in particular, weakly K-injective) resolutions (see \cite[Theorem 4.5]{Sp}), it is enough to prove that if $\calJ$ is weakly K-injective and acyclic, then $Q^{\calO_X}(\calJ)$ is acyclic. 

So, let $\calJ$ be such an object. By Proposition \ref{prop:existenceK-flat}, there exists an object $\calP$ of $\calC \QCoh(X,\calO_X)$ such that $\For^{\calO_X}(\calP)$ is K-flat in $\calC(X,\calO_X)$, and a quasi-isomorphism $\calP \xrightarrow{\qis} Q^{\calO_X}(\calJ)$. Let us denote by $\calH(X,\calO_X)$, respectively $\calH \QCoh(X,\calO_X)$, the homotopy category of $\calC(X,\calO_X)$, respectively $\calC \QCoh(X,\calO_X)$. By adjunction, there is an isomorphism \[ \Hom_{\calH \QCoh(X,\calO_X)}(\calP, Q^{\calO_X}(\calJ)) \, \cong \, \Hom_{\calH(X,\calO_X)}(\For^{\calO_X} (\calP), \calJ). \] The complex of abelian groups $\Hom_{\calO_X}(\For^{\calO_X} (\calP), \calJ)$ is acyclic by \cite[Proposition 5.20]{Sp}. Taking the $0$-th cohomology of this complex, we deduce $\Hom_{\calH(X,\calO_X)}(\For^{\calO_X} (\calP), \calJ) = 0$, hence $\Hom_{\calH \QCoh(X,\calO_X)}(\calP, Q^{\calO_X}(\calJ))=0$. In particular, the quasi-isomor\-phism $\calP \to Q^{\calO_X}(\calJ)$ is homotopic to zero, which implies that $Q^{\calO_X}(\calJ)$ is acyclic. \end{proof}

\begin{prop}
\label{prop:equivdgalg}

\begin{enumerate}
\item 
\label{it:prop-equivdgalg-0}
The subcategory $\calD^{\qc}(X,\calA)$ is a triangulated
\- subcategory of $\calD(X,\calA)$.
\item
\label{it:prop-equivdgalg-1}
The functors
\begin{align*} 
\For^{\calA} : & \ \calD \QCoh(X,\calA)
  \to \calD^{\qc}(X,\calA), \\
RQ^{\calA} : & \ \calD^{\qc}(X,\calA)
  \to \calD \QCoh(X,\calA)
\end{align*}
are quasi-inverse equivalences of triangulated categories.
\item 
\label{it:prop-equivdgalg-2}
There are enough K-injectives in the category $\calC
\QCoh(X,\calA)$.

\end{enumerate}

\end{prop}

\begin{proof} 
\eqref{it:prop-equivdgalg-0} This is an easy consequence of Proposition \ref{prop:thmSGA}\eqref{it:prop-thmSGA-2}.

\eqref{it:prop-equivdgalg-1} We have already seen that these functors are
  adjoint. Hence there are natural adjunction morphisms $\For^{\calA}
  \circ RQ^{\calA} \to \Id$ and $\Id \to RQ^{\calA} \circ
  \For^{\calA}$. By Lemma \ref{lem:compatibilityQ}, the following diagram
  commutes: \[ \xymatrix{ \calD^{\qc}(X,\calA)
    \ar@<0.5ex>[rr]^-{RQ^{\calA}} \ar[d]_-{\For} & & \calD \QCoh(X,\calA)
    \ar@<0.5ex>[ll]^-{\For^{\calA}} \ar[d]^-{\For} \\
    \calD^{\qc}(X,\calO_X) \ar@<0.5ex>[rr]^-{RQ^{\calO_X}} & & \calD
    \QCoh(X,\calO_X). \ar@<0.5ex>[ll]^-{\For^{\calO_X}} } \] Hence it follows from Proposition \ref{prop:thmSGA} that for every $\calF$ in $\calD
  \QCoh(X,\calA)$ and $\calG$ in $\calD^{\qc}(X,\calA)$ the
  induced morphisms $\For^{\calA} \circ RQ^{\calA}(\calG) \to \calG$
  and $\calF \to RQ^{\calA} \circ \For^{\calA}(\calF)$ are
  isomorphisms. Statement \eqref{it:prop-equivdgalg-1} follows. 

\eqref{it:prop-equivdgalg-2} We have seen in \S \ref{ss:definitionsdg} that $Q^{\calA}$ sends K-injectives to
K-injectives, and that there are enough K-injectives in the category
$\calC(X,\calA)$. Let $\calF$ be an object of $\calC \QCoh(X,\calA)$,
and let $\calI$ be a K-injective resolution of $\For^{\calA}(\calF)$ in
$\calC(X,\calA)$. Then in $\calD \QCoh(X,\calA)$ we have $\calF \cong
RQ^{\calA} \circ \For^{\calA}(\calF) \cong Q^{\calA}(\calI)$. Moreover, as the dg-module $Q^{\calA}(\calI)$ is K-injective, we can represent this isomorphism
as a quasi-isomorphism $\calF \xrightarrow{\qis} Q^{\calA}(\calI)$ in
$\calC \QCoh(X,\calA)$. This proves \eqref{it:prop-equivdgalg-2}.
\end{proof}

\begin{remark}
It follows from Proposition \ref{prop:equivdgalg}\eqref{it:prop-equivdgalg-1} and diagram \eqref{eq:diaginverseimage} that we have $Lf^* \bigl( \calD^{\qc}(Y,\calA_Y) \bigr) \subset \calD^{\qc}(X,\calA_X)$.
\end{remark}

It follows from statement \eqref{it:prop-equivdgalg-2} of Proposition \ref{prop:equivdgalg} that if $f: (X,\calA_X) \to (Y,\calA_Y)$
is a morphism of dg-schemes, we
can consider the right derived functors \begin{align*} Rf_* :
  \calD(X,\calA_X) \, & \to \, \calD(Y,\calA_Y), \\ Rf^{\qc}_* : \calD
\QCoh(X,\calA_X) \, & \to \, \calD \QCoh(Y,\calA_Y). \end{align*} (Note that, by \cite[Proposition II.5.8(c)]{HARAG}, under our assumptions the direct image of a quasi-coherent $\calO_X$-module is a quasi-coherent $\calO_Y$-module.) As
usual, we denote by
$f_0 : (X,\calO_X) \to (Y,\calO_Y)$ the associated morphism of
ordinary schemes. We also have the associated derived functors
\begin{align*} R(f_0)_* : \calD(X,\calO_X) \, & \to \, \calD(Y,\calO_Y),
  \\ R(f_0)^{\qc}_* : \calD \QCoh(X,\calO_X) \, & \to \, \calD
  \QCoh(Y,\calO_Y). \end{align*} 

\begin{lem} \label{lem:directimage}

\begin{enumerate}
\item 
\label{it:lem-directimage-1}
We have $R(f_0)_* (
\calD^{\qc}(X,\calO_X) ) \subset \calD^{\qc}(Y,\calO_Y)$.
\item 
\label{it:lem-directimage-2}
The following diagram commutes:
\[
\xymatrix{
\calD \QCoh(X,\calO_X) \ar[rr]^-{R(f_0)^{\qc}_*} \ar[d]_-{\For^{\calO_X}} & & \calD \QCoh(Y,\calO_Y) \ar[d]^-{\For^{\calO_Y}} \\
\calD(X,\calO_X) \ar[rr]^-{R(f_0)_*} & & \calD(Y,\calO_Y).
}
\]

\end{enumerate}

\end{lem}

\begin{proof}
\eqref{it:lem-directimage-1} is proved in \cite[Proposition II.2.1]{H}. (See \cite[Proposition 3.9.2]{Li} for a proof under less restrictive hypotheses.) Let us deduce \eqref{it:lem-directimage-2}. By Proposition \ref{prop:equivdgalg}\eqref{it:prop-equivdgalg-1}, it is enough to prove that
\[ R(f_0)^{\qc}_* \circ RQ^{\calO_X} \cong RQ^{\calO_Y} \circ R(f_0)_*. \] It is known that
\begin{equation}
\label{isomdirectimage}
(f_0)^{\qc}_* \circ Q^{\calO_X} \cong
  Q^{\calO_Y} \circ (f_0)_*
\end{equation}
(see \cite[p.~188, Lemma
3.4]{SGA6}). Moreover, $Q^{\calO_X}$ sends K-injectives to K-injectives, and
$(f_0)_*$ sends K-injectives to weakly K-injectives (see
\cite[Proposition 5.15]{Sp}), which are split on the right for $Q^{\calO_Y}$ (see the
proof of Lemma \ref{lem:compatibilityQ}). Hence the result
follows from isomorphism (\ref{isomdirectimage}).
\end{proof}

Now we deduce analogues of these properties for $\calA$-dg-modules.

\begin{prop} \label{prop:compatibilitydirectimage}

\begin{enumerate}
\item 
\label{it:prop-compatibilitydirectimage-1}
 We have $R(f)_* ( \calD^{\qc}(X,\calA_X) ) \subset
\calD^{\qc}(Y,\calA_Y)$.
\item 
\label{it:prop-compatibilitydirectimage-2}
The following diagram
commutes:
\[
\xymatrix{
\calD \QCoh(X,\calA_X) \ar[rr]^-{Rf^{\qc}_*} \ar[d]_-{\For^{\calA_X}} & & \calD \QCoh(Y,\calA_Y)
  \ar[d]^-{\For^{\calA_Y}} \\
\calD(X,\calA_X) \ar[rr]^-{Rf_*} & &
  \calD(Y,\calA_Y).
}
\]
\item 
\label{it:prop-compatibilitydirectimage-3}
The following diagram also commutes:
\[
\xymatrix{
\calD \QCoh(X,\calA_X)
  \ar[rr]^-{Rf^{\qc}_*} \ar[d]_-{\For_X} & & \calD \QCoh(Y,\calA_Y)
  \ar[d]^-{\For_Y} \\
\calD \QCoh (X,\calO_X) \ar[rr]^-{R(f_0)^{\qc}_*} & &
  \calD \QCoh(Y,\calO_Y).
}
\]

\end{enumerate}

\end{prop}

\begin{proof} Recall that the following diagram commutes:
\begin{equation}
\label{eq:diagdirectimage}
\vcenter{
\xymatrix{
    \calD(X,\calA_X)
    \ar[rr]^-{Rf_*} \ar[d]_-{\For_X} & & \calD(Y,\calA_Y)
    \ar[d]^-{\For_Y} \\
\calD(X,\calO_X) \ar[rr]^-{R(f_0)_*} & &
    \calD(Y,\calO_Y),
}
}
\end{equation}
see \cite[Corollary 1.5.3]{R2}. Then \eqref{it:prop-compatibilitydirectimage-1} follows from Lemma \ref{lem:directimage}\eqref{it:lem-directimage-1}. 

Statement \eqref{it:prop-compatibilitydirectimage-2} can be proved similarly to Lemma \ref{lem:directimage}\eqref{it:lem-directimage-2}. For this we use the fact that a weakly
  K-injective $\calA_Y$-dg-module is also weakly K-injective as an
  $\calO_Y$-dg-module (see \cite[Proposition 5.15(b)]{Sp}) hence it is split on the right for the functor $Q^{\calO_Y}$ (see the proof of Lemma
  \ref{lem:compatibilityQ}), which implies that it is also split on the right for the functor $Q^{\calA_Y}$. 

Finally, consider the diagram of \eqref{it:prop-compatibilitydirectimage-3}. By definition there is a natural morphism of functors $\For_Y \circ Rf^{\qc}_* \to R(f_0)^{\qc}_* \circ \For_X$. The fact that it is an isomorphism follows from diagram \eqref{eq:diagdirectimage}, statement \eqref{it:prop-compatibilitydirectimage-2}, and Lemma \ref{lem:directimage}\eqref{it:lem-directimage-2}. \end{proof}

Because of these compatibility results, we will not write the
supscript ``$\qc$'' if no confusion can arise.

\begin{cor} \label{cor:compositiondirectimage}

Let $f : (X,\calA_X) \to (Y,\calA_Y)$ and $g : (Y,\calA_Y) \to
(Z,\calA_Z)$ be morphisms of dg-schemes. Then we have isomorphisms
\[ R(g \circ f)_* \cong Rg_* \circ Rf_* \quad \text{and} \quad R(g
\circ f)^{\qc}_* \cong Rg^{\qc}_* \circ Rf^{\qc}_*.\]

\end{cor}

\begin{proof} The first isomorphism is proved in \cite[Corollary 1.5.3]{R2}. The second one follows, using Proposition \ref{prop:compatibilitydirectimage}\eqref{it:prop-compatibilitydirectimage-2}. \end{proof}

\subsection{Adjunction}

\begin{prop}

\begin{enumerate}
\item 
\label{it:adjunction-1}
The functors
\begin{align*}
Lf^* : \calD(Y,\calA_Y) \, & \to \, \calD(X,\calA_X), \\ Rf_* : \calD(X,\calA_X) \, & \to \, \calD(Y,\calA_Y)
\end{align*}
are adjoint.
\item 
\label{it:adjunction-2}
Similarly, the functors
\begin{align*}
Lf_{\qc}^* : \calD \QCoh(Y,\calA_Y) \, & \to \, \calD \QCoh(X,\calA_X), \\
Rf^{\qc}_* : \calD \QCoh(X,\calA_X) \, & \to \, \calD \QCoh(Y,\calA_Y) 
\end{align*}
are adjoint.

\end{enumerate}

\end{prop}

\begin{proof} These results follow from general properties of derived functors, see \cite[Lemma 13.6]{Ke}. But they can also be proved directly, as follows. First, \eqref{it:adjunction-1} is proved (by direct methods) in \cite[Theorem 1.6.2]{R2}. Let us
  deduce \eqref{it:adjunction-2}. Let $\calF$ be an object of $\calD
  \QCoh(Y,\calA_Y)$, and $\calG$ an object of $\calD
  \QCoh(X,\calA_X)$. First we have:
\begin{align*}
\Hom_{\calD \QCoh(X,\calA_X)} (Lf_{\qc}^* \calF, \calG) \, & \cong \, \Hom_{\calD(X,\calA_X)} ( \For^{\calA_X} \circ Lf_{\qc}^* \calF, \For^{\calA_X} \, \calG) \\
& \cong \, \Hom_{\calD(X,\calA_X)} ( Lf^* \circ \For^{\calA_Y} \, \calF, \For^{\calA_X} \, \calG).
\end{align*}
Here the first isomorphism follows from Proposition \ref{prop:equivdgalg}\eqref{it:prop-equivdgalg-1}, and the second one
  from diagram \eqref{eq:diaginverseimage}. Hence, by
  \eqref{it:adjunction-1} we have
\[
\Hom_{\calD \QCoh(X,\calA_X)}
    (Lf_{\qc}^* \calF, \calG) \, \cong \, \Hom_{\calD(Y,\calA_Y)} (\For^{\calA_Y} \, \calF, Rf_* \circ \For^{\calA_X} \, \calG).
\]
Then, using
  Proposition \ref{prop:compatibilitydirectimage}\eqref{it:prop-compatibilitydirectimage-2} and again
  Proposition \ref{prop:equivdgalg}\eqref{it:prop-equivdgalg-1}, we deduce 
\begin{align*}
\Hom_{\calD \QCoh(X,\calA_X)} (Lf_{\qc}^* \calF, \calG) \, & \cong \,
    \Hom_{\calD(Y,\calA_Y)} (\For^{\calA_Y} \, \calF, \For^{\calA_Y} \circ Rf^{\qc}_* \calG) \\ & \cong \, \Hom_{\calD \QCoh(Y,\calA_Y)} (\calF,
    Rf^{\qc}_* \calG).
\end{align*}
These isomorphisms are functorial.
\end{proof}

\subsection{Projection formula}

In this subsection we generalize the classical projection formula (see \cite[Proposition II.5.6]{H}) to dg-schemes.

\begin{lem} \label{lem:projformulaf0}

\begin{enumerate}
\item 
\label{it:lem-projformulaf0-1}
The functors $R(f_0)_*$ and $Rf_*$ commute with filtered direct limits.
\item 
\label{it:lem-projformulaf0-2}
For $\calF$ in $\calD(X,\calO_X)$ and $\calG$ in
$\calD^{\qc}(X,\calO_X)$ we have a functorial isomorphism
\[
R(f_0)_*(\calF \, \lotimes_{\calO_X} \, L(f_0)^* \calG) \ \cong \
(R(f_0)_* \calF) \, \lotimes_{\calO_Y} \, \calG.
\]

\end{enumerate}

\end{lem}

\begin{proof} 
\eqref{it:lem-projformulaf0-1} The case of $R(f_0)_*$ can be proved as in \cite[Corollary 1.7.5]{R2}. (In \emph{loc.}~\emph{cit.}, ``direct sum'' can be replaced by ``filtered direct limit'' without any trouble.) Then the case of $Rf_*$ follows, using diagram \eqref{eq:diagdirectimage}.

Assertion \eqref{it:lem-projformulaf0-2} is proved in \cite[Proposition 3.9.4]{Li}.
\end{proof}

\begin{prop}[Projection formula] \label{prop:projectionformula}

For $\calF \in \calD(X,\calA_X)$ and $\calG \in
\calD^{\qc}(Y,\calA_Y)$ we have a functorial isomorphism \[Rf_* (\calF
\, \lotimes_{\calA_X} \, Lf^* \calG) \ \cong \ (Rf_* \calF) \,
\lotimes_{\calA_Y} \, \calG.\]

\end{prop}

\begin{proof} By the same arguments as in \cite[Proposition 5.6]{H} or in \cite[\S 3.4.6]{Li}, there
is a morphism of functors \[ (Rf_* \calF) \, \lotimes_{\calA_Y} \, \calG \ \to \
  Rf_* (\calF \, \lotimes_{\calA_X} \, Lf^* \calG).\] We want to prove
  that it is an isomorphism. First we can assume that $\calG$ is the
  image of an object of $\calC \QCoh(Y,\calA_Y)$ (see Proposition
  \ref{prop:equivdgalg}\eqref{it:prop-equivdgalg-1}). Then, by the proof of
  Proposition \ref{prop:existenceK-flat}\eqref{it:prop-existenceK-flat-2}, $\calG$ is quasi-isomorphic to
  the direct limit of
  some $\calA_Y$-dg-modules $\calP_p$ such that, for any $p$,
  $\calP_p$ has a finite filtration with subquotients of the form
  $\calA_Y \otimes_{\calO_Y} \calG_0$ with $\calG_0$ in
  $\calC \QCoh(Y,\calO_Y)$, K-flat in $\calC(Y,\calO_Y)$. Hence, by Lemma \ref{lem:projformulaf0}\eqref{it:lem-projformulaf0-1}, we can
  assume that $\calG = \calA_Y \otimes_{\calO_Y} \calG_0$ for such a
  $\calG_0$. Then $\calG$ is K-flat, we have in $\calD(X,\calO_X)$ \[ Lf^* \calG \ \cong \ f^* \calG \ \cong \ \calA_X \otimes_{\calO_X} (f_0)^* \calG_0,\] and the
  $\calO_X$-dg-module $(f_0)^* \calG_0$ is K-flat. Hence \[\calF \,
  \lotimes_{\calA_X} \, Lf^* \calG \ \cong \ \calF \otimes_{\calO_X} (f_0)^*
  \calG_0 \ \cong \ \calF \, \lotimes_{\calO_X} \, (f_0)^* \calG_0.\] Similarly
  we have
\[
(Rf_* \calF) \, \lotimes_{\calA_Y} \, \calG \ \cong \ (Rf_* \calF)
  \, \lotimes_{\calO_Y} \, \calG_0.
\]
Hence the result follows from Lemma \ref{lem:projformulaf0}\eqref{it:lem-projformulaf0-2} and the compatibility between $Rf_*$ and
  $R(f_0)_*$, see diagram \eqref{eq:diagdirectimage}.
\end{proof}

\subsection{Quasi-isomorphisms}

\begin{prop} \label{prop:qis}

Let $f : (X,\calA_X) \to (Y,\calA_Y)$ be a morphism of dg-schemes such
that $f_0$ is a closed embedding, and the induced morphism $\calA_Y
\to (f_0)_* \calA_X$ is a quasi-isomorphism of dg-algebras. Then the
functors \begin{align*} Rf_* : \calD(X,\calA_X) \, & \to \,
  \calD(Y,\calA_Y), \\ Lf^* : \calD(Y,\calA_Y) \, & \to \,
  \calD(X,\calA_X) \end{align*} are quasi-inverse equivalences of
triangulated categories. They induce equivalences of triangulated categories \[ \calD^{\qc}(X,\calA_X) \, \cong \, \calD^{\qc}(Y,\calA_Y), \quad \calD \QCoh(X,\calA_X) \, \cong \, \calD \QCoh(Y,\calA_Y). \]

\end{prop}

\begin{proof} One can factor $f$ as the composition \[ (X,\calA_X)
  \xrightarrow{f_1} (Y,(f_0)_* \calA_X) \xrightarrow{f_2} (Y,
  \calA_Y). \] Hence, using Corollaries \ref{cor:compositiondirectimage}
  and \ref{cor:compositioninverseimage}, it is sufficient to prove the
  result for $f_1$ and
  $f_2$. The case of $f_2$ is treated in \cite[Proposition 1.5.6]{R2}. (The
  proof is similar to the usual case of a dg-algebra, see
  \cite{BL}.) And for $f_1$, the (non-derived) functors
  $(f_1)_*$ and $(f_1)^*$ are already (exact) equivalences. \end{proof}

\subsection{Derived fiber product and base change} \label{ss:basechange}

In this subsection we generalize the usual flat base change theorem (\cite[Proposition II.5.12]{H}, \cite[Theorem 3.10.3]{Li}) to dg-schemes. One of the main advantages of considering dg-schemes is that, in this generality, one can replace fiber products by \emph{derived} fiber products, and then get rid of the flatness (or ``independence'') assumption.

First, a morphism of dg-schemes $f: (X,\calA_X) \to (Y,\calA_Y)$ is said to be \emph{smooth} if the underlying morphism of schemes $f_0 : X \to Y$ is smooth, and $\calA_X$ is K-flat over $(f_0)^* \calA_Y$.

Let now $\bX = (X,\calA_X)$, $\bY=(Y,\calA_Y)$, $\bZ=(Z,\calA_Z)$ be dg-schemes and let $f: \bX \to \bZ$, $g: \bY \to \bZ$ be morphisms. As in the case of ordinary schemes (\cite[p.~87]{HARAG}), one can easily define the fiber product dg-scheme \[ \bX \times_{\bZ} \bY. \] (If $X_0$, $Y_0$ and $Z_0$ are affine, the fiber product is given by the tensor product of dg-algebras.) Assume now that the morphisms $f_0 : X_0 \to Z_0$ and $g_0 : Y_0 \to Z_0$ are quasi-projective\footnote{More generally, for the definition of the derived fiber product one only has to assume that $f$ factors through a closed embedding $X_0 \hookrightarrow X_0'$ where the morphism $X_0' \to Z_0$ is smooth, and similarly for $g$. In the proof of Proposition \ref{prop:basechange}, we will also need the assumption that $f_0$ is of finite type.}. Then one can factor these morphisms as compositions \[ \bX \overset{f_1}{\hookrightarrow} \bX ' \xrightarrow{f_2} \bZ, \qquad \bY \overset{g_1}{\hookrightarrow} \bY ' \xrightarrow{g_2} \bZ,\] where $f_1$ and $g_1$ are quasi-isomorphic closed embeddings (i.e.~satisfy the assumptions of Proposition \ref{prop:qis}), and $f_2$ and $g_2$ are smooth. (See \cite[Theorem 2.7.6]{CK} for the existence of such factorizations.) Then one can ``define'' the dg-scheme \[ \bX \, \Rtimes_{\bZ} \, \bY \] as $\bX ' \times_{\bZ} \bY$, or equivalently $\bX \times_{\bZ} \bY '$, or equivalently $\bX ' \times_{\bZ} \bY '$. More precisely, this dg-scheme is defined only ``up to quasi-isomorphism.'' However, the categories $\calD(\bX \, \Rtimes_{\bZ} \, \bY)$ and $\calD \QCoh(\bX \, \Rtimes_{\bZ} \, \bY)$ are well defined, thanks to Proposition \ref{prop:qis}. These are the only objects we are going to use. To give a more precise definition of $\bX \, \Rtimes_{\bZ} \, \bY$, one would have to consider a ``derived category of dg-schemes'' as in \cite[\S 2.2]{CK}; we will not do this here.

There are natural projections $p_1 : \bX \, \Rtimes_{\bZ} \, \bY \to \bX$, $p_2 : \bX \, \Rtimes_{\bZ} \, \bY \to \bY$, represented by the morphisms of dg-schemes $\bX \times_{\bZ} \bY ' \to \bX$, $\bX ' \times_{\bZ} \bY \to \bY$. Note that the following diagram commutes, where $i : \bX \times_{\bZ} \bY ' \to \bX ' \times_{\bZ} \bY '$ is the morphism induced by $f_1$, and $p_1':\bX ' \times_{\bZ} \bY ' \to \bX '$ is induced by $g_2$: \[ \xymatrix{ \calD(\bX \times_{\bZ} \bY ') \ar[rr]^{Ri_*} \ar[d]_{R(p_1)_*} & & \calD(\bX ' \times_{\bZ} \bY ') \ar[d]^{R(p_1')_*} \\ \calD(\bX) \ar[rr]^{R(f_1)_*} & & \calD(\bX ') } \] (see Corollary \ref{cor:compositiondirectimage}). As $R(f_1)_*$ and $Ri_*$ are equivalences of categories with respective inverses $L(f_1)^*$ and $Li^*$ (by Proposition \ref{prop:qis}), we also have $R(p_1)_* \cong L(f_1)^* \circ R(p_1')_* \circ Li^*$. Similar results apply for the inverse image functors, and for $\bY$ instead of $\bX$. It follows from these remarks that, even if $p_1$ and $p_2$ are ``not well defined'' as morphisms of dg-schemes (because their source is not well defined), the associated direct and inverse image functors are well defined (i.e.~are compatible with the natural equivalences between the categories of dg-sheaves associated to the different realizations of $\bX \, \Rtimes_{\bZ} \, \bY$ as a dg-scheme).

\begin{prop}[Base change theorem] \label{prop:basechange}

Consider the diagram \[ \xymatrix{ \bX \, \Rtimes_{\bZ} \, \bY \ar[r]^-{p_2}
  \ar[d]_-{p_1} & \bY \ar[d]^-{g} \\ \bX \ar[r]^-{f} & \bZ. } \] Then
for $\calF$ in $\calD \QCoh(\bX)$ we have a functorial isomorphism
\[ Lg^* \circ Rf_* \, \calF \ \cong \ R(p_2)_* \circ L(p_1)^* \, \calF. \]

\end{prop}

\begin{proof} As usual, there is a morphism of functors from the left hand side to the right hand side of the isomorphism we are trying to prove. (See \cite[Proposition 3.7.2]{Li} for the similar statement for ordinary schemes.) What we have to show is that it is an isomorphism.

Using resolutions as above, one can assume that $f$ and $g$ are smooth. In this case, $\bX \, \Rtimes_{\bZ} \, \bY$ is simply $\bX \times_{\bZ} \bY$. Let $\widetilde{g}: (Y, (g_0)^* \calA_Z) \to (Z,\calA_Z)$, $\widetilde{p_1} : (X \times_Z Y, (p_{1,0})^* \calA_X) \to (X,\calA_X)$ and $\widetilde{p_2} : (X \times_Z Y, (p_{1,0})^* \calA_X) \to (Y, (g_0)^* \calA_Z)$ be the natural morphisms of dg-ringed spaces. By definition, and the ordinary flat base change theorem (\cite[Proposition II.5.12]{H} or \cite[Proposition 3.9.5]{Li}) we have  \begin{align*} Lg^* \circ Rf_* \, \calF \ & \cong \ \calA_Y
    \otimes_{(g_0)^* \calA_Z} (\widetilde{g})^* \circ Rf_* \, \calF \\ & \cong \
    \calA_Y \otimes_{(g_0)^* \calA_Z} R(\widetilde{p_2})_* \circ (\widetilde{p_1})^* \,
    \calF. \end{align*} On the other hand, by definition of the fiber
  product we have \[ \calA_{X \times_Z Y} \cong (p_{1,0})^* \calA_X
  \otimes_{(p_{1,0})^* (f_0)^* \calA_Z} (p_{2,0})^* \calA_Y. \] Hence, by the projection formula (Proposition \ref{prop:projectionformula}),
  \begin{align*} R(p_2)_* \circ L(p_1)^* \, \calF \ & \cong \ R(p_2)_*
    (\calA_{X \times_Z Y} \otimes_{(p_{1,0})^* \calA_X} (\widetilde{p_1})^* \, \calF)
    \\ & \cong \ \calA_Y \otimes_{(g_0)^* \calA_Z} R(\widetilde{p_2})_* \circ
    (\widetilde{p_1})^* \, \calF. \end{align*} This concludes the
  proof. \end{proof}

\subsection{Compatibility of projection and base change} \label{ss:compatibility}

To finish this section, we observe that one can prove some compatibility result for the isomorphisms of Propositions \ref{prop:projectionformula} and \ref{prop:basechange}. This is similar to \cite[Proposition 3.7.3]{Li}, and left to the interested reader.

\section{Convolution and geometric actions} \label{sec:convolution}

In this section we present the formalism of functors on derived
categories of coherent sheaves arising from ``integral kernels.'' We fix two dg-schemes $\bX=(X,\calA_X)$ and $\bY=(Y,\calA_Y)$, and consider a morphism of dg-schemes $f:\bX \to \bY$ such that $f_0$ is quasi-projective.

\subsection{Convolution} \label{ss:convolution}

In this subsection we will consider the category \[ \calK_{\bX,\bY} \, := \, \calD \QCoh(\bX \, \Rtimes_{\bY} \, \bX). \] An important particular case is when $\bX=X$ and $\bY=Y$ are ordinary schemes, such that
\begin{equation}\label{Torvan}
\Tor_{\neq 0}^{f^{-1} \calO_Y}(\calO_X,\calO_X)=0.
\end{equation}
Then the dg-scheme $X \, \Rtimes_{Y} \, X$ ``is'' the ordinary scheme $X\times_Y X$, and $\calK_{X,Y} \cong \calD\QCoh (X \times_Y X)$.

\begin{prop} \label{prop:defconvolution}

There is a natural convolution product on $\calK_{\bX,\bY}$, which endows this category with a monoidal structure.

\end{prop}

\begin{proof}[Sketch of proof] 
We only give the definition of the product. Its main properties are easy to check using the projection formula (Proposition \ref{prop:projectionformula}), the base change theorem (Proposition \ref{prop:basechange}) and their compatibility (\S \ref{ss:compatibility}). A similar construction has been considered in \cite{MRHec} in a special case.

For simplicity, using the constructions of \S \ref{ss:basechange}, one can assume that the morphism $\bX \to \bY$ is smooth, so that all the derived fiber products become ordinary fiber products. For $(i,j) \in \{(1,2), (1,3), (2,3) \}$, we denote by $q_{i,j} : X_0 \times_{Y_0} X_0 \times_{Y_0} X_0 \to X_0 \times_{Y_0} X_0$ the projection on the $i$-th and $j$-th factors, and by $\bZ_{i,j}$ the following dg-scheme: \[ \bZ_{i,j}\ := \ \bigl( X_0 \times_{Y_0} X_0 \times_{Y_0} X_0, \ (q_{i,j})^* \calA_{\bX \times_{\bY} \bX} \bigr). \] There is a natural morphism of dg-schemes $p_{i,j} : \bZ_{i,j} \to \bX \times_{\bY} \bX$, and associated functors $R(p_{i,j})_*$, $L(p_{i,j})^*$. (In the general case, when $\bX \to \bY$ is not assumed to be smooth, the dg-scheme $\bZ_{i,j}$ can be defined as the derived fiber product \[ \bigl( \bX \, \Rtimes_{\bY} \, \bX \bigr) \, \Rtimes_{X_0 \Rtimes_{Y_0} X_0} \, \bigl( X_0 \, \Rtimes_{Y_0} \, X_0 \, \Rtimes_{Y_0} \, X_0 \bigr), \] where the morphism $X_0 \, \Rtimes_{Y_0} \, X_0 \, \Rtimes_{Y_0} \, X_0 \to X_0 \, \Rtimes_{Y_0} \, X_0$ is the derived version of $q_{i,j}$.)

Let also $q_2 : X_0 \times_{Y_0} X_0 \times_{Y_0} X_0 \to X_0$ be the projection on the second factor, and consider the sheaf of dg-algebras $q_2^* \calA_X$ on $X_0 \times_{Y_0} X_0 \times_{Y_0} X_0$. There is a natural bifunctor \[ (- \otimes_2 -) : \ \left\{ \begin{array}{ccc} \calC \QCoh(\bZ_{1,2}) \times \calC \QCoh(\bZ_{2,3}) & \to & \calC \QCoh(\bZ_{1,3}) \\ (\calF,\calG) & \mapsto & \calF \otimes_{q_2^* \calA_X} \calG. \end{array} \right. \] (Here we forget about the action of the middle copy of $\calA_X$ on $\calF \otimes_{q_2^* \calA_X} \calG$.) This bifunctor has a derived bifunctor between the corresponding derived categories, which we denote by $(- \lotimes_2 -)$. It can be computed using K-flat resolutions.

Then the convolution product can be defined as follows for $\calF$, $\calG$ in $\calK_{\bX,\bY}$: 
\[
\calF \star \calG \ := \ R(p_{1,3})_* \, \bigl( L(p_{1,2})^* \calG \, \lotimes_2 \, L(p_{2,3})^* \calF \bigr).
\]
As explained above, the basic properties of this product can be proved by copying the usual proofs for schemes. In particular, the unit object for this product is the direct image of the structure sheaf under the diagonal embedding $\Delta : \bX \to \bX \, \Rtimes_\bY \, \bX$.
\end{proof}

\begin{remark}
Suppose that $\bX$ and $\bY$ are ordinary schemes, which satisfy the Tor vanishing assumption \eqref{Torvan}.
Then the assertion of Proposition \ref{prop:defconvolution} does not involve dg-schemes. However, its proof (the construction of the monoidal structure on $\calK_{X,Y}$) makes use of  the triple product 
$X \, \Rtimes_Y \, X \, \Rtimes_Y \, X$, which reduces to an ordinary scheme only under a stronger Tor vanishing assumption, which does not hold in the examples of interest to us.
\end{remark}

\subsection{Action by convolution}
\label{ss:actionconvolution}

Our main interest in the category $\calK_{\bX,\bY}$ comes from the following result.

\begin{prop} \label{prop:action}

\begin{enumerate}
\item 
\label{it:prop-action-1}
The category $\calD \QCoh(\bX)$ carries a natural convolution action of the monoidal category $\calK_{\bX,\bY}$.
\item 
\label{it:prop-action-2}
Assume that $\bX$ and $\bY$ are ordinary schemes, that $f$ is proper, and that $X$ is a smooth scheme. Let 
$\calK^{\mathrm{coh}}_{X,Y}\subset \calK_{X,Y}$ be the full subcategory
which consists of complexes with a finite number of nonzero
cohomology sheaves, each of which is a coherent sheaf on $X\times_Y X$. Then $\calK^{\mathrm{coh}}_{X,Y}$ is a monoidal subcategory, and its action preserves the full subcategory $\calD^b \Coh(X) \subset \calD \QCoh(X)$.

\end{enumerate}

\end{prop}

\begin{proof}[Sketch of proof]
We only give the definition of the convolution action, leaving the details to the interested reader. By definition there are morphisms of dg-schemes $p_1,p_2 : \bX \, \Rtimes_\bY \, \bX \to \bX$. For $\calF$ in $\calK_{\bX,\bY}$, the associated functor is given by: \[ \left\{ \begin{array}{ccc} \calD \QCoh(\bX) & \to & \calD \QCoh(\bX) \\ \calG & \mapsto & R(p_2)_* \, \bigl( \calF \, \lotimes_{\bX \Rtimes_\bY \bX} \, L(p_1)^* \calG \bigr) \end{array} \right. , \] where $\lotimes_{\bX \Rtimes_\bY \bX}$ stands for the (derived) tensor product over $\calA_{\bX \Rtimes_\bY \bX}$.
\end{proof}

\begin{remark} We will mainly apply Proposition \ref{prop:action} in the situation of \eqref{it:prop-action-2}. Notice that this action factors through another monoidal category which
is simpler to define, namely the category $\calD^b \Coh_{X\times_Y X}(X\times X)$ which is the full
subcategory in $\calD^b \Coh(X \times X)$ consisting of complexes set-theoretically supported on
$X\times_Y X$; more precisely, the action of statement \eqref{it:prop-action-2} is the composition of the action of the monoidal category  $\calD^b \Coh_{X\times_Y X}(X\times X)$ on $\calD^b \Coh(X)$ and the functor of
``direct image under closed embedding'' $\calK_{X,Y}^{\mathrm{coh}} \to \calD^b \Coh_{X\times_Y X}(X\times
X)$.  The motivation for introducing the more complicated category $\calK_{X,Y}$ 
is the base change construction of Proposition \ref{prop:basechangeaction} below. \end{remark}

\begin{prop}
\label{prop:basechangeaction}

Let $\bY ' \to \bY$ be a morphism of dg-schemes, whose underlying morphism of schemes is quasi-projective, and set $\bX ':=\bX \, \Rtimes_{\bY} \, \bY '$.

\begin{enumerate}
\item 
\label{it:prop-basechangeaction-1}
The pull-back functor $\calK_{\bX,\bY} \to \calK_{\bX',\bY'}$ is monoidal. 
\item 
\label{it:prop-basechangeaction-2}
The pull-back functor $\calD \QCoh(\bX) \to \calD \QCoh(\bX')$ is compatible with the actions of $\calK_{\bX,\bY}$, where the action on $\calD \QCoh(\bX')$ is the composition of the action of $\calK_{\bX',\bY'}$ given by Proposition {\rm \ref{prop:action}}\eqref{it:prop-action-1} and the monoidal functor of \eqref{it:prop-basechangeaction-1}.

\end{enumerate}

\end{prop}

\begin{proof} Observe that there is a natural quasi-isomorphism of dg-schemes \[ \bX' \, \Rtimes_{\bY'} \, \bX' \ \overset{\qis}{\cong} \ (\bX \, \Rtimes_{\bY} \, \bX) \, \Rtimes_{\bY} \, \bY'. \] Similarly, let $\bZ_{i,j}$, $(i,j) \in \{(1,2),(2,3),(1,3) \}$, be defined as in the proof of Proposition \ref{prop:defconvolution}, and let $\bZ_{i,j}'$ be the dg-scheme defined similarly for $\bX', \bY'$ instead of $\bX,\bY$. Then we have
\[
\bZ_{i,j}' \ \overset{\qis}{\cong} \ \bZ_{i,j} \, \Rtimes_{\bY} \, \bY'.
\] 

Using these remarks, the proposition follows from the base change theorem. Indeed, the following diagram is cartesian (in the derived sense):
\[
\xymatrix{ \bZ_{1,3}' \ar[d]_-{R(p_{1,3}')_*} \ar[r] & \bZ_{1,3} \ar[d]^-{R(p_{1,3})_*} \\ \bX' \, \Rtimes_{\bY'} \, \bX' \ar@<-0.5ex>[r] & \bX \, \Rtimes_{\bY} \, \bX. }
\]
Applying Proposition \ref{prop:basechange} and the compatibility of the tensor product with inverse images, one obtains \eqref{it:prop-basechangeaction-1}. Statement \eqref{it:prop-basechangeaction-2} can be proved similarly.\end{proof}

\begin{remark}
In most examples relevant for us $\bX$, $\bY$ and $\bY'$ will be ordinary schemes, and we will have $\Tor_{\neq 0}^{(\pi \circ f')^{-1} \calO_Y}((f')^{-1}\calO_{Y'}, (\pi')^{-1} \calO_X)=0$ for the morphisms $\pi : Y' \to Y$, $f'=X' \to Y'$, $\pi' : X' \to X$, so that $X'$ is also an ordinary scheme. Moreover, $X'$ will be regular, and $f:X \to Y$ will be proper, so that we are in the situation of Proposition \ref{prop:action}\eqref{it:prop-action-2}.
\end{remark}

\subsection{Geometric actions}

Motivated by the constructions of \S\S \ref{ss:convolution}, \ref{ss:actionconvolution}, we introduce the following notion.

\begin{defin}

A (weak) \emph{geometric action} of a group $\Gamma$ on a dg-scheme 
$\bX$ over $\bY$ is a homomorphism
from $\Gamma$ to the group of isomorphism classes of invertible objects in the monoidal category
$\calK_{\bX,\bY}$.  

Under the assumptions of Proposition \ref{prop:action}\eqref{it:prop-action-2}, a weak geometric action is called \emph{finite} if 
its image is contained in $\calK^{\mathrm{coh}}_{X,Y}$.

\end{defin}

According to Proposition \ref{prop:basechangeaction}, a weak geometric action induces a usual weak action
of $\Gamma$ by auto-equivalences of $\calD \QCoh(\bX')$ for any $\bX'$ as in that proposition.

\begin{remark}
The above constructions certainly admit a variant when an algebraic group $G$ acts on $\bX,\bY$ and we work with categories of equivariant sheaves, or perhaps more generally for stacks. Our elementary approach of Section \ref{sec:dgschemes} is not adapted to these settings, however. (See Section \ref{sec:equivariant-version} below for some partial results in this direction.)
\end{remark}

Now we come back to the setting of Theorem \ref{thm:existenceaction-R}. Using the terminology introduced above, Theorem \ref{thm:existenceaction} gives statement \eqref{it:thm-actiongeneral-1} of the following theorem. Statement \eqref{it:thm-actiongeneral-2} follows from Proposition \ref{prop:basechangeaction}. (The compatibility for the direct image functor follows from the compatibility for the inverse image functor by adjunction.)

\begin{thm}
\label{thm:actiongeneral}

Let $R=\bbZ[\frac{1}{n_G}]$, where $n_G$ is defined in \S {\rm \ref{ss:statement-R}}.

\begin{enumerate}

\item
\label{it:thm-actiongeneral-1}
There exists a natural finite geometric action of the group $\bB_{\aff}$ on $\wfrakg_R$ (respectively $\wcalN_R$) over $\frakg^*_R$.
\item
\label{it:thm-actiongeneral-2}
For any morphism of dg-schemes $\bX \to \frakg^*_R$ whose underlying morphism of schemes is quasi-projective, there exist geometric actions of $\bB_{\aff}$ on $\calD \QCoh(\wfrakg_R \, \Rtimes_{\frakg_R^*} \, \bX)$ and $\calD \QCoh(\wcalN_R \, \Rtimes_{\frakg_R^*} \, \bX)$. Moreover, the direct and inverse image functors
\[
\xymatrix{
\calD \QCoh(\wfrakg_R) \ar@<0.5ex>[rr]^-{Lp^*} & & \calD \QCoh(\wfrakg_R \, \Rtimes_{\frakg^*_R} \, \bX) \ar@<0.5ex>[ll]^-{Rp_*}
}
\]
for the projection $p: \wfrakg_R \, \Rtimes_{\frakg^*_R} \, \bX \to \wfrakg_R$ commute with these actions, and similarly for $\wcalN_R$.

If the assumptions of Proposition {\rm \ref{prop:action}} are satisfied, then the action is finite.

\end{enumerate}

The same results hold if $R$ is replaced by an algebraically closed field $\bk$ satisfying the assumptions of Theorem {\rm \ref{thm:existenceaction}}.

\end{thm}

\begin{remark} \label{rk:dualnilpotent} Note that the existence of the action for $\wcalN$ cannot be obtained from that for $\wfrakg$ using Proposition \ref{prop:basechangeaction}. Indeed the fiber product $\wfrakg \times_{\frakg^*} \calN'$ is not reduced, hence not isomorphic to $\wcalN$. Here $\calN':=G \cdot (\frakg/\frakb)^* \subset \frakg^*$ is the ``dual nilpotent cone.''\end{remark}

\subsection{Examples} \label{ss:examples}

Now we get back to the notation of Sections \ref{sec:braidrelations} and \ref{sec:kernels}. For simplicity we assume that $\bk$ is of characteristic $0$. In particular, in this case there exists an isomorphism of $G$-modules $\frakg \cong \frakg^*$, which identifies the ``dual nilpotent cone'' $\calN'$ of Remark \ref{rk:dualnilpotent} with the usual nilpotent cone $\calN$.

Consider a nilpotent element $\chi \in \frakg^*$. Let $\mathbb{S}$ be the corresponding Slodowy slice (see \cite{Sl, Gi}; here we follow the notation of \cite[\S 1.2]{Gi}). We also consider the (scheme theoretic) fiber product \[ \widetilde{\calS}:= \wcalN \times_{\frakg^*} \mathbb{S}. \] By \cite[Proposition 2.1.2]{Gi}, this scheme is a smooth variety, of dimension $2 \dim(\calB) - \dim(G \cdot \chi)$ (i.e.~twice the dimension of the associated Springer fiber).

\begin{lem}

The dg-scheme \[ \wcalN \, \Rtimes_{\frakg^*} \, \mathbb{S} \] is concentrated in degree $0$, i.e.~is quasi-isomorphic to the variety $\widetilde{\calS}$.

\end{lem}

\begin{proof} This follows from a simple dimension-counting: one observes that the codimension of $\wcalN \times_{\frakg^*} \mathbb{S}$ in $\wcalN$ is $\dim(G \cdot \chi)$, which is exactly the codimension of $\mathbb{S}$ in $\frakg^*$. The result follows, using a Koszul complex argument. (See \S \ref{ss:endproof} for similar arguments.)
\end{proof}

Using this lemma and Theorem \ref{thm:actiongeneral}, one deduces the following.

\begin{cor} \label{cor:actionSlodowy}

There exists a natural finite geometric action of $\bB_{\aff}$ on $\widetilde{\calS}$ over $\mathbb{S}$. Moreover, the inverse and direct image functors \[ \xymatrix{ \calD^b \Coh(\wcalN) \ar@<0.5ex>[r] & \calD^b \Coh(\widetilde{\calS}) \ar@<0.5ex>[l] } \] are compatible with this action and the one given by Theorem {\rm \ref{thm:existenceaction}}.

\end{cor}

Assume now that $G$ is simple, and that $\chi$ is subregular. Then, according to Brieskorn and Slodowy,
\[
\mathbb{S} \, \cap \, \calN
\]
is a Kleinian singularity, and $\widetilde{\calS}$ is its minimal resolution (see \cite{Sl}). In this case, the braid group action of Corollary \ref{cor:actionSlodowy} is related to spherical twists and mirror symmetry, and has been extensively studied by several mathematicians (see \cite{ST, Bd, IUU, BT}). In particular, it is proved in \cite{BT} (extending results of \cite{ST}) that, if $G$ is simply laced, the restriction of the action on $\widetilde{\calS}$ to $\bB$ is faithful. It follows easily that the actions on $\calD^b \Coh(\wcalN)$ and $\calD^b \Coh(\wfrakg)$ are also faithful. Similarly, in \cite{IUU} it is proved that if $G$ is of type $\mbfA$, then the restriction of the action on $\widetilde{\calS}$ to $\bB_{\aff}^{\mathrm{Cox}}$ is faithful. Again, it follows that the actions on $\calD^b \Coh(\wcalN)$ and $\calD^b \Coh(\wfrakg)$ are also faithful.

\section{Equivariant version}
\label{sec:equivariant-version}

We were not able to extend the theory of dg-sheaves on dg-schemes to the equivariant setting, i.e.~to the situation where the dg-scheme is endowed with an action of algebraic group, and we consider dg-sheaves which are equivariant for this action. (In particular, our proof of Proposition \ref{prop:equivdgalg} does not extend to this setting.) In this section we provide direct arguments to prove an extension of Theorem \ref{thm:existenceaction-R-equivariant} in the spirit of Theorem \ref{thm:actiongeneral}. This statement is used in \cite{BM}.

This section is independent of Sections \ref{sec:kernels}, \ref{sec:dgschemes} and \ref{sec:convolution}. In particular, the assumptions on schemes in Sections \ref{sec:dgschemes} and \ref{sec:convolution} are not in order anymore. We use the general theory of unbounded derived categories of equivariant quasi-coherent sheaves over schemes as developped e.g.~in \cite[\S 1.5]{VV}.

\subsection{Statement}

In this section we let $R$ be either $\bbZ[\frac{1}{n_G}]$, where $n_G$ is defined in \S \ref{ss:statement-R}, or an algebraically closed field which satisfies the assumptions of Theorem \ref{thm:existenceaction}.

As in \S \ref{ss:notation}, we consider the varieties $\frakg^*_R$, $\wcalN_R$, $\wfrakg_R$ over $R$. Let $\calG$ be a group with a fixed homomorphism to $G_R \times_R (\Gm)_R$. We denote by $\calB \calC_{\calG}$, respectively $\calB \calC_{\calG}'$, the category of affine noetherian schemes $S$ endowed with a $\calG$-action and an equivariant morphism to $\frakg^*_R$ such that the natural morphism of dg-schemes
\[
S \times_{\frakg^*_R} \wfrakg_R \to S \, \Rtimes_{\frakg^*_R} \, \wfrakg_R, \quad \text{ respectively } \quad S \times_{\frakg^*_R} \wcalN_R \to S \, \Rtimes_{\frakg^*_R} \, \wcalN_R,
\]
is a quasi-isomorphism. Morphisms are assumed to be equivariant and over $\frakg^*_R$. For such an $S$, we set $\widetilde{S}:= S \times_{\frakg^*_R} \wfrakg_R$, respectively $\widetilde{S}':= S \times_{\frakg^*_R} \wcalN_R$.

\begin{thm}
\label{thm:action-equivariant}

For any $S$ in $\calB \calC_{\calG}$, respectively $\calB \calC_{\calG}'$, there exists an action of $\bB_{\aff}$ on the category $\calD^b \Coh^{\calG}(\widetilde{S})$, respectively $\calD^b \Coh^{\calG}(\widetilde{S}')$, such that for any $\calG$-equivariant morphism $S_1 \to S_2$, the action commutes with the direct and inverse image functors
\[
\xymatrix{
\calD^b \Coh^{\calG}(\wS_1) \ar@<0.5ex>[r] & \calD^b \Coh^{\calG}(\wS_2), \ar@<0.5ex>[l]
}
\]
respectively
\[
\xymatrix{
\calD^b \Coh^{\calG}(\wS_1') \ar@<0.5ex>[r] & \calD^b \Coh^{\calG}(\wS_2'). \ar@<0.5ex>[l]
}
\]

These actions are also compatible with the change of equivariance functors (for a group morphism $\calG' \to \calG$ over $G_R \times_R (\Gm)_R$).

\end{thm}

The proof of Theorem \ref{thm:action-equivariant} will be completed in \S \ref{ss:proof-equivariant}. We give the details only in the case of $\wfrakg_R$; the case of $\wcalN_R$ is similar.

For $S$ in $\calB \calC_{\calG}$, we write $S=\mathrm{Spec}(A_S)$, where $A_S$ is a $\calG$-equivariant $\rmS_R(\frakg_R)$-algebra. If $f: S \to \frakg^*_R$ is the structure morphism, we denote by $\widetilde{f} : \widetilde{S} \to \wfrakg_R$ the morphism obtained by base change.

\subsection{Quasi-isomorphisms of equivariant dg-schemes}

As explained a\-bove, it is not easy to develop a general theory of equivariant dg-schemes. However, part of the theory is easy to adapt. (See also \cite[Section 1]{MRHec} for such results.)

Let $\bX:=(X,\calA)$, where $X$ is a scheme and $\calA$ is a graded-commutative, non-positively graded, quasi-coherent sheaf of $\calO_X$-dg-algebras. We assume furthermore that an algebraic group $H$ acts on $X$, that $\calA$ is $H$-equivariant, and that the multiplication and differential are equivariant. We call such objects $H$-equivariant dg-schemes. We assume furthermore that every $H$-equivariant quasi-coherent sheaf on $X$ is a quotient of an $H$-equivariant quasi-coherent sheaf which is flat as an $\calO_X$-module. (See \cite[Remark 1.5.4]{VV} for comments on this assumption.)

We denote by $\calC \QCoh^H(\bX)$ the category of $H$-equivariant, quasi-coherent sheaves of $\calA$-dg-modules on $X$, and by $\calD \QCoh^H(\bX)$ the associated derived category. It follows easily, as in \cite[Theorem 1.3.3]{R2}, that there are enough objects in the category $\calC \QCoh^H(\bX)$ which are K-flat as $\calA$-dg-modules. Then one can adapt the proof of \cite[Theorem 10.12.5.1]{BL} to prove the following.

\begin{prop}
\label{prop:qis-equivalence-equivariant}

Let $(X,\calA)$ and $(X,\calA')$ be two $H$-equivariant dg-sche\-mes over the same ordinary $H$-scheme $X$ which satisfies the assumption above, and let $\phi : \calA \to \calA'$ be an $H$-equivariant quasi-isomorphism of dg-algebras. Then the extension and restriction functors induce equivalences of categories
\[
\calD \QCoh^H(X,\calA) \ \cong \ \calD \QCoh^H(X,\calA').
\]

\end{prop}

\subsection{Definition of the kernels}

From now on, for any $S$ in $\calB \calC_{\calG}$ we fix a $\calG$-equivariant graded-commutative $\mathrm{S}(\frakg_R)$-dg-algebra $D_S$ which is K-flat as an $\rmS_R(\frakg_R)$-dg-module, and a quasi-isomorphism of $\calG$-equivariant $\rmS_R(\frakg_R)$-dg-algebras $D_S \xrightarrow{\qis} A_S$. (For existence, see the arguments of \cite[Proof of Theorem 2.6.1]{CK}.)

Fix some $S$ as above, and denote by $f: S \to \frakg_R^*$ the associated morphism. By assumption, the derived tensor product
\[
\calO_{\wfrakg_R} \, \lotimes_{\mathrm{S}(\frakg_R)} \, A_S
\]
is concentrated in degree $0$. Hence the morphism of dg-algebras
\[
\calO_{\wfrakg_R} \otimes_{\mathrm{S}(\frakg_R)} D_S \to \calO_{\wfrakg_R} \otimes_{\mathrm{S}(\frakg_R)} A_S
\]
is a quasi-isomorphism. The morphism $\wf: \wS \to \wfrakg_R$ is affine, hence the functor $\wf_*$ is exact, and identifies the category $\Coh^{\calG}(\wS)$ with the category of $\calG$-equivariant coherent $\wf_* \calO_{\wS}$-modules on $\wfrakg_R$. Moreover, there is an isomorphism
\[
\wf_* \calO_{\wS} \ \cong \ \calO_{\wfrakg_R} \otimes_{\mathrm{S}(\frakg_R)} A_S.
\]
We deduce from these remarks, using Proposition \ref{prop:qis-equivalence-equivariant}, that there is a natural equivalence of categories
\begin{equation}
\label{eq:equivalence-wS}
\calD \QCoh^{\calG}(\wS) \ \cong \ \calD \QCoh^{\calG}(\wfrakg_R, \calO_{\wfrakg_R} \otimes_{\mathrm{S}(\frakg_R)} D_S).
\end{equation}
Similar remarks apply to the schemes $\wS \times_R \wfrakg_R$, $\wfrakg_R \times_R \wS$, $\wS \times_R \wS$.

Recall that, for any $t \in \scS$, we have defined the kernel $\calO_{Z_{t,R}}$ in the category $\Coh^{G_R \times (\Gm)_R}(\wfrakg_R \times \wfrakg_R)$. Consider the object 
\[
L(\widetilde{f} \times \Id_{\wfrakg_R})^* \calO_{Z_{t,R}} \ \text{ in } \calD^b \Coh^{\calG}(\widetilde{S} \times \wfrakg_R).
\]

\begin{lem}
\label{lem:definition-KSt}

For any $S$ in $\calB \calC_{\calG}$, there exists an object $\calK_t^S$ in $\calD^b \Coh^{\calG}(\wS \times_R \wS)$ and an isomorphism
\[
R(\Id_{\wS} \times \wf)_* \calK_t^S \ \cong \ L(\widetilde{f} \times \Id_{\wfrakg})^* \calO_{Z_{t,R}}
\]
in $\calD^b \Coh^{\calG}(\wS \times_R \wfrakg)$. Moreover, these objects can be chosen in such a way that: 

\begin{enumerate}
\item the objects $\calK^S_t$ are compatible with change of equivariance functors;
\item if $g: S_1 \to S_2$ is a morphism in the category $\calB \calC_{\calG}$, with morphism $\widetilde{g}: \wS_1 \to \wS_2$ obtained by change change, there is an isomorphism
\begin{equation}
\label{eq:isom-base-change-wS}
R(\Id_{\wS_1} \times \widetilde{g})_* \calK_t^{S_1} \ \cong \ L(\widetilde{g} \times \Id_{\wS_2})^* \calK^{S_2}_t
\end{equation}
in $\calD^b \Coh^{\calG}(\wS_1 \times_R \wS_2)$.
\end{enumerate}

\end{lem}

\begin{proof}
By the variant of equivalence \eqref{eq:equivalence-wS} for $\wS \times_R \wfrakg_R$, we have an equivalence of categories
\begin{multline*}
\calD \QCoh^{\calG}(\wS \times_R \wfrakg_R) \ \cong \\ \calD \QCoh^{\calG} \bigl( \wfrakg_R \times_R \wfrakg_R, \calO_{\wfrakg_R \times \wfrakg_R} \otimes_{\rmS(\frakg_R) \otimes_R \rmS(\frakg_R)} (D_S \otimes_R \rmS(\frakg_R)) \bigr).
\end{multline*}
Under this equivalence, the object $L(\widetilde{f} \times \Id_{\wfrakg})^* \calO_{Z_{t,R}}$ corresponds to
\[
\calO_{Z_{t,R}} \otimes_{\rmS(\frakg_R) \otimes_R \rmS(\frakg_R)} (D_S \otimes_R \rmS(\frakg_R)).
\]
To prove the first assertion of the lemma, it is enough to check that this dg-module can be endowed with the structure of an $\calO_{\wfrakg_R \times \wfrakg_R} \otimes_{\rmS(\frakg_R) \otimes_R \rmS(\frakg_R)} (D_S \otimes_R D_S)$-dg-module. However, by definition $Z_{t,R}$ is a subscheme of $\wfrakg_R \times_{\frakg^*_R} \wfrakg_R$. Hence it suffices to extend the action of $D_S$ to an action of $D_S \otimes_R D_S$ which factors through the multiplication map $D_S \otimes_R D_S \to D_S$.

Let us consider a morphism $S_1 \to S_2$ in $\calB \calC_{\calG}$, associated with an algebra morphism $A_2 \to A_1$. Here, for simplicity, we write $A_1$ for $A_{S_1}$ and $A_2$ for $A_{S_2}$. Similarly, we write $D_1$ and $D_2$ for $D_{S_1}$ and $D_{S_2}$. Then $A_1$ is a $D_1 \otimes_{\rmS(\frakg_R)} D_2$-dg-algebra. (Here, the product in $D_1 \otimes_{\rmS(\frakg_R)} D_2$ is defined using the usual sign rule.) Hence one can choose a $D_1 \otimes_{\rmS(\frakg_R)} D_2$-dg-algebra $\overline{D}_1$ which is graded-commutative and K-flat as a $D_1 \otimes_{\rmS(\frakg_R)} D_2$-dg-module (hence also automatically as a $\mathrm{S}(\frakg_R)$-dg-module), and a quasi-isomorphism $\overline{D}_1 \xrightarrow{\qis} A_1$. Then, we have a quasi-isomorphism $D_1 \to \overline{D}_1$, and a dg-algebra morphism $D_2 \to \overline{D_1}$. As above, we have equivalences of categories
\begin{multline*}
\calD \QCoh^{\calG} \bigl( \wfrakg_R \times_R \wfrakg_R, \calO_{\wfrakg_R \times \wfrakg_R} \otimes_{\rmS(\frakg_R) \otimes_R \rmS(\frakg_R)} (\overline{D}_1 \otimes_R \overline{D}_1) \bigr) \ \cong \\ 
\calD \QCoh^{\calG} \bigl( \wfrakg_R \times_R \wfrakg_R, \calO_{\wfrakg_R \times \wfrakg_R} \otimes_{\rmS(\frakg_R) \otimes_R \rmS(\frakg_R)} (D_1 \otimes_R D_1) \bigr) \\
\cong \ \calD \QCoh^{\calG}(\wS_1 \times_R \wS_1).
\end{multline*}
In the first category, the object $\calK^{S_1}_t$ defined above corresponds to the dg-module
\[
\calO_{Z_{t,R}} \otimes_{\rmS(\frakg_R) \otimes_R \rmS(\frakg_R)} (\overline{D}_1 \otimes_R \rmS(\frakg_R)),
\]
endowed with a $\overline{D}_1 \otimes_R \overline{D}_1$-action which factors through the multiplication map. Similarly, there is an equivalence
\begin{multline*}
\calD \QCoh^{\calG}(\wS_1 \times_R \wS_2) \ \cong \\ \calD \QCoh^{\calG} \bigl( \wfrakg_R \times_R \wfrakg_R, \calO_{\wfrakg_R \times \wfrakg_R} \otimes_{\rmS(\frakg_R) \otimes_R \rmS(\frakg_R)} (\overline{D}_1 \otimes_R D_2) \bigr).
\end{multline*}
Under this equivalence and the preceding one, the functor $R(\Id_{\wS} \times \widetilde{g})_*$ is simply a restriction of scalars functor induced by the dg-algebra morphism $D_2 \to \overline{D}_1$. One can describe similarly the functor $L(\widetilde{g} \times \Id_{\wS_2})^*$ as an extension of scalars functor, and then isomorphism \eqref{eq:isom-base-change-wS} is clear.

The change of equivariance can be treated similarly.\end{proof}

The object $\calK^S_t$ defined in Lemma \ref{lem:definition-KSt} will be the kernel for the action of $T_t$. However, as $\widetilde{S}$ is not a smooth scheme in general, it is not obvious that the associated convolution functor restricts to an endo-functor of $\calD^b \Coh^{\calG}(\wS)$. This will be proved in the following proposition. Here, we extend the notation of \S \ref{ss:convolution1} and set, for any $H$-scheme $X$ and any $\calF$ in $\calD \QCoh^{H}(X \times_R X)$,
\[
F^{\calF}_X: \left\{ 
\begin{array}{ccc}
\calD \QCoh^{H}(X) & \to & \calD \QCoh^H(X) \\
\calH & \mapsto & R(p_2)_* (\calF \, \lotimes_{X \times X} \, L(p_1)^* \calH)
\end{array} \right. .
\]

\begin{prop}
\label{prop:action-Ts-wS}

\begin{enumerate}
\item 
\label{it:prop-action-Ts-wS-1}
Let $t \in \scS$ and $S$ in $\calB \calC_{\calG}$. Let $f: S \to \wfrakg_R$ be the structure morphism, and $\wf:\wS \to \wfrakg_R$ the morphism obtained by base change. Then the following diagram commutes:
\[
\xymatrix@C=2cm{
\calD \QCoh^{\calG}(\wS) \ar[d]_-{F_{\wS}^{\calK^S_t}} \ar[r]^-{R(\wf)_*} & \calD \QCoh^{\calG}(\wfrakg_R) \ar[d]^-{F_{\wfrakg}^{\calO_{Z_{t,R}}}} \\
\calD \QCoh^{\calG}(\wS) \ar[r]^-{R(\wf)_*} & \calD \QCoh^{\calG}(\wfrakg_R). \\
}
\]
\item 
\label{it:prop-action-Ts-wS-2}
The functor $F_{\wS}^{\calK^S_t}$ stabilizes the subcategory $\calD^b \Coh^{\calG}(\wS)$.

\end{enumerate}

\end{prop}

\begin{proof}
\eqref{it:prop-action-Ts-wS-1} This follows easily from the projection formula. Indeed, for $\calH$ an object of $\calD \QCoh^{\calG}(\wS)$, we have isomorphisms
\begin{align*}
R(\wf)_* F_{\wS}^{\calK^S_t}(\calH) \ & \cong \ R(\wf)_* R(p_2^{\wS})_* (\calK^S_t \, \lotimes_{\wS \times \wS} \, L(p_1^{\wS})^* \calH) \\
& \cong \ R(q_2)_* R(\Id_{\wS} \times \wf)_* (\calK^S_t \, \lotimes_{\wS \times \wS} \, L(p_1^{\wS})^* \calH) \\
& \cong \ R(q_2)_* R(\Id_{\wS} \times \wf)_* (\calK^S_t \, \lotimes_{\wS \times \wS} \, L(\Id_{\wS} \times \wf)^* L(q_1)^* \calH).
\end{align*}
Here, $q_1: \wS \times_R \wfrakg_R \to \wS$ and $q_2 : \wS \times_R \wfrakg_R \to \wfrakg_R$ are the projections. Now, by the projection formula, we obtain
\[
R(\wf)_* F_{\wS}^{\calK^S_t}(\calH) \ \cong \ R(q_2)_* \bigl( (R(\Id_{\wS} \times \wf)_* \calK^S_t) \, \lotimes_{\wS \times \wfrakg_R} \, L(q_1)^* \calH \bigr).
\]
By the definition of $\calK^S_t$ (see Lemma \ref{lem:definition-KSt}), we have $R(\Id_{\wS} \times \wf)_* \calK^S_t \cong L(\wf \times \Id_{\wfrakg_R})^* \calO_{Z_{t,R}}$. Hence, using again the projection formula,
\begin{align*}
R(\wf)_* F_{\wS}^{\calK^S_t}(\calH) \ & \cong \ R(p_2^{\wfrakg_R})_* R(\wf \times \Id_{\wfrakg}) \bigl( (L(\wf \times \Id_{\wfrakg})^* \calO_{Z_{t,R}}) \lotimes_{\wS \times \wfrakg_R} L(q_1)^* \calH \bigr) \\
& \cong \ R(p_2^{\wfrakg_R})_* \bigl( \calO_{Z_{t,R}} \lotimes_{\wfrakg_R \times \wfrakg_R} R(\wf \times \Id_{\wfrakg})_* L(q_1)^* \calH \bigr).
\end{align*}
Now it is clear that $R(\wf \times \Id_{\wfrakg})_* L(q_1)^* \calH \cong L(p_1^{\wfrakg_R})^* R(\wf)_* \calH$. This concludes the proof.

\eqref{it:prop-action-Ts-wS-2} The fact that $F_{\wS}^{\calK^S_t}$ stabilizes $\calD^b \QCoh^{\calG}(\wS)$ follows from \eqref{it:prop-action-Ts-wS-1} and the similar claim for the functor $F_{\wfrakg}^{\calO_{Z_{t,R}}}$, which is clear since $\wfrakg_R$ is smooth. Then the claim follows from the fact that the projection $\mathrm{Supp}(\calK^S_t) \to \wS$ is proper.
\end{proof}

\subsection{Proof of Theorem \ref{thm:action-equivariant}}
\label{ss:proof-equivariant}

Now we can prove Theorem \ref{thm:action-equivariant}. In this subsection, we denote by the same symbol ``$\star$'' the convolution functors
\begin{align*}
\calD \QCoh^{\calG}(\wS \times_R \wS) \times \calD \QCoh^{\calG}(\wS \times_R \wS) \ & \to \ \calD \QCoh^{\calG}(\wS \times_R \wS), \\
\calD \QCoh^{\calG}(\wS \times_R \wfrakg_R) \times \calD \QCoh^{\calG}(\wS \times_R \wS) \ & \to \ \calD \QCoh^{\calG}(\wS \times_R \wfrakg_R), \\
\calD \QCoh^{\calG}(\wfrakg_R \times_R \wfrakg_R) \times \calD \QCoh^{\calG}(\wS \times_R \wfrakg_R) \ & \to \ \calD \QCoh^{\calG}(\wS \times_R \wfrakg_R)
\end{align*}
defined as in \S \ref{ss:convolution1}.

For any $t \in \scS$, the action of the generator $T_t \in \bB_{\aff}$ on the category $\calD^b \Coh^{\calG}(\wS)$ is defined as the functor
\[
F_{\wS}^{\calK^S_t} : \calD^b \Coh^{\calG}(\wS) \to \calD^b \Coh^{\calG}(\wS)
\]
(see Proposition \ref{prop:action-Ts-wS}\eqref{it:prop-action-Ts-wS-2}). Through the composition $\wS \to \wfrakg_R \to \calB_R$, we can consider $\wS$ as a scheme over $\calB_R$, hence we have line bundles $\calO_{\wS}(x)$ for any $x \in \bbX$. We define the action of $\theta_x$ as the functor
\[
F_{\wS}^{\calO_{\Delta \wS}(x)} : \calD^b \Coh^{\calG}(\wS) \to \calD^b \Coh^{\calG}(\wS).
\]
Note that this functor is just the twist by the line bundle $\calO_{\wS}(x)$. Hence the claims of Proposition \ref{prop:action-Ts-wS} are also true for the kernel $\calO_{\Delta \wS}(x)$.

To prove that these functors associated with the generators induce an action of $\bB_{\aff}$, one easily checks that it is enough to prove the following claims:

\begin{enumerate}
\item If $t=s_{\alpha}$, there are isomorphisms
\[
\calK^S_t \star \calK^S_t(-\rho,\rho-\alpha) \ \cong \ \calK^S_t(-\rho,\rho-\alpha) \star \calK^S_t \ \cong \ \calO_{\Delta \wS}
\]
in $\calD \QCoh^{\calG}(\wS \times_R \wS)$.
\item The kernels $\calK^S_t$ ($t \in \scS$) satisfy the finite braid relations in the monoidal category $(\calD \QCoh^{\calG}(\wS \times_R \wS), \star)$.
\end{enumerate}

We will only give the proof of the second claim; the first one can be obtained similarly. Our proof is copied from \cite[\S 4]{RAct}. To fix notations, we take $r,t \in \scS$ whose braid relation is $T_r T_t T_r = T_t T_r T_t$. The other cases are similar. We have to prove that there is an isomorphism
\[
\calK^S_r \star \calK^S_t \star \calK^S_r \ \cong \ \calK^S_t \star \calK^S_r \star \calK^S_t,
\]
or equivalently using claim (1) that there is an isomorphism
\begin{multline}
\label{eq:braid-relation-wS}
\calK^S_t(-\rho,\rho-\alpha) \star \calK^S_r(-\rho,\rho-\beta) \star \calK^S_t(-\rho,\rho-\alpha) \star \calK^S_r \star \calK^S_t \star \calK^S_r \\ \cong \ \calO_{\Delta \wS},
\end{multline}
where $t=s_{\alpha}$, $r=s_{\beta}$. As $\wf$ is an affine morphism, it is enough to prove that 
\[
R^i(\Id_{\wS} \times \wf)_* \bigl( \calK^S_t(-\rho,\rho-\alpha) \star \calK^S_r(-\rho,\rho-\beta) \star \calK^S_t(-\rho,\rho-\alpha) \star \calK^S_r \star \calK^S_t \star \calK^S_r \bigr) = 0
\]
if $i>0$ and that there is an isomorphism of $(\Id_{\wS} \times \wf)_* \calO_{\wS \times \wS}$-modules
\begin{multline*}
R^0(\Id_{\wS} \times \wf)_* \bigl( \calK^S_t(-\rho,\rho-\alpha) \star \calK^S_r(-\rho,\rho-\beta) \star \calK^S_t(-\rho,\rho-\alpha) \star \calK^S_r \star \calK^S_t \star \calK^S_r \bigr) \\ \cong \ (\Id_{\wS} \times \wf)_* \calO_{\Delta \wS}.
\end{multline*}

For any $\calM$ in $\calD \QCoh^{\calG}(\wS \times_R \wS)$ and $\calM'$ in $\calD \QCoh^{\calG}(\wfrakg_R \times_R \wfrakg_R)$ we have isomorphisms
\[
\calO_{\Gamma_{\wf}} \star \calM \ \cong \ R(\Id_{\wS} \times \wf)_* \calM, \quad \calM' \star \calO_{\Gamma_{\wf}} \ \cong \ L(\wf \times \Id_{\wfrakg})^* \calM',
\]
where $\Gamma_{\wf}$ is the graph of $\wf$. (See \cite[Lemma 1.2.3]{RAct} for the first isomorphism, and the proof of \cite[Corollary 4.3]{RAct} for the second one.) In particular, using Lemma \ref{lem:definition-KSt}, for $u=s_{\gamma}$ we obtain isomorphisms
\begin{align*}
\calO_{\Gamma_{\wf}} \star \calK^S_u \ & \cong \ \calO_{Z_{u,R}} \star \calO_{\Gamma_{\wf}}, \\
\calO_{\Gamma_{\wf}} \star \calK^S_u(-\rho,\rho-\gamma) \ & \cong \ \calO_{Z_{u,R}}(-\rho,\rho-\gamma) \star \calO_{\Gamma_{\wf}}.
\end{align*}
Hence, convolving the left hand side of \eqref{eq:braid-relation-wS} on the left with $\calO_{\Gamma_{\wf}}$ we obtain
\begin{multline*}
\bigl( \calO_{Z_{t,R}}(-\rho,\rho-\alpha) \star \calO_{Z_{r,R}}(-\rho,\rho-\beta) \star \calO_{Z_{t,R}}(-\rho,\rho-\alpha) \star \calO_{Z_{r,R}} \star \calO_{Z_{t,R}} \star \calO_{Z_{r,R}} \bigr) \\ \star \calO_{\Gamma_{\wf}}.
\end{multline*}
Now we have proved in Corollary \ref{cor:braidrelations} that there is an isomorphism
\begin{multline*}
\calO_{Z_{t,R}}(-\rho,\rho-\alpha) \star \calO_{Z_{r,R}}(-\rho,\rho-\beta) \star \calO_{Z_{t,R}}(-\rho,\rho-\alpha) \star \calO_{Z_{r,R}} \star \calO_{Z_{t,R}} \star \calO_{Z_{r,R}} \\ \cong \calO_{\Delta \wfrakg_R}.
\end{multline*}
in $\calD^b \Coh^{G_R \times (\Gm)_R}(\wfrakg_R \times_R \wfrakg_R)$. The result follows.

Compatibility with direct images can be proved as in Proposition \ref{prop:action-Ts-wS}\eqref{it:prop-action-Ts-wS-1}, using equation \eqref{eq:isom-base-change-wS}. Then, compatibility with inverse images follows by adjunction. Finally, compatibility with change of scalars functors is clear by construction and Lemma \ref{lem:definition-KSt}.

\end{document}